\documentclass[11pt]{book}
\usepackage{amsfonts,phd}

\newcounter{realpagenumber}
\begin{document}

\pagenumbering{roman}
\pagestyle{empty}

\begin{flushleft} \hfill \end{flushleft}
\vspace{14 mm}
\begin{center}
{\Huge {\bf Stability and Invariants of \vspace{2.5 mm} \\
            Hilsum-Skandalis Maps}} \vspace{18 mm} \\
{\LARGE {\em Stabiliteit en \vspace{1.5 mm} \\
             Invarianten van Hilsum-Skandalis-Afbeeldingen}}
             \vspace{1 mm} \\
{\sl (met een samenvatting in het Nederlands)}
\vfill

{\Large {\em Proefschrift}} \vspace{4 mm} \\
{\large {\em ter verkrijging van de graad van doctor aan
de Universiteit Utrecht op gezag van de Rector Magnificus,
Prof. dr. J. A. van Ginkel ingevolge het besluit van het
College van Decanen in het openbaar te verdedigen op
maandag 10 juni 1996 des middags te 16.15 uur
\vspace{20 mm} \\
door \vspace{4 mm} \\
{\LARGE {\sc Janez Mr\v{c}un}} \vspace{2 mm} \\
geboren op 13 april 1966, te Ljubljana.}}
\end{center}

\newpage

\noindent
{\sl promotoren:} \vspace{0.8 mm} \\
{\large Prof. dr. I. Moerdijk, Faculteit der Wiskunde en Informatica, \\
\indent \hspace{33.05 mm} Universiteit Utrecht}
\vspace{0.8 mm} \\
{\large Prof. dr. D. Siersma, \hspace{-0.8 mm}
Faculteit der Wiskunde en Informatica, \\
\indent \hspace{33.05 mm} Universiteit Utrecht}

\vfill

\noindent {\small
\hrule \vspace{1 mm} \noindent
Mathematics Subject Classifications (1991): 58H05, 22A22, 57R30, 57S25
\vspace{1 mm}
\hrule
\vspace{6 mm}
\noindent
CIP-DATA KONINKLIJKE BIBLIOTHEEK, DEN HAAG \vspace{2 mm} \\
Mr\v{c}un, Janez \vspace{2 mm} \\
Stability and invariants of Hilsum-Skandalis maps / Janez
Mr\v{c}un. - \\
Utrecht : Universiteit Utrecht, Faculteit Wiskunde en Informatica \\
Thesis Universiteit Utrecht. - With ref. - With summary in Dutch. \\
ISBN 90-393-1314-8 \\
Subject headings: foliations / topological groupoids / Reeb stability.

\vspace{6 mm}
\noindent
\copyright\ 1996 Janez Mr\v{c}un \\
Printed and bound by Drukkerij Elinkwijk bv, Utrecht}


\chapter*{Preface \markboth{PREFACE}{PREFACE}}
\pagestyle{myheadings}
\startXchapterskip

For the last four years, my work and life in Holland have been
influenced by several people; here I can
mention only few of them. First, I am indebted to my supervisor
professor Ieke Moerdijk: he guided my work carefully,
by pointing out interesting problems during
many discussions, and by broadening my mathematical knowledge
in general. I would like to express my gratitude to professor
P. Molino, professor D. Siersma and in particular to professor W. T.
van Est for their interest in my research.
Also, I would like to thank the Mathematical Department of Utrecht
University for the nice and peaceful working environment, and
Dennis Hesseling for correcting my Dutch in the {\em Samen\-vatting}.
I am grateful to the Mathematical Institute of Ljubljana University
for their support, especially to professor A. Suhadolc and professor
P. \v{S}emrl: without it I would not have been able to start
my studies in Utrecht.

Finally, I would like to thank my parents, my brother and my sister
for their care. A lot I owe to Andreja and her encouragements
all along. And it would be hard to imagine my life for the past few
years without a friend as our dog Tim is.


\chapter*{Contents \markboth{CONTENTS}{CONTENTS}}
\startXchapterskip

\small

\PHDstarttoc

\normalsize


\chapter*{Introduction \markboth{INTRODUCTION}{INTRODUCTION}}
\addcontentsline{toc}{chapter}{Introduction}
\setcounter{realpagenumber}{\value{page}}
\pagenumbering{arabic}
\setcounter{page}{\value{realpagenumber}}
\startXchapterskip

In order to understand foliations, it is of great importance
to investigate the structures obtained by transporting a
foliation along a map that need not be transversal to the leaves.
If not before, this became evident after the Haefliger's
proof of non-existence of analytic foliations of codimension one
on $S^{3}$. A. Haefliger was the first who
considered a foliation on $M$ not solely as a partition of $M$, but
as a cohomology class of cocycles on $M$ with values in the
pseudogroup of diffeomorphisms of $\RRR^{q}$ -- an instance of
what is nowadays known as a Haefliger
structure on a topological space $X$ with values in an \'{e}tale
groupoid $\GG$ \cite{Hae}. Such structures are not only fully
functorial in $X$, but can also be seen as generalized
foliations on $X$. Indeed, a Haefliger $\GG$-structure on $X$
induces a partition of $X$ into leaves equipped with their
holonomy groups, and this partition is locally given by the
fibers of a continuous map.

This is not the only direction from which topological
groupoids enter foliation theory. The holonomy groupoid of
a foliation is a good choice to represent the transversal structure
of the foliation \cite{Con82,Win},
in contrast with the space of leaves which is
often useless. By reducing this groupoid to a complete transversal
one obtains a Morita equivalent \'{e}tale groupoid.
Moreover, a foliation on $M$ can be canonically represented
by a Haefliger structure on $M$ with values in the reduced holonomy
groupoid of the foliation -- the structure which plays the role
of the quotient projection from $M$ to the transversal structure
of the foliation.
Indeed, a Haefliger $\GG$-structure on $X$ should be seen as
a map from $X$ to $\GG$, as a special case of a Hilsum-Skandalis
map between topological groupoids \cite{Hae84,HilSka,Pra}.
These maps were introduced in \cite{HilSka} as the maps between
the transversal structures of foliations. A Hilsum-Skandalis map between
$\HH$ and $\GG$ is an isomorphism class of principal
$\GG$-$\HH$-bibundles. These maps can be composed: they form a
category in which two topological groupoids are isomorphic
if and only if they are Morita equivalent.
Furthermore, the category of topological spaces is a full subcategory
of the category of Hilsum-Skandalis maps, and the Hilsum-Skandalis
maps from a topological space $X$ to an \'{e}tale groupoid $\GG$ are
precisely the Haefliger $\GG$-structures on $X$.
The Hilsum-Skandalis maps arose independently in topos theory:
in \cite{Moe4} I. Moerdijk showed that the Hilsum-Skandalis
maps between suitable groupoids $\HH$ and $\GG$ are in a natural
bijective correspondence with the geometric morphisms between the
classifying toposes of $\HH$ and $\GG$ (see also \cite{Bun}).

The first goal of this thesis is to show that any Hilsum-Skandalis
map from $\HH$ to $\GG$ can be seen as a generalized foliation on
the topological groupoid $\HH$: it induces a partition of the space
of objects $\HH_{0}$
of $\HH$ into $\HH$-invariant subsets -- the leaves.
A leaf comes equipped with its leaf topology, a right $\HH$-action
and its holonomy group. We show that if $\GG$ is \'{e}tale, then the
holonomy group of a leaf $L$ is the image
of a homomorphism (the holonomy homomorphism)
from the fundamental group $\pi_{1}(\HH(L))$
of the groupoid $\HH(L)$
associated to the $\HH$-action on $L$.

In Chapter \ref{chapSta} we study the stability of the
leaves of Hilsum-Skandalis maps with compact orbit spaces.
In its simplest form, the Reeb stability theorem \cite{Reeb}
says that a compact leaf of a foliation on a manifold with finite
holonomy group has a neighbourhood of compact leaves. In particular,
any compact leaf with finite fundamental group has a neighbourhood of
compact leaves. The Haefliger-Reeb-Ehresmann stability theorem
\cite{Hae} extends this result to
Haefliger structures on topological spaces.
On the other hand, the Reeb stability theorem was generalized by
W. P. Thurston in \cite{Thu}. We generalize the
Haefliger-Reeb-Ehresmann stability theorem as well as
the Reeb-Thurston stability theorem to Hilsum-Skandalis maps.
The proofs are eventually reduced to the simple case
of a continuous map between topological spaces.
To generalize the Reeb-Thurston stability theorem, we need to
introduce the linear holonomy homomorphism $d\cH_{L}$
of a leaf $L$ of a Hilsum-Skandalis map $E$ from $\HH$ to an \'{e}tale
$\eCe$-groupoid $\GG$ modeled on a Banach space.
We define $d\cH_{L}$ as a representation
of the fundamental group $\pi_{1}(\HH(L))$ 
on the tangent space of $\GG_{0}$ at a suitable point.
If the map $E$ satisfies some
mild conditions, we prove the following: Let $L$ be a leaf of $E$
such that the orbit space of $L$ is compact. Assume that
\begin{enumerate}
\item [(i)]   $\Ker d\cH_{L}$ is finitely generated,
\item [(ii)]  $\Hom(\Ker d\cH_{L},\RRR)=\{0\}$, and
\item [(iii)] the image of $d\cH_{L}$ is finite.
\end{enumerate}
Then the holonomy group of $L$ is finite, and there exists an
$\HH$-invariant neighbourhood of $L$ in $\HH_{0}$
which is a union of leaves of $E$ with compact orbit spaces.
This generalizes the Reeb-Thurston stability theorem. To facilitate
the proof, we introduce the notion of a differential category.
We show that these categories naturally arise from the affine manifolds.

In Chapter \ref{chapEquFol} we illustrate these results in some
special cases. We show that an equivariant foliation $\cF$ on a manifold
$M$ equipped with an action of a group $G$ can be represented
by a Hilsum-Skandalis map from the groupoid $G(M)$ associated to the
action, to the groupoid of germs of diffeomorphisms of $\RRR^{q}$.
We establish the relation between
the holonomy groups of the leaves of this Hilsum-Skandalis
map and the equivariant holonomy groups of the leaves of $\cF$.
As a consequence of the stability theorem for Hilsum-Skandalis maps
and under some mild conditions, we prove that a leaf of $\cF$ with
finite equivariant holonomy group and compact orbit space has a
$G$-saturated neighbourhood of leaves of $\cF$ with compact orbit
spaces. If in addition the foliation $\cF$ is of codimension one,
we obtain the following generalization of the Reeb-Thurston stability
theorem: Assume that $L$ is a leaf of $\cF$ with finitely generated
equivariant fundamental group such that
\begin{enumerate}
\item [(i)]   the foliation is equivariantly transversely orientable
              around $L$,
\item [(ii)]  the orbit space of $L$ is compact, and
\item [(iii)] the first equivariant cohomology group of $L$ with
              coefficients in $\RRR$ is trivial.
\end{enumerate}
Then the equivariant holonomy group of $L$ is trivial, and there
exists a $G$-saturated neighbourhood of $L$ which is a union of
leaves of $\cF$ with compact orbit spaces.

In Chapter \ref{chapInvHilSkaMap} we introduce some
algebraic invariants of Hilsum-Skandalis maps.
First, we use the category of Hilsum-Skandalis maps
to define higher homotopy groups of topological
groupoids. Since these groups are functorial with respect to the
(pointed) Hilsum-Skandalis maps, they are invariant under Morita
equivalence. For a suitable groupoid $\GG$, these homotopy
groups coincide with the higher homotopy groups of the
classifying space and of the classifying topos of $\GG$ \cite{Moe1}.
In particular, the first homotopy group of a suitable
effective groupoid $\GG$ is exactly the fundamental group
of $\GG$ as described by W. T. van Est \cite{Est}.

Next, we introduce singular homology groups of
topological groupoids. These groups are much simpler than
the homology groups of the classifying spaces associated to
topological groupoids. The singular homology groups
are functorial with respect to the
Hilsum-Skandalis maps and hence invariant under Morita equivalence.
In particular, we can define the $n$-th singular homology
group $H_{n}((M,\cF);\ZZZ)$ of the transverse structure of a
foliation $\cF$ on a manifold $M$, and obtain a homomorphism
$$ H_{n}(M;\ZZZ)\lra H_{n}((M,\cF);\ZZZ)\;.$$
For an \'{e}tale groupoid $\GG$ we give a particularly simple description
of the singular homology groups of $\GG$ in terms of singular
simplexes in the space of objects and in the space of morphisms of
$\GG$. Moreover, we prove that the effect-functor induces
isomorphisms between the singular homology groups
of $\GG$ and of its effect $\Eff(\GG)$.

Finally, we consider the algebra $\eCrc(\GG)$ of $\eCrc$-functions on
an object-separated finite-dimensional \'{e}tale $\eCr$-groupoid $\GG$,
introduced by A. Connes \cite{Con82}. To a Hilsum-Skandalis
map $E$ between such groupoids $\HH$ and $\GG$ we associate a
$\eCrc(\GG)$-$\eCrc(\HH)$-bimodule $\eCrc(E)$, and prove that
$\eCrc$ becomes a functor from the category of Hilsum-Skandalis maps
between separated finite-dimensional \'{e}tale $\eCr$-groupoids to the
category of locally unital bimodules over algebras
with local units. In particular, this yields that
$\eCr$-Morita equivalent separated finite-dimensional \'{e}tale
$\eCr$-groupoids have Morita equivalent algebras -- a more
general version of this result was announced
in \cite{BryNis}.

\pagestyle{headings}
\chapter{Background}
\label{chapBac}
\startchapterskip

To start with, we recall definitions and some well-known
facts about foliations \cite{Cam,EhrRee,Reeb,Thu},
topological groupoids \cite{Bro,Con82,Est,Moe1,Mol,Win}
and Haefliger structures \cite{Bot,Hae}.

\section{Foliations}  \label{secFol}

Let $M$ be a $\eCe$-manifold of
dimension $n$ (without boundary).
A {\em foliation} $\cF$ on $M$ of codimension $q$
($0\leq q\leq n$) is a maximal $\eCe$-atlas
$(\varphi_{i}:U_{i}\ra\RRR^{n})$ on $M$ for
which the change of coordinates diffeomorphisms
$\varphi_{ij}=\varphi_{i}\com\varphi_{j}^{-1}|_{\varphi_{j}
 (U_{i}\cap\, U_{j})}$
are locally of the form
$$ \varphi_{ij}(x,y)=(\varphi^{(1)}_{ij}(x,y),
   \varphi^{(2)}_{ij}(y)) $$
with respect to the decomposition $\RRR^{n}=\RRR^{n-q}
\times\RRR^{q}$. The connected components of
$\varphi_{i}^{-1}(\RRR^{n-q}\times\{y\})$, $y\in\RRR^{q}$,
are called the {\em plaques} of $\cF$ in $U_{i}$. They
decompose $U_{i}$ into $(n-q)$-dimensional submanifolds, and a
change of coordinates preserves this decomposition. Therefore
the plaques globally amalgamate into connected
$\eCe$-manifolds of dimension $n-q$, called the
{\em leaves} of $\cF$, which decompose $M$. Each leaf is
injectively immersed into $M$, but it may not be embedded.

The {\em space of leaves} $M/\cF$ of $\cF$ is
obtained as the quotient of $M$ by identifying two points
if they lie on the same leaf.

Let $L$ be a leaf of $\cF$ and let $\sigma$ be a
path in $L$. Let $T$ and $T'$ be transversal sections
of $\cF$ (i.e. $\eCe$-immersions $\RRR^{q}\ra M$ which are
transversal to the leaves of $\cF$) such that
$T(0)=\sigma(0)$ and $T'(0)=\sigma(1)$. 
The foliation on a small neighbourhood of 
$\sigma([0,1])$ provides a $\eCe$-diffeomorphism between small
neighbourhoods of the origin in $\RRR^{q}$ which fixes the
origin, and the germ at $0$ of this diffeomorphism
$$ \Hol_{T',T}(\sigma)\in\Dif{q} $$
is called the {\em holonomy} of $\sigma$ with respect to the
transversal sections $T$ and $T'$.
If two paths in $L$ are homotopic in $L$ with fixed end-points,
they have the same holonomy, hence we can regard $\Hol_{T',T}$
as a function defined on the homotopy classes of the paths in
$L$ from $T(0)$ to $T'(0)$.
If $\sigma'$ is another path in $L$ with $\sigma'(0)=
\sigma(1)$ and $T''$ a transversal section of $\cF$
with $T''(0)=\sigma'(1)$, then
$$ \Hol_{T'',T'}(\sigma')\com\Hol_{T',T}(\sigma)  
   =\Hol_{T'',T}(\sigma'\sigma)\;.$$
In particular, we have a homomorphism
$$ \Hol_{T}=\Hol_{T,T}:\pi_{1}(L,T(0))\lra\Dif{q}\;,$$
called the {\em holonomy homomorphism} of $L$ (with respect
to $T$). The image of $\Hol_{T}$ is called the {\em holonomy group}
of $L$. The holonomy homomorphism and the holonomy group
are determined uniquely with the leaf $L$ up to a conjugation
in the group $\Dif{q}$.
If holonomy group of $L$ is a subgroup of $\Difp{q}$,
i.e. if all the germs in $\Hol_{T}(\pi_{1}(L,T(0)))$ preserve
the orientation of $\RRR^{q}$, we say that $\cF$ is
{\em transversely orientable around} $L$.

By composing the holonomy homomorphism of $L$ with the homomorphism
$$ d:\Dif{q}\lra GL_{q}(\RRR)\;,$$
which sends a germ of a diffeomorphism to its differential at the
origin, one gets the {\em linear holonomy homomorphism}
$$ d\Hol_{T}:\pi_{1}(L,T(0))\lra GL_{q}(\RRR) $$
of $L$ with respect to $T$. The image of $d\Hol_{T}$ is called
the {\em linear holonomy group} of $L$, and is determined
up to a conjugation by the leaf $L$.

The simplest form of the (local) Reeb stability theorem,
proved by G. Reeb in 1957 \cite{Reeb}, is the following:

\begin{theo}  \label{Itheo1}
Let $\cF$ be a foliation on a finite-dimensional Hausdorff
$\eCe$-manifold $M$.
If $L$ is a compact leaf of $\cF$ with finite holonomy
group, then $L$ has an arbitrary small neighbourhood of
compact leaves.

In particular, if $L$ is a compact leaf of $\cF$ with finite
fundamental group, then $L$ has an arbitrary small
neighbourhood of compact leaves.
\end{theo}

The second part of this theorem was generalized by W. P.
Thurston \cite{Thu} as follows:

\begin{theo}  \label{Itheo2}
Let $\cF$ be a foliation on a finite-dimensional Hausdorff
$\eCe$-manifold $M$. If $L$ is a compact leaf of $\cF$ with
trivial linear holonomy group, and if $H^{1}(L;\RRR)=0$,
then $L$ has an arbitrary small 
neighbourhood of compact leaves.

In particular, if $\cF$ is of codimension one and $L$ is a compact
leaf of $\cF$ such that $\cF$ is transversely orientable around
$L$ and $H^{1}(L;\RRR)=0$, then $L$ has an arbitrary small
neighbourhood of compact leaves.
\end{theo}

A foliation of codimension $q$ on a $\eCe$-manifold $M$ of
dimension $n$ can be 
represented by a family of $\eCe$-submersions
$(s_{i}:V_{i}\ra\RRR^{q})$, where $(V_{i})$ is an open cover of $M$,
which satisfy the following condition: For any $x\in V_{i}\cap V_{j}$
there exists an open neighbourhood $W\subset V_{i}\cap V_{j}$ of $x$
and a $\eCe$-diffeomorphism $s_{ij}^{W}:s_{j}(W)\ra s_{i}(W)$ such that
$$ s_{i}|_{W}=s_{ij}^{W}\com s_{j}|_{W}\;.$$
The leaves of the foliation are now obtained as the amalgamations
of the connected components of the fibers of the submersions.

With a fixed $W$, the diffeomorphism $s_{ij}^{W}$ is clearly uniquely
determined. In particular, for any $x\in V_{i}\cap V_{j}$ the germ
of $s_{ij}^{W}$ at $s_{j}(x)$ does not depend on the choice of $W$.
We will denote this germ by $c_{ij}(x)$.

Denote by $\GmqCe$ the set of all germs of
(locally defined)
$\eCe$-diffeomorphisms of $\RRR^{q}$. If $\gm=\germ_{x}f$ for a
diffeomorphism $f:U\ra V$, where $U$ and $V$ are open subsets of
$\RRR^{q}$ and $x\in U$, we define the domain $\dom\gm$ of $\gm$ to be
$x$ and the codomain $\cod\gm$ to be $f(x)$. This gives the domain
and the codomain functions
$$ \dom,\cod:\GmqCe\lra\RRR^{q}\;.$$
For any $x\in\RRR^{q}$ there is a germ $\uni x=1_{x}$ of the
identity diffeomorphism at $x$. This gives the function
$$ \uni:\RRR^{q}\lra\GmqCe\;.$$
If $\gm',\gm\in\GmqCe$ are such that $\dom\gm'=\cod\gm$, we define
the composition $\cmp(\gm',\gm)=\gm'\com \gm$ to be the germ at $\dom\gm$
of the composition of diffeomorphisms representing $\gm'$ and $\gm$.
This gives a function
$$ \cmp:\GmqCe\times_{\RRR^{q}}\GmqCe\lra
   \GmqCe\;.$$
Finally, for any $\gm\in\GmqCe$ there is the inverse 
$\inv\gm=\gm^{-1}$, given as the germ at $\cod\gm$ of the inverse of a
diffeomorphism representing $\gm$. Hence we have a function
$$ \inv:\GmqCe\lra\GmqCe\;.$$
This gives $\GmqCe$ a structure of a category in which
all the morphisms are invertible.

There is a natural structure of a (non-Hausdorff) 
$\eCe$-manifold of dimension $q$ on $\GmqCe$,
which can be characterized as follows:
For any $\eCe$-diffeomorphism $f:U\ra V$ between open subsets
of $\RRR^{q}$, the map $\bar{f}:U\ra\GmqCe$, given by
$$ \bar{f}(x)=\germ_{x}f\;,$$
is a $\eCe$-diffeomorphism onto its (open) image.
With this structure, the maps
$\dom$, $\cod$, $\uni$, $\cmp$ and $\inv$ becomes local
$\eCe$-diffeomorphisms.

Now the germ $c_{ij}(x)$ defined above is an element of
$\GmqCe$, and this gives submersions
$$ c_{ij}:V_{i}\cap V_{j}\lra\GmqCe\;.$$
Moreover,
\begin{enumerate}
\item [(i)]   $c_{ii}(U_{i})\subset\uni(\RRR^{q})$,
\item [(ii)]  $\dom\com c_{ij}=\dom\com c_{jj}|_{V_{i}\cap\, V_{j}}$,
              $\cod\com c_{ij}=\cod\com c_{ii}|_{V_{i}\cap\, V_{j}}$,
              and
\item [(iii)] $\cmp\com (c_{ij},c_{jk})|_{V_{i}\cap\,
              V_{j}\cap\, V_{k}}=
              c_{ik}|_{V_{i}\cap\, V_{j}\cap\, V_{k}}$.
\end{enumerate}
Note that the submersions $(s_{i})$ can be recovered from the family
$(c_{ij})$ as $s_{i}=\dom\com c_{ii}=\cod\com c_{ii}$.
A family of $\eCe$-maps
$(c_{ij}:V_{i}\cap V_{j}\lra\GmqCe)$,
which satisfies the conditions (i), (ii) and (iii) above, is called a
$\GmqCe$-{\em cocycle of class} $\eCe$ on the open cover
$(V_{i})$ of $M$. A $\GmqCe$-cocycle of class $\eCe$
on $(V_{i})$ represents a foliation on $M$ precisely if
the maps of the cocycle are submersions.

\begin{ex}  \label{Iex3}  \rm
(1) Let $M$ and $N$ be finite-dimensional
$\eCe$-manifolds and $s:M\ra N$ a
submersion. Then $s$ defines a foliation $\cF$ on $M$. The
codimension of $\cF$ is the dimension of $N$.
The leaves are the connected components of the fibers of $s$,
they are all embedded, and they all have trivial holonomy groups.

(2) Let $\cF$ be a foliation on a $\eCe$-manifold $M$
of dimension $n$.
A diffeomorphism of $M$ preserves $\cF$ if it maps leaves into
leaves. Let $G$ be a properly discontinuous group of 
$\eCe$-diffeomorphisms of $M$ which preserve $\cF$.
Then $M/G$ is a $\eCe$-manifold of dimension $n$
such that the quotient projection
$$ p:M\lra M/G $$
is a covering projection of class $\eCe$. Moreover, there
is the induced foliation $\cF/G$ on $M/G$ such that $p$ maps
leaves of $\cF$ to leaves of $\cF/G$.

(3) As a special case of (1) and (2), take
$M=\RRR^{2}$, $N=\RRR$, $s=\pr_{2}:\RRR^{2}\ra\RRR$, and let
$G$ be the infinite cyclic group generated by the diffeomorphism
$f$ of $\RRR^{2}$ given by
$$ f(x,y)=(x+1,-y)\;.$$
Clearly $f$ (and hence also $f^{k}$) preserves the foliation $\cF$
on $M$ given by $s$, and $G$ is a properly discontinuous group.
The space $M/G$ is the open M\"{o}bius strip and all the leaves
of $\cF/G$ are diffeomorphic to $S^{1}$. All the leaves of $\cF/G$
have trivial holonomy groups with the
exception of the central one, i.e.
$s^{-1}(0)/G$, which has the holonomy group $\ZZZ/2\ZZZ$.

(4) Take $M=\RRR^{2}$, let $a\in (0,1)\setminus\QQQ$
and define a submersion
$s:\RRR^{2}\ra\RRR$ by
$$ s(x,y)=x-ay\;.$$
Let $G$ be the group
generated by diffeomorphisms $f$ and $f'$ of $\RRR^{2}$
given by $f(x,y)=(x+1,y)$ and $f'(x,y)=(x,y+1)$. The diffeomorphisms
$f$ and $f'$ preserve the foliation $\cF$ given by $s$, the group
$G$ is properly discontinuous, and $M/G\cong T^{2}$. Any leaf
of $\cF/G$ is diffeomorphic to $\RRR$, and is dense in $M/G$
with trivial holonomy group. The space of leaves
of $\cF/G$ has the trivial topology. The foliation
$\cF/G$ is called the Kronecker foliation of the torus. 

(5) Let $M=\RRR^{2}\setminus\{0\}$ and let $\cF$ be given by
the submersion $s=\pr_{2}:M\ra\RRR$.
Let $G$ be the properly discontinuous group generated by the
diffeomorphism $f$ of $M$ given by
$$ f(x,y)=(2x,2y)\;.$$
The diffeomorphism $f$ preserves $\cF$.
The induced foliation $\cF/G$ on the orbit space 
$M/G\cong T^{2}$ has two compact leaves, both diffeomorphic to
$S^{1}$ and with holonomy group $\ZZZ$.
Any other leaf of $\cF/G$ is diffeomorphic to $\RRR$, has trivial
holonomy group, and has both the compact leaves in its closure. 
\end{ex}

\section{Topological Groupoids}  \label{secTopGro}

A {\em groupoid} $\GG$ is a small category in which any morphism is
invertible.
We denote the set of objects of $\GG$ by $\GG_{0}$, and 
the set of all morphisms of $\GG$ by $\GG_{1}$. Since $\GG$ is a
category, we have the domain and the codomain map
$$ \dom,\cod:\GG_{1}\lra\GG_{0}\;,$$
the unit map
$$ \uni:\GG_{0}\lra\GG_{1} $$
and the composition map
$$ \cmp:\GG_{1}\times_{\GG_{0}}\GG_{1}\lra\GG_{1}\;.$$
We write $1_{a}=\uni(a)$ and $g'\com g=\cmp(g',g)$. In addition, we
have the inverse map
$$ \inv:\GG_{1}\lra\GG_{1} $$
which maps a morphism $g$ to its inverse $g^{-1}$. We write
$\GG(a,a')$ for the set of morphisms of $\GG$ with domain $a$ and
codomain $a'$. The group $\GG(a,a)$ is called the {\em vertex group}
of $\GG$ at $a$.

A {\em topological groupoid} \cite{Bro,Hae,Moe4}
is a groupoid $\GG$ such that $\GG_{0}$ and $\GG_{1}$
are topological spaces, the structure maps of $\GG$ (i.e. the maps
$\dom$, $\cod$, $\uni$, $\cmp$ and $\inv$) are continuous, and the
domain map of $\GG$ is open.

In a topological groupoid $\GG$,
the codomain and the compositions map are also open, and
the unit map is an embedding. In particular, we will often
identify $\GG_{0}$ with $\uni(\GG_{0})\subset\GG_{1}$.

Let $r\in\{0,1,\ldots,\infty\}$. 
A $\eCr$-{\em groupoid} is a topological groupoid
$\GG$ such that both $\GG_{0}$ and $\GG_{1}$ are
$\eCe$-manifolds modeled on Banach spaces (over $\RRR$) \cite{Lan},
all the structure maps are of class $\eCr$, and the domain map
is a submersion. A $\eCr$-groupoid $\GG$ is finite-dimensional
if both $\GG_{0}$ and $\GG_{1}$ are finite-dimensional
$\eCr$-manifolds.
We will always assume that our
$\eCr$-manifolds are without boundary, but they may not be
Hausdorff. Observe that the assumption that the domain
map is a submersion implies that
the pull-back (i.e. the fibered product)
$\GG_{1}\times_{\GG_{0}}\GG_{1}$ has a natural
structure of a $\eCr$-manifold. Moreover, the codomain and the
composition map are also submersions, and the unit map is
an immersion.

An {\em \'{e}tale groupoid} is a topological groupoid with the domain
map a local homeomorphism. All the structure maps of an \'{e}tale
groupoid are local homeomorphisms, and all the structure maps of
an \'{e}tale $\eCr$-groupoid are local $\eCr$-diffeomorphisms.

A {\em continuous functor} between topological groupoids
$\phi:\HH\ra\GG$ is a functor such that both
$\phi_{0}:\HH_{0}\ra\GG_{0}$
and $\phi=\phi_{1}:\HH_{1}\ra\GG_{1}$  are continuous.
Analogously, one has the notion of a $\eCr$-{\em functor} between
$\eCr$-groupoids.

We will denote by $\Gpd$ (respectively $\Gpde$) the category of
topological (respectively \'{e}tale) groupoids and continuous functors.
We denote by $\Gpdr$ (respectively $\Gpder$) the category of $\eCr$-groupoids (respectively \'{e}tale $\eCr$-groupoids) and
$\eCr$-functors.

\begin{ex}  \label{Iex4}  \rm
(1) Let $X$ be a topological space. Then $X$ can be viewed as an
\'{e}tale groupoid, with $X_{0}=X_{1}=X$ and with all the structure maps
the identities. The category $\Top$ of topological spaces is a full
subcategory of $\Gpde$.
Any $\eCr$-manifold is an \'{e}tale $\eCr$-groupoid.

(2) Let $G$ be a topological group. Then $G$ is also a topological
groupoid with $G_{1}=G$ and $G_{0}=\{\ast\}$, i.e. an one-point space.
It is \'{e}tale if and only if the group is discrete. The category of
(discrete) groups $\Grp$ is a full subcategory of $\Gpde$.

(3) Let $X$ be a topological space and $G$ a topological group
acting continuously on $X$. Denote the action by
$$ \nu:X\times G\lra X\;.$$
There is a topological groupoid $G(X,\nu)$ (denoted also by $G(X)$)
associated to this action, with $G(X)_{0}=X$ and $G(X)_{1}=X\times G$.
The domain map of $G(X)$ is $\nu$, the codomain map is the second
projection, and the composition is given by
$$ (x',g')\com (x,g)=(x',g'g)\;.$$
The groupoid $G(X)$ is \'{e}tale if and only if $G$ is discrete.

(4) In Section \ref{secFol} we saw that the $\eCe$-manifold
$\GmqCe$ is the space of morphisms of an \'{e}tale
$\eCe$-groupoid, denoted again by $\GmqCe$. In particular,
the space of objects of $\GmqCe$ is $\RRR^{q}$. The
groupoid $\GmqCe$ is called the Haefliger groupoid of
dimension $q$.

More generally, if $X$ is a topological space, we define the \'{e}tale
groupoid $\Gm(X)$ of $X$ in an analogous way, so that the space of
objects of $\Gm(X)$ is $X$ and the space of morphisms of $\Gm(X)$ is
the space of germs of locally defined homeomorphisms of $X$.

If $M$ is a $\eCr$-manifold, the germs of locally defined
$\eCr$-diffeomorphisms of $M$ form an \'{e}tale $\eCr$-groupoid
$\Gm_{\eCr}(M)$, which is an open subgroupoid of $\Gm(M)$.
In particular, $\GmqCe=\Gm_{\eCe}(\RRR^{q})$.

(5) Let $M$ be a finite-dimensional
$\eCe$-manifold and $\cF$ a foliation on $M$.
Then there is a groupoid $\Hol(M,\cF)$, called the
{\em holonomy groupoid} (or the {\em groupoid of leaves})
of $\cF$, defined such that
the space of objects of $\Hol(M,\cF)$ is $M$, and the morphisms of
$\Hol(M,\cF)$ are the holonomy classes of paths inside the leaves of
$\cF$. There is a natural structure of a $\eCe$-manifold on
$\Hol(M,\cF)_{1}$ which makes $\Hol(M,\cF)$ into a $\eCe$-groupoid
(for details, see \cite{Win}). The groupoid of leaves of a foliation
reflects the transversal structure of the foliation, and is a good
replacement for the space of leaves of the foliation which can behave
very badly (see Example \ref{Iex3} (4)).
\end{ex}

Let $\GG$ be a topological groupoid. The {\em space of orbits}
$|\GG|$ of $\GG$ is the coequalizer
$$\CD
  \GG_{1} \cdrr{\cod}{\dom} \GG_{0} \cdr{}{} |\GG|
  \endCD\;.$$
Thus $|\GG|$ is the quotient space obtained from $\GG_{0}$ by
identifying two points whenever there is a morphisms in $\GG_{1}$
between them. Clearly, $|\oo|$ can be extended to a functor
$$ |\oo|:\Gpd\lra\Top\;.$$

Let $\GG$ be an \'{e}tale groupoid. There is a canonical
continuous functor
$$ \eee:\GG\lra\Gm(\GG_{0})\;,$$
called the {\em effect-functor}, which is given by
$$ \eee(g)=\germ_{\dom g}(\cod\com(\dom|_{U})^{-1})\;,$$
where $U$ is a small neighbourhood of $g\in\GG_{1}$.
In particular, $\eee_{0}$ is the identity on $\GG_{0}$.
The map $\eee$ is a local homeomorphism.
The image of $\eee$ is an open \'{e}tale subgroupoid $\Eff(\GG)$ of
$\Gm(\GG_{0})$, called the {\em effect} of $\GG$.
The \'{e}tale groupoid $\GG$ is called {\em effective}
if $\eee$ is injective on morphisms:
in this case $\eee$ provides
an isomorphism $\GG\cong\Eff(\GG)$.
The effective groupoids are also referred to as
$S$-atlases \cite{Est}.
If $\GG$ is an \'{e}tale $\eCr$-groupoid,
then $\Eff(\GG)$ is actually an open
$\eCr$-subgroupoid of $\Gm_{\eCr}(\GG_{0})$. Observe that $\Eff$ can
be extended to a functor
$$ \Eff:\Gpde\lra\Gpde\;,$$
whose image is equivalent to the full subcategory of effective
groupoids, and $\Eff^{2}=\Eff$, i.e. $\Eff$ is a projector.

\begin{ex}  \label{Iex5}  \rm
Let $X$ be a topological space. Denote by $\cT_{X}$ the set of
all transitions of $X$, i.e. of all homeomorphisms between open
subsets of $X$. A {\em pseudogroup} on $X$ is a subset
$\cP$ of $\cT_{X}$ such that
\begin{enumerate}
\item [(i)]   $\id_{U}\in\cP$ for any open subset $U$ of $X$,
\item [(ii)]  if $f,f'\in\cP$, then
              $f'\com f|_{f^{-1}(\dom(f'))}\in\cP$,
\item [(iii)] if $f\in\cP$, then $f^{-1}\in\cP$, and
\item [(iv)]  if $f\in\cT_{X}$ and $(U_{i})$ is an open
              cover of $\dom(f)$ such that $(f|_{U_{i}})\subset\cP$,
              then $f\in\cP$.
\end{enumerate}
For example, $\cT_{X}$ itself is a pseudogroup.
If $M$ is a $\eCr$-manifold, then all
$\eCr$-diffeomorphisms between open subsets of $M$ form a
pseudogroup $\eCr_{M}$ on $M$.

For a pseudogroup $\cP$ on $X$ one can construct the associated
\'{e}tale groupoid $\Gm_{\cP}(X)$ of germs of homeomorphisms in $\cP$,
in an analogous way as we constructed the groupoid $\GmqCe$
in Section \ref{secFol}. This groupoid is obviously effective.
One has $\Gm_{\cT_{X}}(X)=\Gm(X)$
and $\Gm_{\eCr_{M}}(M)=\Gm_{\eCr}(M)$.

Conversely, let $\GG$ be an effective groupoid. A subset
$U$ of $\GG_{1}$ is called elementary if it is open and both
$\dom|_{U}$ and $\cod|_{U}$ are injective. Define
$$ \Psi(\GG)=\{\,\cod\com(\dom|_{U})^{-1}\,|\,U\mbox{ elementary
   subset of }\GG_{1}\,\}\;.$$
It is easy to see that $\Psi(\GG)$ is a pseudogroup on $\GG_{0}$.

For any topological space $X$, these two constructions give a
bijective correspondence between pseudogroups on $X$ and
isomorphism classes of effective groupoids with the space
of objects $X$.
\end{ex}

A continuous functor (respectively a $\eCr$-functor)
$\phi:\HH\ra\GG$ between topological groupoids (respectively
$\eCr$-groupoids) is an {\em essential equivalence} (respectively
a $\eCr$-{\em essential equivalence}) if
\begin{enumerate}
\item [(i)]  the map
             $$ \cod\com\pr_{1}:\GG_{1}\times_{\GG_{0}}\HH_{0}
                \lra\GG_{0} $$
             is an open surjection (respectively a surjective
             $\eCr$-submersion), and
\item [(ii)] the square
             $$\CD
               \HH_{1}     \cdr{\phi_{1}}{}     \GG_{1}       \\
               \cdd{(\dom,\cod)}{}        \cdd{}{(\dom,\cod)} \\
               \HH_{0}\times\HH_{0}
                     \cdr{\phi_{0}\times\phi_{0}}{}
                                         \GG_{0}\times\GG_{0}
             \endCD$$
             is a pull-back of topological spaces (respectively
             of $\eCr$-manifolds).
\end{enumerate}
In (i) above, the pull-back $\GG_{1}\times_{\GG_{0}}\HH_{0}$
is the pull-back of the maps $\dom$ and $\phi_{0}$.
If $\HH$ and $\GG$ are \'{e}tale and $\phi:\HH\ra\GG$ is an essential
equivalence, then both $\phi_{0}$ are $\phi_{1}$ are local
homeomorphisms.

The {\em Morita equivalence} is an equivalence relation between
topological gro\-up\-oids, the smallest such that two topological
groupoids are Morita equivalent whenever
there exists an essential equivalence between them. Similarly,
one has the notion of a $\eCr$-{\em Morita equivalence}.
\vspace{4 mm}

\Rem
An essential equivalence as well as Morita equivalence are
in some literature referred to as weak equivalence. It can be
shown that two topological groupoids $\HH$ and $\GG$ are Morita
equivalent if and only if there exists a topological groupoid $\KK$
and essential equivalences $\phi:\KK\ra\GG$ and $\psi:\KK\ra\HH$
\cite{Moe3}.

\begin{ex}  \label{Iex6}  \rm
(1) Let $X$ be topological space, and let
$(V_{i})$ be an open cover of $X$.
We will now define an \'{e}tale groupoid
$\VV$ associated to the open cover $(V_{i})$
and an essential equivalence $\iota:\VV\ra X$.

Define the space of objects of $\VV$ as the disjoint union
$\VV_{0}=\coprod_{i} V_{i}$, and
let $\iota_{0}:\VV\ra X$ be the obvious local
homeomorphism. Define now $\VV_{1}$, $\iota_{1}$ and
the domain and the codomain map of $\VV$ by the following
pull-back:
$$\CD
  \VV_{1} \cdr{\iota_{1}}{}       X          \\
  \cdd{(\dom,\cod)}{}      \cdd{}{(\id,\id)} \\
  \VV_{0}\times\VV_{0}
         \cdr{\iota_{0}\times\iota_{0}}{}
                              X\times X
  \endCD$$
In fact, we can take $\VV_{1}=\coprod_{i,j}V_{i}\cap V_{j}$.
There is a unique way to define the composition on $\VV$, since
there is at most one morphism between arbitrary two objects of
$\VV$. It is easy to see that $\VV$ is indeed an \'{e}tale groupoid,
and it is obvious that $\iota$ is an essential equivalence.

If $X$ is a $\eCr$-manifold, then $\VV$ is a $\eCr$-groupoid
and $\iota$ is a $\eCr$-essential equivalence.

(2) Let $\cF$ be a foliation of codimension $q$ on a
$\eCe$-manifold $M$ of dimension $n$. Represent $\cF$ by
a $\GmqCe$-cocycle of submersions $(c_{ij})$
on $M$ with respect to an open cover $(V_{i})$, as we
explained in Section \ref{secFol}. Let $\VV$ be the \'{e}tale
$\eCe$-groupoid associated to the cover $(V_{i})$ and let
$\iota:\VV\ra M$ be the $\eCe$-essential equivalence,
as in (1).

Now the coproduct of the submersions $(c_{ij})$ defines
a $\eCe$-functor 
$$ \coprod_{ij}c_{ij}:\VV\lra\GmqCe\;,$$
since we have $\VV_{1}=\coprod_{ij}V_{i}\cap V_{j}$.
This functor is a submersion on morphisms and hence also on
objects. Observe that in fact the $\eCe$-functors
$\VV\ra\GmqCe$ are exactly the $\GmqCe$-cocycles of class $\eCe$
on $M$ with respect to the open cover $(V_{i})$, and that
the functors which are submersions on morphisms correspond
exactly to the cocycles of submersions, hence they represent
foliations on $M$.

(3) Let $\cF$ be a foliation of codimension $q$ on a
$\eCe$-manifold $M$ of dimension $n$. A {\em complete transversal}
of $\cF$ is a $\eCe$-immersion $T:N\ra M$ of a $\eCe$-manifold
$N$ of dimension $q$ which is transversal to the leaves of $\cF$
and such that $T(N)$ meets any leaf of $\cF$ in at least one point.
A complete transversal obviously always exists, and in fact one can
assume that $N$ is a disjoint union of copies of $\RRR^{q}$.

Choose a complete transversal $T:N\ra M$. Now one can define the
{\em holonomy groupoid} $\Hol_{T}(M,\cF)$ {\em reduced to} $T$
such that the space of objects of $\Hol_{T}(M,\cF)$ is $N$ and
the space of morphisms, the domain and the codomain map of
$\Hol_{T}(M,\cF)$ are given by the following pull-back:
$$\CD
  \Hol_{T}(M,\cF)_{1} \cdr{}{}    \Hol(M,\cF)_{1}  \\
  \cdd{(\dom,\cod)}{}          \cdd{}{(\dom,\cod )} \\
  N\times N       \cdr{T\times T}{}  M\times M
  \endCD$$
Here $\Hol(M,\cF)$ is the (unreduced) holonomy groupoid
(Example \ref{Iex4} (5)). In other words, if $y$ and $y'$ are
two points in $N$, then the morphisms of $\Hol_{T}(M,\cF)$
from $y$ to $y'$ are the holonomy classes of paths inside
the leaves of $\cF$ from $T(y)$ to $T(y')$. In particular,
there are no morphisms between $y$ and $y'$ if $T(y)$ and $T(y')$
do not lie on the same leaf. In fact, the holonomy class
of a path $\sigma$ from $T(y)$ to $T(y')$ in a leaf of
$\cF$ can be viewed as the germ of a diffeomorphism
between small neighbourhoods of $y$ and $y'$ in $N$,
which is provided by the foliation on a small neighbourhood of
$\sigma([0,1])$. Hence we can view $\Hol_{T}(M,\cF)$ as an
open subgroupoid of $\Gm_{\eCe}(N)$. In particular,
$\Hol_{T}(M,\cF)$ is an effective $\eCe$-groupoid.

It is obvious that the functor $\Hol_{T}(M,\cF)\ra\Hol(M,\cF)$
is a $\eCe$-essential equivalence.
In particular, if $T':N'\ra M$ is another complete transversal of
$\cF$, then the groupoids
$\Hol_{T}(M,\cF)$ and $\Hol_{T'}(M,\cF)$ are $\eCe$-Morita
equivalent. In fact, the $\eCe$-map
$T''=T\amalg T':N\amalg N'\ra M$ is
also a complete transversal, and the inclusions
$\Hol_{T}(M,\cF)\hookrightarrow\Hol_{T''}(M,\cF)$ and
$\Hol_{T'}(M,\cF)\hookrightarrow\Hol_{T''}(M,\cF)$ are clearly
essential equivalences of class $\eCe$ between \'{e}tale gro\-up\-oids.
In particular, the holonomy groupoid $\Hol(M,\cF)$ is $\eCe$-Morita
equivalent to an effective $\eCe$-groupoid.

(4) Let $\cF$ be a foliation of codimension $q$ on a
$\eCe$-manifold $M$ of dimension $n$. Choose an atlas
$\varphi=(\varphi_{i}:U_{i}\ra\RRR^{n-q}\times\RRR^{q})_{i\in I}$
of surjective charts for $\cF$. Put $s_{i}=\pr_{2}\com\varphi_{i}$,
let $N=\coprod_{i\in I}\RRR^{q}$, and let $T:N\ra M$ be the complete
transversal of $\cF$ given by $T(y,i)=\varphi_{i}^{-1}(0,y)$.
For any $i,j\in I$ choose an open cover $\cW_{ij}$ of
$U_{i}\cap U_{j}$ such that for any $W\in\cW_{ij}$
there exists the $\eCe$-diffeomorphism
$s_{ij}^{W}:s_{j}(W)\ra s_{i}(W)$ with
$s_{ij}^{W}\com s_{j}|_{W}=s_{i}|_{W}$.
(see Section \ref{secFol}). We have the $\eCe$-transition
$f_{ij}^{W}:s_{j}(W)\times\{j\}\ra s_{i}(W)\times\{i\}$ of $N$
defined by
$$ f_{ij}^{W}(y,j)=(s_{ij}^{W}(y),i)\;.$$
Let $\Psi_{\varphi}(M,\cF)$ be the pseudogroup generated by
$\{\,f_{ij}^{W}\,|\, i,j\in I,\, W\in\cW_{ij}\,\}$,
i.e. the smallest subpseudogroup of $\eCe_{N}$
which includes $(f_{ij}^{W})$.
Observe that $\Psi_{\varphi}(M,\cF)$
depends only on the choice of $\varphi$.
Clearly there is a natural $\eCe$-isomorphism
$$ \Gamma_{\Psi_{\varphi}(M,\cF)}(N)\cong\Hol_{T}(M,\cF)\;.$$

Conversely, let $\cP$ be a pseudogroup on a $\eCe$-manifold $N$
of dimension $q$, generated by a countable set of $\eCe$-transitions
$(f_{i})_{i=1}^{\infty}$. Let $V$ be the open subset
of $N\times\RRR\times\RRR$ given by
$$ V=(N\times\RRR\times (0,1))\cup\bigcup_{i=1}^{\infty}
   (\dom(f_{i})\times (i,i+1)\times (0,3))\;.$$
There is a natural foliation on $V$ of codimension $q$
given by the first projection. Let $M$ be the $\eCe$-manifold
of dimension $q+2$, obtained as the quotient of $V$
by identifying
$(y,t,t')$ with $(f_{i}(y),t,t'-2)$
for any $i=1,2,\ldots\,$, $y\in\dom(f_{i})$, $t\in (i,i+1)$ and
$t'\in (2,3)$. Observe that the foliation on $V$ induces
a foliation $\cF$ on $M$ of codimension $q$.
For any $y\in N$ denote by $T(y)\in M$ the equivalence class
of the point $(y,0,1/2)\in V$. The $\eCe$-map
$T:N\ra M$ thus defined is a complete transversal of $\cF$.
It is easy to check that the effective $\eCe$-groupoid
$\Gamma_{\cP}(N)$ is $\eCe$-isomorphic to the reduced
holonomy groupoid $\Hol_{T}(M,\cF)$ of $\cF$.
The idea for this construction was communicated to us by
J. Pradines, who attributed it to G. Hector.
\end{ex}

Let $\HH$ be a topological groupoid and let
$E$ be a topological space. A {\em right $\HH$-action}
on $E$ is a pair $(\xW,\nu)$, where $\xW:E\ra\HH_{0}$ and 
$\nu:E\times_{\HH_{0}}\HH_{1}\ra E$ are continuous maps
(we will write $e\cdot h=\nu(e,h)$) which satisfies
\begin{enumerate}
\item [(i)]   $\xW(e\cdot h)=\dom h$,
\item [(ii)]  $e\cdot 1_{\xW(e)}=e$, and
\item [(iii)] $(e\cdot h)\cdot h'=e\cdot(h\com h')$
\end{enumerate}
for any $e\in E$ and $h,h'\in\HH_{1}$ with $\xW(e)=\cod h$
and $\dom h=\cod h'$. If we have a right $\HH$-action $(\xW,\nu)$
on the space $E$, we say that $\HH$ acts on $E$ with respect to
the map $\xW$. If $\HH$ is a $\eCr$-groupoid and $M$
a $\eCr$-manifold, then a right $\HH$-action $(\xW,\nu)$
on $M$ is of class $\eCr$ if both $\xW$ and $\nu$ are $\eCr$-maps.

A {\em right $\HH$-space} is a triple $(E,\xW,\nu)$, where $E$ is a
topological space and $(\xW,\nu)$ a right $\HH$-action on $E$.
We will often write simply $E$ or $(E,\xW)$ for a right $\HH$-space
$(E,\xW,\nu)$. 

The {\em orbit space} $E/\HH$ of a right $\HH$-space 
$(E,\xW,\nu)$ is the coequalizer
$$\CD
  E\times_{\HH_{0}}\HH_{1} \cdrr{\nu}{\pr_{1}} E \cdr{q_{E}}{} E/\HH
  \endCD\;.$$
In other words, the orbit space is obtained as the quotient of $E$
by identifying two points $e,e'\in E$ if there exists a morphism
$h$ of $\HH$ such that $\cod h=\xW(e)$ and $e\cdot h=e'$.
The map $q_{E}$ is open since the domain map of $\HH$ is open.
A subspace $Y\subset E$ is called {\em $\HH$-invariant} if 
$q_{E}^{-1}(q_{E}(Y))=Y$.

A right $\HH$-space $E$ is called {\em $\HH$-compact}
(respectively $\HH$-{\em Hausdorff}, $\HH$-{\em connected})
if the associated orbit space $E/\HH$ is
compact (respectively Hausdorff, connected).
An $\HH$-{\em connected component}
of $E$ is the inverse image $q_{E}^{-1}(Z)$ of a connected component 
$Z$ of $E/\HH$.

Any topological groupoid $\HH$ acts canonically on the right on its 
space of objects $\HH_{0}$, and $\HH_{0}/\HH=|\HH|$.
Conversely, if $(E,\xW,\nu)$ is a right $\HH$-space, there is the
associated topological groupoid $\HH(E)$ with $\HH(E)_{0}=E$, 
$\HH(E)_{1}=E\times_{\HH_{0}}\HH_{1}$, $\dom=\nu$, $\cod=\pr_{1}$
and with the obvious composition.
If $\HH$ is \'{e}tale, $\HH(E)$ is also \'{e}tale. If $\HH$ is a
$\eCr$-groupoid, $M$ a $\eCr$-manifold and $(\xW,\nu)$ a right
$\HH$-action of class $\eCr$ on $M$, then $\HH(E)$ is a
$\eCr$-groupoid. 
Analogously one can consider left actions.

Assume that $\GG$ is
a topological groupoid and $X$ is a topological space. A
$\GG$-bundle over $X$ is a topological space $E$, equipped with
a left $\GG$-action $(p,\mu)$
(which we write by $g\cdot e=\mu(g,e)$)
and a continuous map $\xW:E\ra X$
such that $\GG$ acts along the fibers of $\xW$, i.e.
$$ \xW(g\cdot e)=\xW(e) $$
for any $g\in\GG_{1}$ and $e\in E$ with $\dom g=p(e)$.
We denote such a $\GG$-bundle by $(E,p,\xW,\mu)$ or $(E,p,\xW)$, or
simply by $E$. If $\GG$ is a $\eCr$-groupoid and $M$ a
$\eCr$-manifold, a  $\GG$-bundle $(E,p,\xW,\mu)$ over $M$ is of
class $\eCr$ if $E$ is a $\eCr$-manifold, $\xW$ a $\eCr$-map and
$(p,\mu)$ a $\GG$-action of class $\eCr$.

A $\GG$-equivariant map between $\GG$-bundles $(E,p,\xW)$
and $(E',p',\xW')$ over $X$ is a map $\alpha:E\ra E'$ satisfying
$\xW'\com\alpha=\xW$, $p'\com\alpha=p$ and
$$ \alpha(g\cdot e)=g\cdot\alpha(e) $$
for all $g\in\GG_{1}$ and $e\in E$ with $\dom g=p(e)$.
Two $\GG$-bundles over $X$ are isomorphic if there exists
a $\GG$-equivariant homeomorphism between them.
If $E$ is a $\GG$-bundle over $X$, we will often write the
isomorphism class of $E$ again by $E$.

A $\GG$-bundle $(E,p,\xW,\mu)$ over $X$ is called {\em transitive}
if the map $\xW$ is an open surjection, and the map
$$ (\mu,\pr_{2}):\GG_{1}\times_{\GG_{0}}E\lra E\times_{X}E $$
is a surjection, i.e. if $\GG$ acts transitively along the
fibers of $\xW$.
A $\GG$-bundle $(E,p,\xW,\mu)$ over $X$ is called {\em principal}
if the map $\xW$ is an open surjection, and the map
$(\mu,\pr_{2})$ is a homeomorphism
(hence the action of $\GG$ is free and transitive along the fibers
of $\xW$).
If $\GG$ is an \'{e}tale groupoid and $(E,p,\xW)$ a
principal $\GG$-bundle over $X$,
then the map $\xW$ is a local homeomorphism.
If $\GG$ is a $\eCr$-groupoid and $M$ a
$\eCr$-manifold, a  $\GG$-bundle $(E,p,\xW,\mu)$ over $M$ of
class $\eCr$ is $\eCr$-{\em principal} if $\xW$ is a surjective
submersion and the map $(\mu,\pr_{2})$ is a $\eCr$-diffeomorphism.
If $\GG$ is an \'{e}tale $\eCr$-groupoid and $(E,p,\xW)$ a
$\eCr$-principal $\GG$-bundle over $M$,
then the map $\xW$ is a local $\eCr$-diffeomorphism.

Let $\GG$ be a topological groupoid.
Then any $\GG$-equivariant map between principal $\GG$-bundles
is a homeomorphism. If $\GG$ is a $\eCr$-groupoid, any
$\GG$-equivariant map of class $\eCr$ between $\eCr$-principal
$\GG$-bundles is a diffeomorphism of class $\eCr$.

\begin{ex}  \label{Iex6a}  \rm
Let $X$ and $Y$ be topological spaces. An action of $Y$ on $X$ is
just a continuous map $X\ra Y$. If $p:X\ra Y$ is continuous map,
then $(X,p,id_{X})$ is a principal $Y$-bundle over $X$.

Conversely, let $(E,p,\xW)$ be a principal $Y$-bundle over $X$. Then
$\xW$ is a homeomorphism. In particular, $\xW$ is a $Y$-equivariant map
between $(E,p,\xW)$ and $(X,p\com \xW^{-1},id_{X})$. We can thus identify
the isomorphism classes of principal $Y$-bundles over $X$ with the
continuous maps from $X$ to $Y$.
\end{ex}

\section{Haefliger Structures}  \label{secHaeStr}

In Section \ref{secFol} we defined the notion of a $\GmqCe$-cocycle
of class $\eCe$ on an open cover of a $\eCe$-manifold $M$.
This notion can be
easily generalized by replacing $\GmqCe$ by an \'{e}tale groupoid and
$M$ by a topological space \cite{Hae}.

Let $\GG$ be an \'{e}tale groupoid and $X$ a topological space.
Let $\cU=(U_{i})_{i\in I}$ be an open cover of $X$.
A $\GG$-{\em cocycle} on $\cU$ is a family of continuous maps
$$ c=(c_{ij}:U_{i}\cap U_{j}\lra\GG_{1})_{i,j\in I} $$
such that
\begin{enumerate}
\item [(i)]   $c_{ii}(x)\in\GG_{0}\subset\GG_{1}$ for any $x\in 
              U_{i}$,
\item [(ii)]  $\dom c_{ij}(x)=c_{jj}(x)$, $\cod c_{ij}(x)=c_{ii}(x)$
              for any $x\in U_{i}\cap U_{j}$, and
\item [(iii)] $c_{ij}(x)\com c_{jk}(x)=c_{ik}(x)$ for any
              $x\in U_{i}\cap U_{j}\cap U_{k}$.
\end{enumerate}
In particular, this implies that $c_{ij}(x)=c_{ji}(x)^{-1}$ for
$x\in U_{i}\cap U_{j}$. In \cite{Hae}, a $\GG$-cocycle on $\cU$
is referred to as 1-cocycle on $\cU$ with values in $\GG$.
Denote by $Z^{1}(\cU,\GG)$ the set of all $\GG$-cocycles on $\cU$.

Two $\GG$-cocycles $c, c'\in Z^{1}(\cU,\GG)$ are
{\em cohomologious} if there exists a family of continuous maps
$$ b=(b_{i}:U_{i}\lra\GG_{1})_{i\in I} $$
which intertwines $c$ and $c'$, i.e.
$\dom\com b_{i}=c_{ii}$, $\cod\com b_{i}=c'_{ii}$ and
$$ c'_{ij}(x)\com b_{j}(x)=b_{i}(x)\com c_{ij}(x)
   \;\;\;\;\;\;\;\; x\in U_{i}\cap U_{j}\;.$$
This is an equivalence relation on $Z^{1}(\cU,\GG)$, and the
set of equivalence classes (called the cohomology classes) of
$\GG$-cocycles on $\cU$ is denoted by
$$ H^{1}(\cU,\GG)\;.$$

Now assume that $\cV=(V_{k})_{k\in K}$ is another open cover of $X$,
finer than $\cU$. Hence we can choose a function
$\tau:K\ra I$ such that $V_{k}\subset U_{\tau(k)}$ for any $k\in K$.
This function induces a function
$$ \tau^{\#}:Z^{1}(\cU,\GG)\lra Z^{1}(\cV,\GG) $$
which is given by restriction, i.e.
$$ \tau^{\#}(c)=
   (c_{\tau(k)\tau(l)}|_{V_{k}\cap\, V_{l}})_{k,l\in K}\;.$$
This function maps cohomologious cocycles into cohomologious cocycles,
hence induces a function
$$ \cV/\cU:H^{1}(\cU,\GG)\lra H^{1}(\cV,\GG)\;.$$
This function does not depend on the choice of $\tau$. Indeed,
if $\upsilon:K\ra I$ is another function with $V_{k}\subset
U_{\upsilon(k)}$, and if $c\in Z^{1}(\cU,\GG)$, then
the family
$$ (c_{\upsilon(k)\tau(k)}|_{V_{k}}:V_{k}\lra\GG_{1})_{k\in K} $$
obviously intertwines $\tau^{\#}(c)$ and
$\upsilon^{\#}(c)$.

Observe that we have $\cU/\cU=\id$ and
$\cW/\cV\com\cV/\cU=\cW/\cU$ for a refinement $\cW$ of $\cV$.
The set of all open covers of $X$ is partially ordered, in the
sense that $\cV\leq\cU$ if $\cV$ is a refinement of $\cU$.
This order is filtered. Hence we can define the filtered colimit
$$ H^{1}(X,\GG)=\lim_{\ra_{\,\cU}}H^{1}(\cU,\GG)\;.$$
In other words, the elements of $H^{1}(X,\GG)$ are the
equivalence 
classes of $\GG$-cocycles on open covers of $X$, where two
cocycles $c\in Z^{1}(\cU,\GG)$ and 
$c'\in Z^{1}(\cV,\GG)$ are equivalent if there exists
a common refinement $\cW$ of $\cU$ and $\cV$ such that
the restrictions of $c$ and $c'$ on $\cW$ are cohomologious.
The elements of $H^{1}(X,\GG)$ are called the
{\em Haefliger $\GG$-structures} on $X$.

Let $c$ be a $\GG$-cocycle on an open cover $\cU=(U_{i})_{i\in I}$
of $X$. We will define a principal $\GG$-bundle $(\Sigma(c),p,\xW)$ 
over $X$ associated to the cocycle $c$. Put
$$ \bar{\Sigma}(c)=\GG_{1}\times_{\GG_{0}}\coprod_{i\in I}U_{i}=
   \{\,(g,x,i)\,|\,g\in\GG_{1},\,i\in I,\,x\in U_{i},\,
   \dom g=c_{ii}(x)\,\}\;.$$
Now define $\Sigma(c)$ to be the quotient space of the space
$\bar{\Sigma}(c)$ by identifying $(g,x,i)$ with
$(g\com c_{ij}(x),x,j)$ for all
$x\in U_{i}\cap U_{j}$ and $g\in\GG_{1}$ with $\dom g=c_{ii}(x)$.
This is clearly an equivalence relation on $\bar{\Sigma}(c)$, and
one can easily see that the quotient projection
$\bar{\Sigma}(c)\ra \Sigma(c)$ is a local homeomorphism. We denote
by $[g,x,i]\in \Sigma(c)$ the equivalence class of
$(g,x,i)\in\bar{\Sigma}(c)$.
The maps $\xW:\Sigma(c)\ra X$ and $p:\Sigma(c)\ra\GG_{0}$ are
given by
$$ \xW([g,x,i])=x\;\;\;\;\mbox{and}\;\;\;\; p([g,x,i])=\cod g\;.$$
The action of $\GG$ on $\Sigma(c)$ with
respect to $\xW$ and along the fibers of $p$ is given by the
composition in $\GG$, i.e.
$$ g'\cdot [g,x,i]=[g'\com g,x,i]\;.$$
It is easy to check that $(\Sigma(c),p,\xW)$ is indeed a principal
$\GG$-bundle over $X$. Thus we get a map $\Sigma$ which assigns
the principal $\GG$-bundle $\Sigma(c)$ over $X$ to a
$\GG$-cocycle $c$ on an open cover of $X$.

\begin{prop}  \label{Iprop7}
Let $\GG$ be an \'{e}tale groupoid and $X$ a topological space.
The map $\Sigma$ induces a bijective correspondence between the
Haefliger $\GG$-structures on $X$ and the isomorphism classes
of principal $\GG$-bundles over $X$.
\end{prop}
\Proof
First we have to prove that $\Sigma(c)$ depends up to an
isomorphism only on
the Haefliger $\GG$-structure represented by $c$.

To see this, assume first that $c$ is a $\GG$-cocycle on an open
cover $\cU=(U_{i})_{i\in I}$ of $X$, and $\cV=(V_{k})_{k\in K}$
is a refinement of $\cU$. Choose $\tau:K\ra I$ with $V_{k}\subset
U_{\tau(k)}$. Then we define
$\alpha:\Sigma(\tau^{\#}(c))\ra \Sigma(c)$ by
$$ \alpha([g,x,k])=[g,x,\tau(k)]\;.$$
This is clearly a well-defined continuous map, and since it is
also $\GG$-equivariant, it is a homeomorphism.

Next, assume that $c,c'\in Z^{1}(\cU,\GG)$ are cohomologious.
Let $b$ be a family of maps which intertwines $c$ and $c'$.
Then we define $\beta:\Sigma(c')\ra \Sigma(c)$ by
$$ \beta([g,x,i])=[g\com b_{i}(x),x,i]\;.$$
Again, this is a $\GG$-equivariant homeomorphism.
This proves that $\Sigma$ induces a map from $H^{1}(X,\GG)$ to 
the set of isomorphism classes of principal $\GG$-bundles over
$X$. We will denote the induced map again by $\Sigma$.

To prove that the correspondence is indeed bijective,
we will construct the inverse for $\Sigma$. Let $(E,p,\xW,\mu)$
be a principal $\GG$-bundle over $X$. Since $\xW$ is a local
homeomorphism, we can choose an open cover
$\cU=(U_{i})_{i\in I}$
of $X$ and continuous sections $t_{i}:U_{i}\ra E$ of $\xW$.
Since $E$ is principal, the map $(\mu,\pr_{2}):\GG_{1}
\times_{\GG_{0}}E\ra E\times_{X}E$ is a homeomorphism, with
the inverse of the form $(\theta,\pr_{2})$. Now we define
$$ c=(c_{ij}:U_{i}\cap U_{j}\lra\GG_{1})_{i,j\in I} $$
by
$$ c_{ij}(x)=\theta(t_{i}(x),t_{j}(x))\;.$$
It is clear that the family $c$ is a $\GG$-cocycle on $\cU$,
and that this construction provides an inverse for $\Sigma$.
\eop

Let $\GG$ be an \'{e}tale $\eCr$-groupoid and $M$ a
$\eCr$-manifold. Then a $\GG$-cocycle on an open cover
$\cU$ of $M$ is of class $\eCr$ if all the maps in the
cocycle are of class $\eCr$.  A Haefliger $\GG$-structure
is of class $\eCr$ if it can be represented by a cocycle
of class $\eCr$. Observe that the families
which intertwine $\GG$-cocycles of class $\eCr$ are necessarily
of class $\eCr$. In particular, any cocycle which represents
a Haefliger $\GG$-structure of class $\eCr$ is itself of
class $\eCr$. If $c$ is a $\GG$-cocycle of class $\eCr$,
then $\Sigma(c)$ is clearly a $\eCr$-principal $\GG$-bundle
over $M$.

\begin{cor}  \label{Icor8}
Let $\GG$ be an \'{e}tale $\eCr$-groupoid and let $M$ be a
$\eCr$-manifold. Then the map $\Sigma$ induces a bijective 
correspondence between the Haefliger $\GG$-structures
of class $\eCr$ on $M$ and the isomorphism classes
of $\eCr$-principal $\GG$-bundles over $M$.
\end{cor}
\Proof
The proof goes exactly
the same as the proof of Proposition \ref{Iprop7}, since
all the constructions in that proof have a natural
$\eCr$-structure.
\eop

\begin{ex}  \label{Iex9}  \rm
(1) Let $M$ be a finite-dimensional $\eCe$-manifold.
As we have seen in Section \ref{secFol},
a foliation $\cF$ on $M$ can be represented by a
$\GmqCe$-cocycle of $\eCe$-submersions on an open cover of
$M$. It is not difficult to verify that two such cocycles
of submersions determine the same foliation precisely if
they determine the same Haefliger $\GmqCe$-structure on $M$.
In other words, the foliations on $M$ are in a natural
bijective correspondence with the Haefliger $\GmqCe$-structures
of $\eCe$-submersions on $M$, i.e. with those Haefliger 
$\GmqCe$-structures on $M$ which can be represented by cocycles
of $\eCe$-submersions. Note that again any cocycle which 
represents a Haefliger $\GmqCe$-structure of $\eCe$-submersions
is itself a cocycle of $\eCe$-submersions.

We say that a $\eCe$-principal $\GG$-bundle $(E,p,\xW)$
over $M$ is {\em submersive}
if the map $p$ is a submersion.
One can easily see from the construction that $\Sigma(c)$ is
a submersive $\GG$-bundle precisely if $c$ is a $\GG$-cocycle
of $\eCe$-submersions on an open cover of $M$.

With this, Corollary \ref{Icor8} implies that the foliations
on $M$ are in a natural bijective correspondence with the
isomorphism classes of the $\eCe$-principal submersive
$\GmqCe$-bundles over $M$.

(2) Let $M$ be a $\eCe$-manifold of dimension $n$ and $\cF$ a
foliation on $M$ of codimension $q$. We can choose a $\eCe$-atlas
$(\varphi_{i}:U_{i}\ra\RRR^{n-q}\times\RRR^{q})_{i\in I}$
representing $\cF$ such that each $\varphi_{i}$ is surjective.
Put $N=\coprod_{i\in I}\RRR^{q}$ and define $T:N\ra M$ by
$$ T(y,i)=\varphi_{i}^{-1}(0,y)\;.$$
Clearly $T$ is a complete transversal. Recall that
$\Hol_{T}(M,\cF)$ is an effective $\eCe$-groupoid, i.e.
it can be viewed as an open subgroupoid of $\Gm_{\eCe}(N)$.
Let $s_{i}:U_{i}\ra N$ be the submersion given by
$s_{i}(x)=(\pr_{2}(\varphi_{i}(x)),i)$. For any
$x\in U_{i}\cap U_{j}$ there exists an open neighbourhood
$W\subset U_{i}\cap U_{j}$ of $x$ and a diffeomorphism
$s^{W}_{ij}:s_{j}(W)\ra s_{i}(W)$ with $s_{i}|_{W}=
s^{W}_{ij}\com s_{j}|_{W}$. Just like in Section \ref{secFol}
we obtain now a $\Gm_{\eCe}(N)$-cocycle $(c_{ij})$ on $(U_{i})$
given by $c_{ij}(x)=\germ_{s_{j}(x)}s^{W}_{ij}$.
But observe that $c_{ij}(x)$ in fact belongs to 
$\Hol_{T}(M,\cF)\subset\Gm_{\eCe}(N)$.
In other words, we obtain a $\Hol_{T}(M,\cF)$-cocycle on $(U_{i})$,
hence the corresponding Haefliger $\Hol_{T}(M,\cF)$-structure
of $\eCe$-submersions on $M$, or equivalently, the corresponding
submersive $\eCe$-principal $\Hol_{T}(M,\cF)$-bundle on $M$.
\end{ex}

Let $\GG$ be an \'{e}tale groupoid and let $\ccc$ be a
Haefliger $\GG$-structure over a topological space $X$.
Motivated by Example \ref{Iex9} (1), one can define the leaves and
the holonomy groups of the leaves of $\ccc$, and view $\ccc$ as a
generalized kind of a foliation on the topological space $X$.

Let $\cU=(U_{i})_{i\in I}$ be an open cover
of $X$ and $c\in Z^{1}(\cU,\GG)$ such that $c$ represents $\ccc$.
Then we define the plaques of $c$ in $U_{i}$ to be the
connected components of the fibers of $c_{ii}$. Since $\GG$
is \'{e}tale, the maps $c_{ij}$ guarantee that the plaques globally
amalgamate into connected subsets of $X$, called the leaves of $c$,
which forms a partition of $X$.
More precisely, the sets of the form
$$ c_{ii}^{-1}(a)\cap U\;,$$
where $a\in\GG_{0}$, $i\in I$
and $U$ is an open subset of $X$, form a basis for a (finer) topology
on $X$, called the leaf topology associated to $c$.
The connected components
of $X$ equipped with the leaf topology are the leaves of $c$.
It is easy to see that the partition of $X$ on leaves depends only
on the Haefliger $\GG$-structure $\ccc$.

In \cite{Hae}, A. Haefliger defined the holonomy group of a leaf
of a Haefliger $\GG$-structure on $X$. Then he proved the
Haefliger-Reeb-Ehresmann stability theorem for such a structure,
which says that under some conditions, a compact leaf with finite
holonomy group of a Haefliger $\GG$-structure has an arbitrary
small neighbourhood of compact leaves. This theorem
generalizes the Reeb stability theorem \ref{Itheo1},
in the sense of Example \ref{Iex9} (1). We shall discuss this
in more detail in Chapters \ref{chapHilSkaMap} and
\ref{chapSta}, where we will generalize the
Haefliger-Reeb-Ehresmann stability theorem, as well as the
Thurston-Reeb stability theorem (Theorem \ref{Itheo2}),
to Hilsum-Skandalis maps between topological groupoids.
As we will show, these maps can be viewed as generalized
foliations on topological groupoids.

\section{Fundamental Group of a Topological Groupoid}
\label{secFunGro}

\begin{dfn}  \label{Idfn10}  \rm
Let $\HH$ be a topological groupoid and $(E,\xW)$ a right
$\HH$-space. An $\HH$-{\em path} in $(E,\xW)$ from a point
$e\in E$ to a point $e'\in E$
consists of a sequence
$(\sigma_{i})_{i=0}^{n}$ of paths in $E$ and a sequence
$(h_{i})_{i=1}^{n}$ 
of elements of $\HH_{1}$ such that $\sigma_{0}(0)=e$,
$h_{i}\in\HH(\xW(\sigma_{i-1}(1)),\xW(\sigma_{i}(0)))$,
$\sigma_{i}(0)\cdot h_{i}=\sigma_{i-1}(1)$ for any
$i=1,2,\ldots n$, and $\sigma_{n}(1)=e'$. We will
denote such an $\HH$-path as
$$ \sigma_{n}\cdot h_{n}\cdot\ldots\cdot h_{1}
   \cdot\sigma_{0}\;. $$
An $\HH$-{\em loop} in $(E,\xW)$ with a base-point
$e_{0}\in E$ is an $\HH$-path in $(E,\xW)$ from $e_{0}$ to
$e_{0}$.
\end{dfn} 
\Rem
Observe that the $\HH$-paths can be concatenated in the
obvious way. Analogously, one can consider $\GG$-paths in a
left $\GG$-space, for a topological groupoid $\GG$.

\begin{prop}  \label{Iprop11}
Let $\HH$ be a topological groupoid and $(E,\xW)$ a locally
path-connected right $\HH$-space. Then $e,e'\in E$ lie in the
same $\HH$-connected component of $E$ if and only if there
exists an $\HH$-path from $e$ to $e'$.
\end{prop}
\Proof
Denote by $q=q_{E}:E\ra E/\HH$ the quotient map on the orbit
space of $E$.
Any $\HH$-path in $E$ from $e$ to $e'$ clearly induces a
path in $E/\HH$ from $q(e)$ to $q(e')$. 

Conversely, assume that $e$ and $e'$ lie in the same 
$\HH$-connected component of $E$. Hence $q(e)$ and 
$q(e')$ lie in the same connected component of $E/\HH$.
Since $E$ is locally path-connected, so is
$E/\HH$ and thus there exists a path
$\sigma:[0,1]\ra E/\HH$ from 
$q(e)$ to $q(e')$. Moreover, we can choose an open cover
$(U_{j})_{j\in J}$ of $E$, consisting of the path-connected
sets.

Now $(q(U_{j}))_{j\in J}$ is
an open cover of $E/\HH$, and by compactness of
$\sigma([0,1])$
we can choose a partition 
$0=t_{0}<t_{1}<\ldots<t_{n}<t_{n+1}=1$ and 
$j_{0},j_{1},\ldots,j_{n}\in J$
such that
$$ \sigma([t_{i},t_{i+1}])\subset q(U_{j_{i}}) $$
for all $i=0,1,\ldots,n$, with $e\in U_{j_{0}}$ and
$e'\in U_{j_{n}}$. We can find $e_{i},e'_{i}\in U_{j_{i}}$
with $q(e_{i})=\sigma(t_{i})$
and $q(e'_{i})=\sigma(t_{i+1})$, for all $i=0,1,\ldots,n$.
In particular, we can take $e_{0}=e$ and $e'_{n}=e'$. 
Since $U_{j_{i}}$ is path-connected,
there exists a path $\sigma_{i}$ in $U_{j_{i}}$ from
$e_{i}$ to $e'_{i}$ $(i=0,1,\ldots,n)$, and since 
$q(e_{i})=q(e'_{i-1})$, there exists
$h_{i}\in\HH(\xW(e'_{i-1}),\xW(e_{i}))$ with
$e_{i}\cdot h_{i}=e'_{i-1}$, 
for any $i=1,2,\ldots,n$. This gives an $\HH$-path
from $e$ to $e'$.
\eop

Let $\HH$ be a topological groupoid such that $\HH_{0}$
is locally path-connected. Since $\HH_{0}$ is a right
$\HH$-space, we can speak about $\HH$-paths in 
$\HH_{0}$. Let $b_{0}\in\HH_{0}$. An $\HH$-{\em loop} in
$\HH_{0}$ with the base-point $b_{0}$
thus consists of a sequence
$(\sigma_{i})_{i=0}^{n}$ of paths in $\HH_{0}$ and a
sequence $(h_{i})_{i=1}^{n}$ of elements of $\HH_{1}$
with $\sigma_{0}(0)=\sigma_{n}(1)=b_{0}$ and 
$h_{i}\in\HH(\sigma_{i-1}(1),\sigma_{i}(0))$ for all
$i=1,2,\ldots,n$.
We denote this $\HH$-loop by
$$ \sigma_{n}\cdot h_{n}\cdot\ldots\cdot h_{1}
   \cdot\sigma_{0}\;. $$
Denote by $\Omega(\HH,b_{0})$ the set of all
$\HH$-loops in $\HH_{0}$ with the base-point $b_{0}$.

We will now define an equivalence relation on
$\Omega(\HH,b_{0})$, which we will call simply
{\em equivalence}. It is the smallest equivalence
relation such that
\begin{itemize}
\item [(i)]  an $\HH$-loop 
             $$ \sigma_{n}\cdot h_{n}\cdot\ldots\cdot
                h_{1}\cdot\sigma_{0} $$
             is equivalent to the $\HH$-loop
             $$ \sigma_{n}\cdot h_{n}\cdot\ldots\cdot
                (h_{i+1}\com h_{i})\cdot\sigma_{i-1}\cdot
                \ldots\cdot h_{1}\cdot\sigma_{0}\;,$$
             if $\sigma_{i}$ is a constant path for some
             $0<i<n$, and
\item [(ii)] an $\HH$-loop
             $$ \sigma_{n}\cdot h_{n}\cdot\ldots\cdot
                h_{1}\cdot\sigma_{0} $$
             is equivalent to the $\HH$-loop
             $$ \sigma_{n}\cdot h_{n}\cdot\ldots
                \cdot h_{i+1}\cdot(\sigma_{i}\sigma_{i-1})
                \cdot\ldots\cdot h_{1}\cdot\sigma_{0}\;,$$
            if $h_{i}\in\HH_{0}$ for some $1\leq i\leq n$.
\end{itemize}
Here $\sigma_{i}\sigma_{i-1}$ denotes the usual
concatenation of $\sigma_{i-1}$ and $\sigma_{i}$.

A {\em deformation} of an $\HH$-loop
$\sigma_{n}\cdot h_{n}\cdot\ldots\cdot 
h_{1}\cdot\sigma_{0}$ to an $\HH$-loop
$\sigma'_{n}\cdot h'_{n}\cdot\ldots\cdot
h'_{1}\cdot\sigma'_{0}$
consists of homotopies
$$ D_{i}:[0,1]^{2}\lra\HH_{0} $$
from $D_{i}(0,\oo)=\sigma_{i}$ to $D_{i}(1,\oo)=\sigma'_{i}$
$(i=0,1,\ldots,n)$ and paths
$$ d_{i}:[0,1]\lra\HH_{1} $$
from $h_{i}$ to $h'_{i}$ $(i=1,2,\ldots,n)$
which satisfy 
\begin{enumerate}
\item [(a)] $\dom\com d_{i}=D_{i-1}(\oo,1)$ and
            $\cod\com d_{i}=D_{i}(\oo,0)$ for all
            $i=1,2,\ldots,n$, and
\item [(b)] $D_{0}([0,1],0)=D_{n}([0,1],1)=\{b_{0}\}$.
\end{enumerate}

Two $\HH$-loops in $\Omega(\HH,b_{0})$ are
{\em homotopic} if one can pass from one to another
in a sequence of deformations and equivalences.
With the multiplication induced by the concatenation,
the homotopy classes
of $\HH$-loops in $\HH_{0}$ with the base-point $b_{0}$
form a group
$$ \pi_{1}(\HH,b_{0})\;,$$
called the {\em fundamental group} of the topological
groupoid $\HH$ with the base-point $b_{0}$.
\vspace{4 mm}

\Rem
Proposition \ref{Iprop11} implies that the group
$\pi_{1}(\HH,b_{0})$ depends only on the $\HH$-connected
component $Z$ of $\HH_{0}$ with $b_{0}\in Z$.
In particular, if $|\HH|$ is connected, the fundamental
group of $\HH$ does not depend (up to an isomorphism) on
the choice of the base-point, and hence we shall denote
its isomorphism class by $\pi_{1}(\HH)$.

If $\HH$ is a locally path-connected topological space
$X$, this group is exactly the classical fundamental group
of $X$. Moreover, if $\HH$ is a discrete group
then $\pi_{1}(\HH)=\HH$. 

If $\HH$ is a suitable \'{e}tale groupoid,
$\pi_{1}(\HH)$ coincide with the 
fundamental group of the classifying space and of
the classifying topos of $\HH$ (see \cite{Moe1,Seg}).
If the groupoid $\HH$ is effective,
a similar (but obviously equivalent)
definition was given in \cite{Mol}. If in addition $\HH_{0}$
is simply-connected (but not necessarily connected),
it is not difficult to see that $\pi_{1}(\HH)$
is exactly the fundamental group of
$\HH$ as described by W. T. van Est \cite{Est}.
The construction of the fundamental group in $\cite{Est}$ is
extremely natural: it is the vertex group of the discrete
groupoid $\pi_{0}(\HH)$. Recall that
$\pi_{0}(\HH)_{0}=\pi_{0}(\HH_{0})$ and that
the morphisms of $\pi_{0}(\HH)$ are generated by
$\pi_{0}(\HH_{1})$ with respect to the partially defined
composition induced by
$\pi_{0}(\cmp):\pi_{0}(\HH_{1}\times_{\HH_{0}}\HH_{1})
\ra\pi_{0}(\HH_{1})$.
\vspace{4 mm}

Let $\GG$ be another topological groupoid
with $\GG_{0}$ locally path-connected, and
let $\phi:\HH\ra\GG$ be a continuous functor.
Denote $a_{0}=\phi_{0}(b_{0})$.
The functor $\phi$ gives a function
$$ \phi_{\#}:\Omega(\HH,b_{0})\lra\Omega(\GG,a_{0}) $$
defined by
$$ \phi_{\#}(\sigma_{n}\cdot h_{n}\cdot\ldots
             \cdot h_{1}\cdot\sigma_{0})=
   (\phi_{0}\com\sigma_{n})\cdot\phi(h_{n})\cdot\ldots
   \cdot\phi(h_{1})\cdot(\phi_{0}\com\sigma_{0})\;.$$
This function maps homotopic $\HH$-loops to homotopic
$\GG$-loops, hence it induces a map
$$ \phi_{\ast}:\pi_{1}(\HH,b_{0})\lra\pi_{1}(\GG,a_{0})\;,$$
which is clearly a homomorphism of groups.

If $\HH$ is a topological groupoid and
$(E,\xW)$ is a locally path-connected right $\HH$-space,
we define the {\em fundamental group} of $(E,\xW)$
with a base point $e_{0}\in E$ to be
the fundamental group $\pi_{1}(\HH(E),e_{0})$ of the
associated topological groupoid $\HH(E)$.
If $(E',\xW')$ is another locally path-connected right
$\HH$-space and $\alpha:E\ra E'$ an $\HH$-equivariant map
(i.e. $\xW'\com\alpha=\xW$ and $\alpha(e\cdot h)=\alpha(e)\cdot h$),
then $\alpha$ induces a continuous functor
$$ \HH(\alpha):\HH(E)\lra\HH(E') $$
in the obvious way. Hence we get the homomorphism of groups
$$ \alpha_{\ast}=\HH(\alpha)_{\ast}:
   \pi_{1}(\HH(E),e_{0})\lra\pi_{1}(\HH(E'),\alpha(e_{0}))\;.$$
Note that the $\HH$-loops in $\HH(E)_{0}$ are precisely the
$\HH$-loops in the right $\HH$-space $E$.

If $\alpha$ is a covering projection,
the unique path lifting property for $\alpha$
clearly implies the unique $\HH$-path lifting property.
Moreover, homotopic $\HH$-paths clearly
lift to homotopic $\HH$-paths.
In particular, $\alpha_{\ast}$ is in this case injective.
\vspace{4 mm}

\Not
Let $X$ be a topological space. For $x,x'\in X$
we will denote by $\bpi_{1}(X)(x,x')$ the set of homotopy
classes (with fixed end-points) of paths in $X$ from $x$
to $x'$. This defines a groupoid $\bpi_{1}(X)$ with
$\bpi_{1}(X)_{0}=X$ and with the composition induced by
the concatenation of paths. For example, we have
$\bpi_{1}(X)(x,x)=\pi_{1}(X,x)$.

Any continuous map $f:X\ra Y$ between topological
spaces induces a functor
$$ f_{\ast}=\bpi_{1}(f):\bpi_{1}(X)\lra\bpi_{1}(Y)\;,$$
which is given by the composition of paths with $f$.
\vspace{4 mm}

Let $X$ be a locally path-connected topological space,
equipped with a continuous right action of a discrete group
$G$. For any $g\in G$ we denote by $\hg$ the homeomorphism
$\oo\cdot g:X\ra X$. The groupoid $G(X)$ associated to this
action is \'{e}tale.

The following proposition gives a description of the
fundamental group of $G(X)$ in terms of $G$ and
$\bpi_{1}(X)$:

\begin{prop}  \label{Iprop12}
Let $X$ be a locally path-connected topological space,
let $x_{0}\in X$, and assume that $X$ is equipped with
a continuous right action of a discrete group $G$. Then
$$ \pi_{1}(G(X),x_{0})\cong\{\,(g,\vsig)\,|\,g\in G,\,
   \vsig\in\bpi_{1}(X)(x_{0},x_{0}\cdot g)\,\}\;,$$
with the multiplication given by
$$ (g',\vsig')(g,\vsig)=(g'g,\hg_{\ast}(\vsig')\,\vsig)\;.$$
In particular, there is a short exact sequence of groups
$$ 1\lra\pi_{1}(X,x_{0})\lra\pi_{1}(G(X),x_{0})\lra G\lra 1\;.$$
\end{prop}
\Rem
The group $\pi_{1}(G(X),x_{0})$ is thus isomorphic to the
fundamental group of the space $X\times_{G}\mbf{E}G$, where
$\mbf{E}G$ is the universal $G$-bundle (see for example \cite{AtiBot}).
\vspace{4 mm}

\Proof
Denote $\Upsilon=\{\,(g,\vsig)\,|\,g\in G,\,\vsig\in\bpi_{1}(X)
(x_{0},x_{0}\cdot g)\,\}$. It is easy to check that $\Upsilon$
is indeed a group. We define a function
$$ f:\Omega(G(X),x_{0})\lra\Upsilon $$
by
\begin{eqnarray*}
\lefteqn{
f(\sigma_{n}\cdot g_{n}\cdot\ldots\cdot\sigma_{2}
\cdot g_{2}\cdot\sigma_{1}\cdot g_{1}\cdot\sigma_{0})
        } \hspace{5 mm} \\
& = &
(g_{n}\ldots g_{2}g_{1},
(\hg_{1}\com\hg_{2}\com\ldots\com\hg_{n})_{\ast}(\vsig_{n})
\ldots
(\hg_{1}\com\hg_{2})_{\ast}(\vsig_{2})\,
(\hg_{1})_{\ast}(\vsig_{1})\,\vsig_{0})\;,
\end{eqnarray*}
where $\vsig_{i}$ denotes the homotopy class of $\sigma_{i}$.
It is easy to see that $f$ identifies homotopic
$G$-loops in $X$, hence induces a function
$$ \bar{f}:\pi_{1}(G(X),x_{0})\lra\Upsilon $$
which is clearly a surjective homomorphism.

Next, observe that there is an obvious
deformation between the $G$-loops
$$ \sigma_{n}\cdot g_{n}\cdot\ldots
   \cdot\sigma_{i}\cdot g_{i}\cdot\sigma_{i-1}\cdot
   \ldots\cdot g_{1}\cdot\sigma_{0} $$
and
$$ \sigma_{n}\cdot g_{n}\cdot\ldots
   \cdot\varpi(\sigma_{i}(1))\cdot g_{i}\cdot
   (\hg_{i}\com\sigma_{i})\,\sigma_{i-1}\cdot
   \ldots\cdot g_{1}\cdot\sigma_{0}\;,$$
where $\varpi(\sigma_{i}(1))$ denotes the constant path
with the image point $\sigma_{i}(1)$. This implies that the
$G$-loops
$$ \sigma_{n}\cdot g_{n}\cdot\ldots\cdot\sigma_{2}
   \cdot g_{2}\cdot\sigma_{1}\cdot g_{1}\cdot\sigma_{0} $$
and
$$ \varpi(x_{0})\cdot (g_{n}\ldots g_{2}g_{1})\cdot
   (\hg_{1}\com\hg_{2}\com\ldots\com\hg_{n}\com\sigma_{n})
   \ldots
   (\hg_{1}\com\hg_{2}\com\sigma_{2})\,
   (\hg_{1}\com\sigma_{1})\,\sigma_{0} $$
are homotopic, where $\varpi(x_{0})$ denotes the constant
path with the image point $x_{0}$. In particular, if
$f(\sigma_{n}\cdot g_{n}\cdot\ldots\cdot\sigma_{2}
\cdot g_{2}\cdot\sigma_{1}\cdot g_{1}\cdot\sigma_{0})$ is
the unit in $\Upsilon$, then $g_{n}\ldots g_{2}g_{1}=1$
and $(\hg_{1}\com\hg_{2}\com\ldots\com\hg_{n}\com\sigma_{n})
\ldots(\hg_{1}\com\hg_{2}\com\sigma_{2})\,
(\hg_{1}\com\sigma_{1})\,\sigma_{0}$ is homotopic to the
constant loop in $X$. In other words, the homomorphism
$\bar{f}$ is injective, hence an isomorphism.
\eop

\chapter{Hilsum-Skandalis Maps}
\label{chapHilSkaMap}
\startchapterskip

In Section \ref{secHaeStr} we saw that the Haefliger
$\GG$-structures on a space $X$ correspond exactly to the
isomorphism classes of principal $\GG$-bundles over $X$, for
any \'{e}tale groupoid $\GG$ and topological space $X$.
In particular, if $\GG$ is just a topological space $Y$,
the Haefliger $Y$-structures over $X$ are precisely the
continuous maps from $X$ to $Y$. In this Chapter we show,
among other things, that in fact any Haefliger $\GG$-structure
can be viewed as a map from $X$ to $\GG$.

In 1987 \cite{HilSka} M. Hilsum and G. Skandalis introduced
the notion of a $\eCe$-map between the
holonomy groupoids $\Hol(M,\cF)$ and $\Hol(M',\cF')$
associated to foliations $\cF$ and $\cF'$
on finite-dimensional $\eCe$-manifolds.
By replacing the holonomy groupoids with arbitrary topological
groupoids, one gets the notion of a Hilsum-Skandalis map between
topological groupoids \cite{Hae84,HilSka,Pra}.
Such a map between topological groupoids
$\HH$ and $\GG$ is an isomorphism class of principal
$\GG$-$\HH$-bibundles. The Hilsum-Skandalis maps arose
independently in topos theory \cite{Bun,Moe1}.

The Hilsum-Skandalis maps form a category $\cGpd$, with the
topological groupoids as objects. We will construct a
functor from the category $\Gpd$ of continuous functors between
topological groupoids to the category $\cGpd$, which
is the identity on objects and
sends essential equivalences to isomorphisms.

If $\GG$ is an \'{e}tale groupoid and $X$ a topological
space, the principal $\GG$-$X$-bibundles are exactly the principal
$\GG$-bundles over $X$. In particular, a Haefliger $\GG$-structure
on $X$ is a Hilsum-Skandalis
map from $X$ to $\GG$. We will show that just as a Haefliger
$\GG$-structure on $X$ can be seen as a generalized foliation on $X$,
a Hilsum-Skandalis map $\HH\ra\GG$ between topological groupoids can
be seen as a generalized foliation on the topological groupoid $\HH$.

\section{Bibundles}  \label{secBib}

\begin{dfn}  \label{IIdfn1}
Let $\GG$ and $\HH$ be topological groupoids. A
$\GG$-$\HH$-bibundle is a topological space $E$, equipped with a
left $\GG$-action $(p,\mu)$ and with a right $\HH$-action
$(\xW,\nu)$ such that the two actions commute with each other,
i.e.
\begin{enumerate}
\item [(i)]   $\xW(g\cdot e)=\xW(e)$,
\item [(ii)]  $p(e\cdot h)=p(e)$ and
\item [(iii)] $(g\cdot e)\cdot h=g\cdot(e\cdot h)$,
\end{enumerate}
for any $g\in\GG_{1}$, $e\in E$, $h\in\HH_{1}$ with
$\dom g=p(e)$ and $\xW(e)=\cod h$.
\end{dfn}
\Rem
We will denote such a $\GG$-$\HH$-bibundle by $(E,p,\xW,\mu,\nu)$
or $(E,p,\xW)$, or simply by $E$. A special kind of
$\GG$-$\HH$-bibundles was introduced in \cite{HilSka} for the case
where $\GG$ and $\HH$ are the holonomy groupoids of foliations.
In \cite{Moe4}, $\GG$-$\HH$-bibundles are referred to as
$\GG$-$\HH$-bispaces.
\vspace{4 mm}

Let $\GG$ and $\HH$ be topological groupoids and let $(E,p,\xW,\mu,\nu)$
be a $\GG$-$\HH$-bibundle. Observe that, in particular, $(E,p,\xW,\mu)$
is a $\GG$-bundle over $\HH_{0}$, called the {\em underlying}
$\GG$-bundle 
of the $\GG$-$\HH$-bibundle $E$.
We call $(E,p,\xW,\mu,\nu)$ {\em left transitive}
(or simply {\em transitive}) if the underlying
$\GG$-bundle $(E,p,\xW,\mu)$ is transitive.
We say that $(E,p,\xW,\mu,\nu)$ is {\em left principal} (or simply
{\em principal}) if the underlying
$\GG$-bundle $(E,p,\xW,\mu)$ is principal.

Let $\psi:\GG\ra\GG'$ and $\phi:\HH\ra\HH'$ be continuous functors
between topological groupoids, $(E,p,\xW)$ a $\GG$-$\HH$-bibundle and
$(E',p',\xW')$ a $\GG'$-$\HH'$-bibundle. A continuous map $\alpha:E\ra E'$ 
is called $\psi$-$\phi$-{\em equivariant} if
$\psi_{0}\com p=p'\com\alpha$, $\phi_{0}\com \xW=\xW'\com\alpha$ and
$$ \alpha(g\cdot e\cdot h)=\psi(g)\cdot\alpha(e)\cdot\phi(h)\;,$$
for all $g\in\GG_{1}$, $e\in E$, $h\in\HH_{1}$ with $\dom g=p(e)$,
$\xW(e)=\cod h$.

If $E$ and $E'$ are two $\GG$-$\HH$-bibundles, a continuous map
$\alpha:E\ra E'$ is called $\GG$-$\HH$-{\em equivariant} (or simply
{\em equivariant}) if it is $\id_{\!\GG}$-$\id_{\!\HH}$-equivariant.
Two $\GG$-$\HH$-bibundles are called {\em isomorphic} if there
exists an equivariant homeomorphism between them.
We will often denote the isomorphism class of a $\GG$-$\HH$-bibundle
$E$ again by $E$. Any equivariant map between principal
$\GG$-$\HH$-bibundles is a homeomorphism. An isomorphism class
of principal $\GG$-$\HH$-bibundles is called a {\em Hilsum-Skandalis
map} from $\HH$ to $\GG$.

If $\GG$ and $\HH$ are $\eCr$-groupoids, a $\GG$-$\HH$-bibundle
of class $\eCr$ is a $\GG$-$\HH$-bibundle $E$ such that $E$ is a
$\eCr$-manifold and both actions are of class $\eCr$. Such a bundle
is $\eCr$-principal if the underlying $\GG$-bundle over $\HH_{0}$
is $\eCr$-principal. Two $\GG$-$\HH$-bibundles of class $\eCr$
are {\em isomorphic} if there exists an equivariant
$\eCr$-diffeomorphism between them.
A {\em Hilsum-Skandalis $\eCr$-map} from
$\HH$ to $\GG$ is an isomorphism class of $\eCr$-principal
$\GG$-$\HH$-bibundles.

\begin{ex}  \label{IIex2}  \rm
If $\GG$ is a topological groupoid and $X$ a topological space,
the (principal) $\GG$-$X$-bibundles are exactly the (principal)
$\GG$-bundles over $X$. In particular, if $\GG$ is \'{e}tale, the
Haefliger $\GG$-structures on $X$ are in a natural bijective
correspondence with the Hilsum-Skandalis maps from $X$ to $\GG$
(Proposition \ref{Iprop7}).
If $X$ and $Y$ are topological spaces, then the Hilsum-Skandalis maps
from $X$ to $Y$ are precisely the continuous maps from $X$ to $Y$
(Example \ref{Iex6a}). If $M$ and $N$ are $\eCr$-manifolds,
the Hilsum-Skandalis maps of class $\eCr$ from $M$ to $N$
are precisely the $\eCr$-maps from $M$ to $N$.
\end{ex}

Let $\GG$ and $\HH$ be topological groupoids and $(E,p,\xW)$
a $\GG$-$\HH$-bibundle. Denote by $q_{E}:E\ra E/\HH$ the quotient
projection on the space of orbits $E/\HH$ for the action of $\HH$ on
$E$. Then $\xW$ induces a map
$$ \xW/\HH:E/\HH\lra |\HH|\;,$$
and $p$ factors through $q_{E}$ as 
$$ p=p/\HH\com q_{E}\;.$$
Since the left action of $\GG$ commutes with the action
of $\HH$, it induces a left action of $\GG$ on $E/\HH$ with respect
to $p/\HH$ and along the fibers of $\xW/\HH$. 
Therefore $(E/\HH,p/\HH,\xW/\HH)$ is a $\GG$-bundle over $|\HH|$,
called the {\em associated} $\GG$-bundle to the $\GG$-$\HH$-bibundle
$(E,p,\xW)$. The map $q_{E}$ is $\id_{\!\GG}$-$\rrr$-equivariant, where
$\rrr:\HH\ra |\HH|$ is the canonical continuous functor.
If $(E,p,\xW)$ is transitive, then $(E/\HH,p/\HH,\xW/\HH)$ is also
transitive. However, if $(E,p,\xW)$ is principal, the bundle
$(E/\HH,p/\HH,\xW/\HH)$ may not be principal.

Let $(E,p,\xW,\mu,\nu)$ be a $\GG$-$\HH$-bibundle and let
$(E',p',\xW',\mu',\nu')$ be an $\HH$-$\KK$-bibundle. Then we can define
a $\GG$-$\KK$-bibundle
$$ (E\otimes E',p\otimes p',\xW\otimes \xW',
   \mu\otimes\mu',\nu\otimes\nu')\;,$$
called the {\em tensor product} of $E$ and $E'$,
as follows: Define a right $\HH$-action on the space
$$ E\times_{\HH_{0}}E'=\{\,(e,e')\,|\,\xW(e)=p'(e')\,\}\subset
   E\times E' $$
with respect to the map $\xW\com\pr_{1}=p'\com\pr_{2}$ by
$$ (e,e')\cdot h=(e\cdot h,h^{-1}\cdot e')\;,$$
for any $e\in E$, $e'\in E'$, $h\in\HH_{1}$ with
$\xW(e)=p'(e')=\cod h$. Now we take
$$ E\otimes E'=(E\times_{\HH_{0}}E')/\HH\;.$$
We denote by $e\otimes e'\in E\otimes E'$ the equivalence class of
$(e,e')\in E\times_{\HH_{0}}E'$.
The maps $p\otimes p'$ and $\xW\otimes \xW'$ are given by
$$ (p\otimes p')(e\otimes e')=p(e) $$
and
$$ (\xW\otimes \xW')(e\otimes e')=\xW'(e')\;.$$
The left action of $\GG$ is given by
$$ g\cdot(e\otimes e')=g\cdot e\otimes e'\;,$$
and the right action of $\KK$ is given by
$$ (e\otimes e')\cdot h=e\otimes e'\cdot h\;.$$
It is obvious that $(E\otimes E',p\otimes p',\xW\otimes \xW')$ is indeed
a $\GG$-$\KK$-bibundle. But we claim:

\begin{prop}  \label{IIprop3}
Let $E$ be a $\GG$-$\HH$-bibundle and
$E'$ a $\HH$-$\KK$-bibundle. If $E$ and $E'$ are principal
(respectively transitive), then the $\GG$-$\KK$-bibundle
$E\otimes E'$ is also principal (respectively transitive).
\end{prop}
\Proof
Let $E=(E,p,\xW,\mu,\nu)$ and $E'=(E',p',\xW',\mu',\nu')$.
Denote by $q:E\times_{\HH_{0}}E'\ra E\otimes E'$ the open quotient
projection. Since $\xW$ is an open surjection, the projection
$\pr_{2}:E\times_{\HH_{0}}E'\ra E'$ is an open surjection.
Since $q$ is open and $\xW\otimes \xW'$ is the factorization of
$\xW'\com\pr_{2}$ though $q$, this implies that $\xW\otimes \xW'$
is an open surjection.

Assume that the bibundles $E$ and $E'$ are principal (in the
transitive case the proof is analogous). The map
$(\mu,\pr_{2}):\GG_{1}\times_{\GG_{0}}E\ra E\times_{\HH_{0}}E$
is a homeomorphism, with the inverse of the form $(\theta,\pr_{2})$.
Also, the map
$(\mu',\pr_{2}):\HH_{1}\times_{\HH_{0}}E'\ra E'\times_{\KK_{0}}E'$ is
a homeomorphism, with the inverse of the form $(\theta',\pr_{2})$.
We have to prove that the map
$$ (\mu\otimes\mu',\pr_{2}):\GG_{1}\times_{\GG_{0}}E\otimes E'\lra
   E\otimes E'\times_{\KK_{0}}E\otimes E' $$
is a homeomorphism. But we can define the inverse $\kappa$ 
of $(\mu\otimes\mu',\pr_{2})$ by
$$ \kappa(e\otimes e',e_{1}\otimes e'_{1})=
   (\theta(e\cdot\theta'(e',e'_{1}),e_{1}),e_{1}\otimes e'_{1})\;.$$
It is easy to verify that this is indeed a well-defined inverse.
\eop

\begin{dfn}  \label{IIdfn4}
The category $\cGpd$ of Hilsum-Skandalis maps is given by:
\begin{enumerate}
\item [(i)]   the objects of $\cGpd$ are the topological groupoids,
\item [(ii)]  the morphisms in $\cGpd(\HH,\GG)$ are the
              Hilsum-Skandalis maps from $\HH$ to $\GG$, and
\item [(iii)] the composition in $\cGpd$ is induced by the tensor
              product of principal bibundles.
\end{enumerate}
\end{dfn}
\Rem
Let us emphasize again that in notation we will not make any
distinction between the principal bibundles and their isomorphism
classes, i.e. the associated Hilsum-Skandalis maps. In particular,
we will denote the composition in $\cGpd$ induced by the tensor
product $\otimes$ of bibundles again by $\otimes$.

It is easy to verify that $\cGpd$ is indeed a category. For example,
the tensor product  of bibundles is associative up to
a canonical isomorphism, and
the identity in $\cGpd$ is (the isomorphism class of) the principal
$\GG$-$\GG$-bibundle
$$ (\GG_{1},\cod,\dom) $$
with the obvious actions given by the composition in $\GG$.

We will denote by $\cGpde$ the full subcategory of $\cGpd$ with
objects all \'{e}tale groupoids. The category $\Top$ of
topological spaces is a full subcategory of $\cGpde$ (Example
\ref{IIex2}).

There is a close relation between the category $\cGpd$
and the category of Grothendieck toposes and geometric morphisms between
them \cite{Moe4}.
\vspace{4 mm}

Let $\phi:\HH\ra\GG$ be a continuous functor between topological
groupoids. We define
$$ \angs{\phi}=\GG_{1}\times_{\GG_{0}}\HH_{0}=\{\,(g,b)\,|\,
   \dom g=\phi_{0}(b)\,\}\subset\GG_{1}\times\HH_{0}\;.$$
There are natural actions of $\GG$ and $\HH$ on
$\angs{\phi}$ such that
$$ (\angs{\phi},\cod\com\pr_{1},\pr_{2}) $$
is a principal $\GG$-$\HH$-bibundle. Indeed, we have a left
$\GG$-action $\mu$ on $\angs{\phi}$ with respect to $\cod\com\pr_{1}$
given by
$$ g'\cdot (g,b)=(g'\com g,b) $$
and a right $\HH$-action $\nu$ on $\angs{\phi}$ with respect to 
$\pr_{2}$ given by
$$ (g,b)\cdot h=(g\cdot\phi(h),\dom h)\;.$$
With this, $\angs{\phi}$ is clearly a $\GG$-$\HH$-bibundle.
Observe that $pr_{2}$
is an open surjection since the domain map of $\GG$ is an open
surjection. Moreover, the map
$$ (\mu,\pr_{2}):\GG_{1}\times_{\GG_{0}}\angs{\phi}\lra
   \angs{\phi}\times_{\HH_{0}}\angs{\phi} $$
has a continuous inverse $\kappa$ given by
$\kappa((g_{1},b),(g,b))=(g_{1}\com g^{-1},(g,b))$.
Thus $\angs{\phi}$ is indeed a principal $\GG$-$\HH$-bibundle.

If $\psi:\KK\ra\HH$ is another continuous functor between open
groupoids, there is a continuous $\GG$-$\KK$-equivariant map
$$ \alpha:\angs{\phi\com\psi}\lra\angs{\phi}\otimes\angs{\psi} $$
given by $\alpha(g,c)=(g,\psi_{0}(c))\otimes(1_{\psi_{0}(c)},c)$.
As any equivariant map between principal bibundles, the map $\alpha$
is a homeomorphism. Furthermore, there is an isomorphism
$\angs{\id_{\GG}}\cong (\GG_{1},\cod,\dom)$.
This proves that
$$ \angs{\oo}:\Gpd\lra\cGpd $$
is a functor. On objects, the functor $\angs{\oo}$ is the identity.

\begin{prop}  \label{IIprop5}
Let $(E,p,\xW)$ be a principal $\GG$-$\HH$-bibundle. Then
$E\cong\angs{\phi}$ for some continuous functor $\phi:\HH\ra\GG$ if
and only if the map $\xW$ has a continuous section.
\end{prop}
\Proof
If $s$ is a continuous section of $\xW$, i.e. a
right inverse of $\xW$, then define
$\phi$ with $\phi_{0}=p\com s$ such that
$$ \phi(h)\cdot s(\dom h)=s(\cod h)\cdot h\;.$$
Since $E$ is principal, this equation defines $\phi$
uniquely.
\eop

\begin{prop}  \label{IIprop6}
Let $\phi:\HH\ra\GG$ be a continuous functor between topological
groupoids. Then $\angs{\phi}$ is an isomorphism in $\cGpd$
if and only if $\phi$ is an essential equivalence.
\end{prop}
\Proof
(a) First assume that $\phi$ is an essential equivalence.
Write the structure maps of the associated principal bundle of
$\phi$ by $(\angs{\phi},p,\xW,\mu,\nu)$.
The left $\GG$-action $(p,\mu)$ induces the corresponding right
$\GG$-action $(p,\bar{\mu})$ on $\angs{\phi}$ by
$$ \bar{\mu}(e,g)=\mu(g^{-1},e)\;.$$
Similarly, the right $\HH$-action $(\xW,\nu)$ induces the
corresponding left $\HH$-action $(\xW,\bar{\nu})$ on $\angs{\phi}$.
Clearly, $(\angs{\phi},\xW,p,\bar{\nu},\bar{\mu})$ is a
$\HH$-$\GG$-bibundle. We will
show that this bundle is principal and provides the inverse for
$(\angs{\phi},p,\xW,\mu,\nu)$ in $\cGpd$.

Since $\phi$ is an essential equivalence, the map $p$ is an open
surjection and $(\phi,\dom,\cod):\HH_{1}\ra\GG_{1}
\times_{(\GG_{0}\times\GG_{0})}(\HH_{0}\times\HH_{0})$ is a
homeomorphism. Denote by $\vartheta$ the inverse of
$(\phi,\dom,\cod)$. To prove that
$(\angs{\phi},\xW,p,\bar{\nu},\bar{\mu})$
is principal we still have to prove that the map
$$ (\bar{\nu},\pr_{2}):\HH_{1}\times_{\HH_{0}}\angs{\phi}\lra
   \angs{\phi}\times_{\GG_{0}}\angs{\phi} $$
is a homeomorphism.
In fact, we can define the inverse $\kappa$ of $(\bar{\nu},\pr_{2})$
explicitly
(by using $\angs{\phi}=\GG_{1}\times_{\GG_{0}}\HH_{0}$) by
$$ \kappa((g',b'),(g,b))=
   (\vartheta(g'^{-1}\!\!\com g,(b,b')),(g,b))\;.$$
Hence the bibundle $(\angs{\phi},\xW,p,\bar{\nu},\bar{\mu})$ is principal,
and it is straightforward to check that it provides the inverse for
$(\angs{\phi},p,\xW,\mu,\nu)$ in $\cGpd$.

(b) Assume that $\angs{\phi}$ is an isomorphism in $\cGpd$. Hence
there exists a principal $\HH$-$\GG$-bibundle $(E,p,\xW)$ and
equivariant homeomorphisms $\alpha:\HH_{1}\ra E\otimes\angs{\phi}$
and $\beta:\angs{\phi}\otimes E\ra\GG_{1}$.

First observe that
$E\otimes\angs{\phi}\cong E\times_{\GG_{0}}\HH_{0}$.
By composing the homeomorphism $\alpha$ with the first
projection one gets a continuous map
$$ \delta:\HH_{1}\lra E $$
which satisfies $\delta(h\com h')=h\cdot\delta(h')=
\delta(h)\cdot\phi(h')$. We now define a continuous map
$u:\angs{\phi}\ra E$ (recall that 
$\angs{\phi}=\GG_{1}\times_{\GG_{0}}\HH_{0}$) by
$$ u(g,b)=\delta(1_{b})\cdot g^{-1}\;.$$
Clearly we have
$u(g'\cdot (g,b)\cdot h)=h^{-1}\cdot u(g,b)\cdot g'^{-1}$.

Next we define a continuous map
$v:E\ra\angs{\phi}$ by
$$ v(e)=(\beta((1_{\phi(p(e))},p(e))\otimes e)^{-1},p(e))\;.$$
Since $\beta$ is equivariant,
we have $v(h\cdot e\cdot g)=g^{-1}\cdot v(e)\cdot h^{-1}$.

Now $u\com v$ is an equivariant map of the principal bibundle
$E$ into itself, hence a homeomorphism. Similarly, $v\com u$ is
a homeomorphism. Therefore both $u$ and $v$ are homeomorphisms.
Since $\xW$ is an open surjection and $\xW\com u=\cod\com\pr_{1}$,
this implies that
$\cod\com\pr_{1}:\GG_{1}\times_{\GG_{0}}\HH_{0}\ra\GG_{0}$
is an open surjection. Finally, the composition of the
homeomorphisms $\inv:\HH_{1}\ra\HH_{1}$,
$\alpha:\HH_{1}\ra E\times_{\GG_{0}}\HH_{0}$ and
$$ u^{-1}\times\id_{\!\HH_{0}}:E\times_{\GG_{0}}\HH_{0}\lra
   (\GG_{1}\times_{\GG_{0}}\HH_{0})\times_{\GG_{0}}\HH_{0} $$
is exactly of the form $(\phi,\dom,\cod)$. Thus $\phi$ is indeed
an essential equivalence.
\eop

\begin{cor}  \label{IIcor7}
Two topological groupoids are Morita equivalent if and only if
they are isomorphic in the category $\cGpd$.
\end{cor}
\Proof
One direction is clear. For the other, assume that $(E,p,\xW)$ is
a principal $\GG$-$\HH$-bibundle which is an isomorphism in
$\cGpd$. Let $\KK=\GG\times\HH$, i.e. the product of $\GG$
and $\HH$ given by $\KK_{0}=\GG_{0}\times\HH_{0}$ and
$\KK_{1}=\GG_{1}\times\HH_{1}$ and with the obvious structure
maps. This groupoid acts on $E$ with respect to $(p,\xW)$ from
the right by
$$ e\cdot (g,h)=g^{-1}\cdot e\cdot h\;.$$
Let $\KK(E)$ be the topological groupoid associated
to this action on $E$. Observe that there is a continuous functor
$$ \phi:\KK(E)\lra H\;,$$
given by $\phi_{0}=\xW$ and $\phi_{1}(e,(g,h))=h$, which is
clearly an essential equivalence. Moreover, there is a continuous
functor
$$ \psi:\KK(E)\lra G $$
given by $\psi_{0}=p$ and $\psi_{1}(e,(g,h))=g$. It is easy to
check that
$$ E\otimes\angs{\phi}=\angs{\psi}\;.$$
Now Proposition \ref{IIprop6} implies first that $\angs{\psi}$
is an isomorphism, and second that $\psi$ is an essential
equivalence.
\eop

\Not
Let $(E,p,\xW)$ be a $\GG$-$\HH$-bibundle. If
$\phi:\HH'\ra\HH$ is a continuous functor between topological
groupoids, we denote the structure maps of the
tensor product $\phi^{\ast}E=E\otimes\angs{\phi}$ as
$$ (\phi^{\ast}E,\phi^{\ast}p,\phi^{\ast}\xW)\;.$$
Observe that they fit into the pull-back
$$\CD
  \phi^{\ast}E \cdr{\phi^{E}}{} E \cdr{p}{} \GG_{0} \\
  \cdd{\phi^{\ast}\xW}{}       \cdd{}{\xW}          \\
      \HH'_{0} \cdr{\phi_{0}}{} \HH_{0} 
  \endCD$$
with $\phi^{\ast}p=p\com\phi^{E}$. The map
$\phi^{E}$, which is given by $\phi^{E}(e\otimes(h,b'))=e\cdot h$,
is $\id_{\!\GG}$-$\phi$-equivariant.

Let $\psi:\GG\ra\GG'$ be a continuous functor. We will write the
structure maps of $\angs{\psi}\otimes E$ as
$$ (\angs{\psi}\otimes E,\angs{\psi}\otimes p,
   \angs{\psi}\otimes \xW)\;.$$
Note that $\angs{\psi}\times_{\GG_{0}} E\cong
\GG'_{1}\times_{\GG'_{0}} E$, so we will denote
an element $(g',p(e))\otimes e$ of $\angs{\psi}\otimes E$ simply
by $g'\otimes e$.

Observe that the map $\psi_{E}:E\ra\angs{\psi}\otimes E$, given by
$$ \psi_{E}(e)=1_{\psi_{0}(p(e))}\otimes e\;,$$
is $\psi$-$\id_{\!\HH}$-equivariant.
\vspace{4 mm}

\Rem
Let $\GG$, $\HH$ and $\KK$ be \'{e}tale $\eCr$-groupoids,
let $E$ be a $\eCr$-principal $\GG$-$\HH$-bibundle, and
let $E'$ be a $\eCr$-principal $\HH$-$\KK$-bibundle.
The quotient map
$$ q:E\times_{\HH_{0}}E'\lra E\otimes E' $$
is a local homeomorphism,
and since this is a quotient for an action of
class $\eCr$, there is a natural structure of a $\eCr$-manifold
on $E\otimes E'$ such that the map $q$ is a $\eCr$-diffeomorphism.
With this, $E\otimes E'$ is clearly a $\GG$-$\KK$-bibundle of
class $\eCr$. The same argument as in the proof of
Proposition \ref{IIprop3} now implies
that $E\otimes E'$ is in fact a $\eCr$-principal
$\GG$-$\KK$-bibundle. Therefore a composition of
Hilsum-Skandalis $\eCr$-maps is a Hilsum-Skandalis $\eCr$-map.

Thus we can define the category $\cGpder$ of \'{e}tale $\eCr$-groupoids
and Hilsum-Skandalis $\eCr$-maps between them, with the composition
given by the tensor product.

If $\phi:\HH\ra\GG$ is a $\eCr$-functor, then $\angs{\phi}$ is
clearly a $\eCr$-principal $\GG$-$\HH$-bibundle. Thus we have
a functor
$$ \angs{\oo}:\Gpder\lra\cGpder\;.$$
The same arguments as in the proof of Proposition \ref{IIprop6}
and Corollary \ref{IIcor7}
now show that $\angs{\phi}$ is an isomorphism in $\cGpder$
if and only if $\phi$ is a $\eCr$-essential equivalence,
and that $\GG$ and $\HH$ are $\eCr$-Morita equivalent if and
only if they are isomorphic in $\cGpder$. The category of
$\eCr$-manifolds and $\eCr$-maps between them is a full
subcategory of $\cGpder$.

\begin{ex}  \label{IIex8}  \rm
(1) Let $X$ be a topological space and $\GG$ an \'{e}tale groupoid.
Proposition \ref{Iprop7} gives a natural
bijective correspondence
$$ H^{1}(X,\GG)\cong\cGpd(X,\GG)\;.$$
If $f:X'\ra X$ is a continuous function, one gets a map
$$ f^{\ast}:H^{1}(X,\GG)\lra H^{1}(X',\GG) $$
by composition with $f$ in $\cGpd$. If $\ccc\in H^{1}(X,\GG)$
is represented by $(c_{ij})\in Z^{1}(\cU,\GG)$ for some open
cover $\cU$ of $X$, then the $\GG$-cocycle
$$ (c_{ij}\com f|_{f^{-1}(U_{i})\cap f^{-1}(U_{j})}) $$
represents $f^{\ast}(\ccc)$.
Further, if $\psi:\GG\ra\GG'$ is a continuous functor, there is a map
$$ \psi_{\ast}:H^{1}(X,\GG)\lra H^{1}(X,\GG') $$
given by the composition with $\angs{\psi}$ is $\cGpd$. The
Haefliger $\GG'$-structure $\psi_{\ast}(\ccc)$ is now represented
by $(\psi\com c_{ij})$.

(2) Let $M$ be a $\eCe$-manifold of dimension $n$
and $\cF$ a foliation on $M$ of codimension $q$. 
Let $(\Sigma(\cF),p,\xW)$ be the $\eCe$-principal
$\GmqCe$-bundle
over $M$ associated to $\cF$ (Example \ref{Iex9} (1)).
This bundle gives a Hilsum-Skandalis $\eCe$-map
$$ \Sigma(\cF)\in\cGpde^{1}(M,\GmqCe)\;.$$
Now let $M'$ be another $\eCe$-manifold of dimension $n'$
and $f:M'\ra M$ a $\eCe$-map.
The $\eCe$-principal $\GmqCe$-bundle 
$$ (f^{\ast}\Sigma(\cF),f^{\ast}p,f^{\ast}\xW) $$
over $M'$ does not represent
a classical foliation on $N$, since the map $f^{\ast}p$ may
not be a submersion. We say that $f^{\ast}\Sigma(\cF)$ represents a
{\em foliation with modular singularities} if
$f^{\ast}p$ has locally path-connected fibers. 
It is easy to see that this is the case
if and only if for any foliation chart $\varphi:U\ra\RRR^{n}$
for $\cF$, the inverse images along $f$ of the
plaques of $\cF$ in $U$ are locally path-connected.
The Haefliger-Reeb-Ehresmann stability theorem for Haefliger
structures gives a generalization of the Reeb stability theorem
to the foliations with modular singularities (see Chapter
\ref{chapSta}).

(3) Let $M$ be a $\eCe$-manifold of dimension $n$, $\cF$ a
foliation on $M$ of codimension $q$, and let 
$$ \nu:M\times G\lra M $$
be a $\eCe$-action of a discrete group $G$ on $M$ which preserves
$\cF$. More precisely, for any $g\in G$, the $\eCe$-diffeomorphism
$\hg=\nu(\oo,g):M\ra M$ maps leaves into leaves. We will refer to
such a foliation $\cF$ as $G$-equivariant, or simply equivariant.

Let $(\varphi_{i}:U_{i}\ra\RRR^{n-q}\times\RRR^{q})_{i\in I}$
be the maximal atlas for $\cF$.
If $i\in I$ and $g\in G$, then
$$ \varphi_{i}\com\hg^{-1}|_{U_{i}\cdot g}:U_{i}\cdot g\lra
   \RRR^{n-q}\times\RRR^{q} $$
is also a chart for $\cF$ and hence equal to the chart
$\varphi_{i\cdot g}$ for a uniquely determined $i\cdot g\in I$.

Now take $s_{i}=\pr_{2}\com\varphi_{i}$, and let $c=(c_{ij})$ be the
corresponding $\GmqCe$-cocycle on $(U_{i})$, as described
in Section \ref{secFol}. Let $(\Sigma(c),p,\xW)$ be the
$\eCe$-principal $\GmqCe$-bundle over $M$ associated to this
cocycle. Recall that an element of $\Sigma(c)$ is of the form
$[\gm,x,i]$, where $\gm\in\GmqCe$, $i\in I$ and $x\in U_{i}$
with $\dom\gm=s_{i}(x)$.

The group $G$ acts on $\Sigma(c)$ by
$$ [\gm,x,i]\cdot g=[\gm,x\cdot g, i\cdot g]\;.$$
This action clearly induces an action of the \'{e}tale groupoid $G(M)$
on $\Sigma(c)$, and $\Sigma(c)$ becomes a $\eCe$-principal
$\GmqCe$-$G(M)$-bibundle. Thus we obtain a
Hilsum-Skandalis $\eCe$-map
$$ \Sigma(\cF,\nu)\in\cGpde^{1}(G(M),\GmqCe)\;.$$

(4) Let $M$ be a $\eCe$-manifold of dimension $n$,
and let $r,s,t$ be natural numbers such that $r+s+t=n$.
A {\em nested foliation} on $M$ of type $(r,s,t)$ is a (maximal)
atlas $(\varphi_{i}:U_{i}\ra\RRR^{n})_{i\in I}$ on $M$ such that
the change of coordinates diffeomorphisms
$\varphi_{ij}=\varphi_{i}\com
 \varphi_{j}^{-1}|_{\varphi_{j}(U_{i}\cap\, U_{j})}$
are locally of the form
$$ \varphi_{ij}(x,y,z)=
   (\varphi^{(1)}_{ij}(x,y,z),\varphi^{(2)}_{ij}(y,z),
   \varphi^{(3)}_{ij}(z)) $$
with respect to the decomposition 
$\RRR^{n}=\RRR^{r}\times\RRR^{s}\times\RRR^{t}$.
In particular, a nested foliation determines two foliations
$\cF_{1}$ and $\cF_{2}$ on $M$, the first of codimension $t$ and
the second of codimension $s+t$. Moreover, $\cF_{2}$ restricts 
on each leaf of $\cF_{1}$ to a foliation of codimension $s$. 

We can assume without loss of generality that each $\varphi_{i}$
is surjective. Put $N_{1}=\coprod_{i\in I}\RRR^{t}$ and
$N_{2}=\coprod_{i\in I}\RRR^{s+t}$, and define
$T_{1}:N_{1}\ra M$ and $T_{2}:N_{2}\ra M$ by
$T_{1}(z,i)=\varphi_{i}^{-1}(0,0,z)$ and
$T_{2}((y,z),i)=\varphi_{i}^{-1}(0,y,z)$. Clearly
$T_{1}$ is a complete transversal for $\cF_{1}$, and $T_{2}$ is
a complete transversal for $\cF_{2}$.

We have the canonical projection $\psi_{0}:N_{2}\ra N_{1}$.
It is easy to check that with $\psi_{0}$ one can project each
germ of a diffeomorphism in $\Hol_{T_{2}}(M,\cF_{2})$
to a germ of a diffeomorphism in $\Hol_{T_{1}}(M,\cF_{1})$,
obtaining a $\eCe$-functor 
$$ \psi:\Hol_{T_{2}}(M,\cF_{2})\lra\Hol_{T_{1}}(M,\cF_{1})\;,$$
which is a submersion.
Hence we get a Hilsum-Skandalis $\eCe$-map
$$ \angs{\psi}\in\cGpde^{1}(\Hol_{T_{2}}(M,\cF_{2}),
   \Hol_{T_{1}}(M,\cF_{1}))\;.$$
Moreover, if $E_{1}$ is the $\eCe$-principal
$\Hol_{T_{1}}(M,\cF_{1})$-bundle over $M$ associated to
$\cF_{1}$ and $E_{2}$ the $\eCe$-principal
$\Hol_{T_{2}}(M,\cF_{2})$-bundle over $M$ associated to $\cF_{2}$,
(Example \ref{Iex9} (2)), then
$$ \angs{\psi}\otimes E_{2}=E_{1}\;.$$

(5) The category $\Grp$ of (discrete) groups is a full subcategory of
$\Gpde$. If $G$ and $H$ are groups and $f,f':H\ra G$ homomorphisms,
then $\angs{f}=\angs{f'}$ if and only if $f$ and $f'$ differ by the
conjugation by an element of $G$.
\end{ex}

\section{Leaves and Holonomy of Transitive Bibundles}
\label{secLeaHolTranBib}

Let $\GG$ and $\HH$ be topological groupoids and let
$(E,p,\xW)$ be a transitive $\GG$-$\HH$-bibundle. Let $a\in\GG_{0}$.
The {\em fiber} of $E$ over $a$ is the space
$E_{a}=p^{-1}(a)\subset E$, i.e. the fiber of the map $p$ over $a$.
The fiber $E_{a}$ is clearly $\HH$-invariant and also
$\GG(a,a)$-invariant. So the right action of $\HH$ restricts
to $E_{a}$, and we have the quotient projection
$$ q_{E_{a}}:E_{a}\lra E_{a}/\HH $$
on the orbit space $E_{a}/\HH$.

Let $\tL$ be an $\HH$-connected component
of $E_{a}$, i.e. the inverse image
$q_{E_{a}}^{-1}(Z)$ of a connected component
$Z$ of $E_{a}/\HH$. Again, $\tL$ is $\HH$-invariant, and also
$\xW(\tL)$ is an $\HH$-invariant subspace of $\HH_{0}$. Let $\cH(\tL)$
be the subgroup of those elements of $\GG(a,a)$ which leave
$\tL$ invariant, i.e.
$$ \cH(\tL)=\{\,g\,|\,g\cdot\tL\subset\tL\,\}
   \subset\GG(a,a)\;.$$
The group $\cH(\tL)$ is called the {\em holonomy group}
of the $\HH$-connected component $\tL$
of the fiber of $E$ over $a$.
Since the action of $\GG(a,a)$ on $E_{a}$ induces an action
on $E_{a}/\HH$, it follows that an element $g\in\GG(a,a)$ belongs to
$\cH(\tL)$ if and only if $g\cdot e\in\tL$ for some element $e\in\tL$.
Since $E$ is transitive, this implies that $\cH(\tL)$ acts on
$\tL$ transitively along the fibers of
$$ \xW|_{\tL}:\tL\lra \xW(\tL)\subset\HH_{0}\;.$$
Therefore $\xW|_{\tL}$ factors through the quotient projection
$\tL\ra\tL/\cH(\tL)$ as a continuous bijection
$$ i(\tL):\tL/\cH(\tL)\lra \xW(\tL)\subset\HH_{0}\;.$$
In particular, $i(\tL)$ induces a new topology on $\xW(\tL)$,
called the {\em leaf topology}, which is finer that the one inherited
from the space $\HH_{0}$. 

The action of $\HH$ on $\tL$ induces an action on $\tL/\cH(\tL)$
and $i(\tL)$ is equivariant under this action.
The set $\xW(\tL)$ with the leaf topology and
with the inherited action of $\HH$ from $\HH_{0}$
is thus an $\HH$-connected right $\HH$-space,
isomorphic to $\tL/\cH(\tL)$.

\begin{dfn}  \label{IIdfn9}
Let $\GG$ and $\HH$ be topological groupoids and let
$(E,p,\xW)$ be a transitive $\GG$-$\HH$-bibundle.
A leaf $L$ of $E$ associated to an
$\HH$-connected component $\tL$ of a fiber of $E$ is the
right $\HH$-space $\xW(\tL)$ with the corresponding leaf topology
and with the inherited action of $\HH$ from $\HH_{0}$.
A leaf $L$ is embedded if the leaf topology of $L$
coincide with the topology inherited from the space $\HH_{0}$.
\end{dfn}

\begin{lem}  \label{IIlem10}
Let $\GG$ and $\HH$ be topological groupoids and let $(E,p,\xW)$
be a transitive $\GG$-$\HH$-bibundle.
Let $\tL$ be an $\HH$-connected component
of a fiber $E_{a}$ and $\tL'$ an $\HH$-connected component
of a fiber $E_{a'}$, such that $\xW(\tL)\cap \xW(\tL')\neq\emptyset$.
Then there exists an element $g\in\GG(a,a')$ such that
$g\cdot\tL=\tL'$ and $g\com\cH(\tL)\com g^{-1}=\cH(\tL')$.
In particular, $\xW(\tL)=\xW(\tL')$ and the
leaf topologies of $\xW(\tL)$ and $\xW(\tL')$ coincide.
\end{lem}
\Proof
Let $b\in \xW(\tL)\cap \xW(\tL')$, and choose $e\in\tL$ and $e'\in\tL'$
such that $\xW(e)=\xW(e')=b$. Since $E$ is transitive,
there exists an element $g\in\GG(a,a')$
with $g\cdot e=e'$. Since the action of $g$ maps $E_{a}$ into
$E_{a'}$ and the $\HH$-connected components of $E_{a}$ into the
$\HH$-connected components of
$E_{a'}$, it follows that $g\cdot\tL\subset\tL'$.
Symmetrically, $\tL\supset g^{-1}\cdot\tL'$,
and hence $g\cdot\tL=\tL'$.

Let $g_{1}\in\cH(\tL)$, i.e. $g_{1}\cdot\tL\subset\tL$.
It follows that
$$ (g\com g_{1}\com g^{-1})\cdot\tL'=
   (g\com g_{1}\com g^{-1})\cdot(g\cdot\tL)=
   (g\com g_{1})\cdot\tL\subset g\cdot\tL=\tL'\;, $$
so $g\com\cH(\tL)\com g^{-1}\subset\cH(\tL')$. Again by symmetry
we get $\cH(\tL)\supset g^{-1}\com\cH(\tL')\com g$, and hence
$g\com\cH(\tL)\com g^{-1}=\cH(\tL')$. Thus the $\HH$-equivariant
homeomorphism $(g\cdot\oo):\tL\ra\tL'$ induces an $\HH$-equivariant
homeomorphism
$$ \alpha:\tL/\cH(\tL)\lra\tL'/\cH(\tL') $$
satisfying $i(\tL)=i(\tL')\com\alpha$.
\eop
\Rem
Lemma \ref{IIlem10} and surjectivity of $\xW$
imply that for a given $b\in\HH_{0}$
there exists a unique leaf $L$ of $E$ such that $b\in L$,
with a unique leaf topology.
The leaves of $E$ thus determine a partition of $\HH_{0}$ into 
$\HH$-invariant subsets. A leaf $L$ of $E$ may be associated to
different $\HH$-connected components of different fibers,
but the holonomy groups of these $\HH$-connected components are
isomorphic. This partition clearly depends only on the
isomorphism class of a transitive $\GG$-$\HH$-bibundle.

In this way, a transitive $\GG$-$\HH$-bibundle can be seen
as a generalized kind of foliation on the topological
groupoid $\HH$.

\begin{dfn}  \label{IIdfn11}
Let $\GG$ and $\HH$ be topological groupoids,
and let $(E,p,\xW)$ be a transitive $\GG$-$\HH$-bibundle.
The holonomy group $\cH(L)$ of a leaf $L$ of $E$
is (up to a conjugation) the holonomy group $\cH(\tL)$
of an $\HH$-connected component $\tL$ of a fiber of $E$ such
that $L$ is associated to $\tL$.
\end{dfn}

\begin{ex}  \label{IIex12}  \rm
(1) Let $p:X\ra Y$ be a continuous map between topological
spaces. Viewing $p$ as a principal $Y$-$X$-bibundle,
the leaves of $p$ are the connected
components of the fibers of $p$, and they are all embedded,
with trivial holonomy groups.

(2) Let $\GG$ be an \'{e}tale groupoid and $X$ a topological space.
Let $\ccc$ be a Haefliger $\GG$-structure on $X$, and
let $(E,p,\xW)$ be the corresponding (isomorphism class of a)
principal $\GG$-bundle over $X$. Note that, according to our
definition, the fibers of $E$ (which is a $\GG$-$X$-bibundle)
are the fibers of $p$ and not the fibers of $\xW$.
Now the definition of leaves and holonomy of
$E$ above exactly matches the definition of leaves and holonomy
of the $\GG$-structure $\ccc$ in \cite{Hae}.
In particular, in the special case of a foliation $\cF$ on a
finite-dimensional $\eCe$-manifold $M$ (Example \ref{Iex9})
this definition gives the classical leaves and holonomy of $\cF$.
\end{ex}

\begin{prop}  \label{IIprop13}
Let $\GG$ and $\HH$ be topological groupoids and $(E,p,\xW)$ a
transitive $\GG$-$\HH$-bibundle. If $L$ is a leaf of $E$,
then $L/\HH$ is a leaf of the associated $\GG$-bundle $E/\HH$
over $|\HH|$, and $\cH(L)\cong\cH(L/\HH)$.
\end{prop}
\Proof
Directly by definition, $q_{E}:E\ra E/\HH$ maps the fibers of
$E$ into the fibers of $E/\HH$ and the $\HH$-connected components
of the fibers of $E$ into the connected components of the fibers
of $E/\HH$. Since $q_{E}$ is also $\GG$-equivariant, the holonomy
groups are preserved as well.
\eop

\begin{prop}  \label{IIprop14}
Let $\psi:\GG\ra\GG'$ and $\phi:\HH\ra\HH'$ be continuous functors
between topological groupoids. Let $(E,p,\xW)$ be a transitive  
$\GG$-$\HH$-bibundle, $(E',p',\xW')$ a transitive $\GG'$-$\HH'$-bibundle,
and $\alpha:E\ra E'$ a continuous $\psi$-$\phi$-equivariant map.
If $\tL$ is an $\HH$-connected component of a fiber of $E$,
then $\alpha(\tL)$ lies in an $\HH'$-connected component $\tL'$ of a
fiber of $E'$, and
$$ \psi(\cH(\tL))\subset\cH(\tL')\;.$$ 
In particular, if $L$ is a leaf of $E$, then $\phi_{0}(L)$ lies in
a leaf $L'$ of $E'$ and $\phi_{0}|_{L}:L\ra L'$ is continuous with
respect to the leaf topologies.
\end{prop}
\Proof
Let $a\in\GG_{0}$. Since $\psi_{0}\com p=p'\com\alpha$, we have
$\alpha(E_{a})\subset E'_{\psi_{0}(a)}$. Since $\alpha$ is equivariant,
it induces a map $E_{a}/\HH\ra E'_{\psi_{0}(a)}/\HH'$ and hence maps an
$\HH$-connected component $\tL$ of $E_{a}$ into an
$\HH'$-connected component $\tL'$ of $E'_{\psi_{0}(a)}$. 
Now if $g\in\cH(\tL)$, we have $g\cdot\tL\subset\tL$ and thus
$$ \psi(g)\cdot\alpha(\tL)=\alpha(g\cdot\tL)\subset\alpha(\tL)\;.$$
This implies that $\psi(g)\cdot\tL'\subset\tL'$, 
and therefore $\psi(\cH(\tL))\subset\cH(\tL')$.
Thus $\alpha$ induces a continuous equivariant map
$$ \bar{\alpha}:\tL/\cH(\tL)\lra\tL'/\cH(\tL')\;,$$
which clearly satisfies $i(\tL')\com\bar{\alpha}=\phi_{0}\com i(\tL)$.
\eop

\section{Holonomy of Principal Bibundles}  \label{secHolPriBib}

Let $\GG$ be an \'{e}tale groupoid and $\HH$ a topological
groupoid, and assume that $E$ is a principal
$\GG$-$\HH$-bibundle with locally path-connected fibers
$E_{a}$, $a\in\GG_{0}$.
In this section we demonstrate that in this case
the holonomy group of a leaf $L$ of $E$ is the image of a
homomorphism from the fundamental group of the groupoid
$\HH(L)$. The analogy with foliations on manifolds is
hence complete.

\begin{prop}  \label{IIprop15}
Let $\GG$ be an \'{e}tale groupoid, $\HH$ a topological groupoid
and $(E,p,\xW)$ a principal $\GG$-$\HH$-bibundle.
If $\tL$ is an $\HH$-connected component of a fiber of $E$,
then the holonomy group $\cH(\tL)$ is discrete, and acts properly
discontinuously on $\tL$.
\end{prop}
\Proof
Let $a\in\GG_{0}$ such that $\tL\subset E_{a}$. Since $\GG$ is \'{e}tale,
$\GG(a,a)$ is discrete and hence so is $\cH(\tL)$. Take any $e\in\tL$.
Since $E$ is principal, the map $\xW$ is a local homeomorphism.
Therefore we can find an open neighbourhood $V$ of $e$ in $E$ such that
$\xW|_{V}$ is injective. Put $U=V\cap\tL$. Let $g\in\cH(\tL)$ and $e'\in U$
such that $g\cdot e'\in U$.
Now $\xW(g\cdot e')=\xW(e')$, and by injectivity of $\xW$ on $V\supset U$ it
follows that $g\cdot e'=e'$. But the action of $\GG$ is free, therefore
$g=1_{a}$. Thus $g\cdot U\cap U=\emptyset$ for any 
$g\in\cH(\tL)\setminus\{1_{a}\}$.
\eop

Let $\GG$ be an \'{e}tale groupoid, $\HH$ a topological groupoid
and $(E,p,\xW)$ a principal $\GG$-$\HH$-bibundle with locally
path-connected fibers. Let $e\in E$. Now $e$ lies in a unique
$\HH$-connected component $\tL$ of the fiber $E_{p(e)}$ of $E$.
Let $L$ be the associated leaf of $\tL$, i.e. the unique leaf of
$E$ with $\xW(e)\in L$.

Since $\cH(\tL)$ acts properly discontinuously on $\tL$, 
the $\HH$-equivariant map
$$ \zeta:\tL\lra\tL/\cH(\tL)\cong L $$
is a covering projection. In particular,
$L$ is locally path-connected. By the unique $\HH$-path
lifting property for $\zeta$ we obtain a surjective homomorphism
$$ \cH_{e}:\pi_{1}(\HH(L),\xW(e))\lra\cH(\tL)\subset\GG(p(e),p(e))\;,$$
called the {\em holonomy homomorphism} of the leaf $L$
with respect to the base point $e$.
If $\ell=\sigma_{n}\cdot h_{n}\cdot\ldots\cdot h_{1}\cdot\sigma_{0}$
is an $\HH$-loop in $L$ with the base-point $\xW(e)$
and
$$ \tilde{\sigma}_{n}\cdot h_{n}\cdot\ldots\cdot h_{1}
   \cdot\tilde{\sigma}_{0} $$
is the unique lift of $\ell$ in $\tL$ such that
$\tilde{\sigma}_{0}(0)=e$,
then $\cH_{e}([\ell])$ is the unique element of $\cH(\tL)$ such
that
$$ \cH_{e}([\ell])^{-1}\cdot e=\tilde{\sigma}_{n}(1)\;.$$
Here $[\ell]$ denotes the homotopy class of $\ell$.

If $e'$ is another point in $E$ with $\xW(e')=\xW(e)$,
the homomorphisms $\cH_{e}$ and $\cH_{e'}$ differ by the conjugation
by the uniquely determined $g\in\GG_{1}$ such that $g\cdot e=e'$.
If $e''$ is a point in $E$ such that $\xW(e'')\in L$, the homomorphisms
$\cH_{e}$ and $\cH_{e''}$ are again isomorphic, but now the isomorphism
depends on the choice of an $\HH$-path in $L$ from $\xW(e)$ to $\xW(e'')$.
In other words, each leaf $L$ of $E$ determines (up to an isomorphism)
the {\em holonomy homomorphism} $\cH_{L}$ of $L$, 
$$ \cH_{L}:\pi_{1}(\HH(L))\lra\GG\;.$$

\begin{theo}  \label{IItheo16}
Let $\psi:\GG\ra\GG'$ be a continuous functor between \'{e}tale groupoids,
and let $\phi:\HH\ra\HH'$ be a continuous functor
between topological groupoids. Let $(E,p,\xW)$ be a principal
$\GG$-$\HH$-bibundle, $(E',p',\xW')$ a principal $\GG'$-$\HH'$-bibundle
and $\alpha:E\ra E'$ a continuous $\psi$-$\phi$-equivariant map.
Assume that both $(E,p,\xW)$ and $(E',p',\xW')$ have locally path-connected
fibers. Let $e\in E$, let $L$ be the leaf of $(E,p,\xW)$ with $\xW(e)\in L$
and let $L'$ the leaf of $(E',p',\xW')$ with $\xW'(\alpha(e))\in L'$.
Then $\phi$ restricts to a continuous functor
$\phi|_{L}:\HH(L)\ra\HH'(L')$ and
$$ \psi\com\cH_{e}=\cH_{\alpha(e)}\com(\phi|_{L})_{\ast}\;.$$
Hence $\psi\com\cH_{L}\cong\cH_{L'}\com(\phi|_{L})_{\ast}$.
\end{theo}
\Proof
Put $e'=\alpha(e)$.
Let $\tL$ be the $\HH$-connected component of $E_{p(e)}$ with $e\in\tL$
and $\tL'$ the $\HH'$-connected component of
$E'_{p'(e')}$ with $e'\in\tL'$.
Proposition \ref{IIprop14} implies that $\alpha(\tL)\subset\tL'$
and $\phi_{0}(L)\subset L'$.
If $\ell=\sigma_{n}\cdot h_{n}\cdot\ldots\cdot h_{1}\cdot\sigma_{0}$
is an $\HH$-loop in $L$ with the base point $\xW(e)$,
there is the unique lift
$$ \tilde{\sigma}_{n}\cdot h_{n}\cdot\ldots\cdot 
   h_{1}\cdot\tilde{\sigma}_{0} $$
of $\ell$ in $\tL$ with $\tilde{\sigma}_{0}(0)=e$. Now
$$ (\alpha\com\tilde{\sigma}_{n})\cdot\phi(h_{n})\cdot\ldots\cdot
   \phi(h_{1})\cdot (\alpha\com\tilde{\sigma}_{0}) $$
is an $\HH$-path in $\tL'$ with the initial point $e'$, and is
exactly the lift of the $\HH$-loop
$$ (\phi|_{L})_{\#}(\ell)=
   (\phi_{0}\com\sigma_{n})\cdot\phi(h_{n})\cdot\ldots\cdot
   \phi(h_{1})\cdot(\phi_{0}\com\sigma_{0}) $$
in $L'$ with the base point $\xW'(e')$. Therefore
$$ \psi(\cH_{e}[\ell])^{-1}\cdot e'=\alpha(\cH_{e}[\ell]^{-1}
   \!\!\cdot e)=\alpha(\tilde{\sigma}_{n}(1))=
   \cH_{e'}((\phi|_{L})_{\ast}[\ell])^{-1}\cdot e'\;.$$
\eop

\begin{theo}  \label{IItheo17}
Let $\psi:\GG\ra\GG'$ be a continuous functor between \'{e}tale
groupoids such that $\psi_{0}$ is a local homeomorphism.
Let $\HH$ be a topological groupoid and let
$(E,p,\xW)$ be a principal $\GG$-$\HH$-bibundle
with locally path-connected fibers. Then the tensor product 
$(\angs{\psi}\otimes E,\angs{\psi}\otimes p,\angs{\psi}\otimes \xW)$
is a principal $\GG'$-$\HH$-bibundle with
locally path-connected fibers, and the leaves of
$E$ are precisely the leaves of $\angs{\psi}\otimes E$,
with the same leaf topology. 
\end{theo}
\Proof
Recall that $\angs{\psi}\otimes E$ is the orbit space of
$\angs{\psi}\times_{\GG_{0}} E\cong\GG'_{1}\times_{\GG'_{0}}E$
with a right action of $\GG$. 
Observe that the quotient projection 
$\GG'_{1}\times_{\GG'_{0}}E\ra\angs{\psi}\otimes E$ 
is a local homeomorphism.
Since $\xW$ and $\angs{\psi}\otimes \xW$
are local homeomorphisms
and $\xW=(\angs{\psi}\otimes \xW)\com\psi_{E}$, the map 
$\psi_{E}$ is also a local homeomorphism.

Let $a'\in\GG'_{0}$. The fiber $(\angs{\psi}\otimes E)_{a'}$
is the orbit space of the $\GG$-invariant subspace
$$ \cod^{-1}(a')\times_{\GG'_{0}} E\subset
   \GG'_{1}\times_{\GG'_{0}} E\;.$$
The space $\GG'^{a'}_{1}=\cod^{-1}(a')$ is discrete, and for any
fixed $g'\in\GG'^{a'}_{1}$ we have
$\{g'\}\times_{\GG'_{0}}E\cong
\bigcup_{a\in\psi_{0}^{-1}(\dom g')}E_{a}\subset E$.
Since $\psi_{0}^{-1}(\dom g')$ is also discrete, $E_{a}$ is open in 
$\bigcup_{a\in\psi_{0}^{-1}(\dom g')}E_{a}$, for any
$a\in\psi_{0}^{-1}(\dom g')$.
Since each fiber of $E$ is locally path-connected, this yields that
$\GG'^{a'}_{1}\times_{\GG'_{0}}E$ is locally path-connected,
and hence also $(\angs{\psi}\otimes E)_{a'}$ is locally
path-connected. Moreover, for any $a\in\psi_{0}^{-1}(a')$, the
restriction
$$ \psi_{E}|_{E_{a}}:E_{a}\lra (\angs{\psi}\otimes E)_{a'} $$
is a local homeomorphism. Let $a\in\psi_{0}^{-1}(a')$, 
let $\tL$ be an $\HH$-connected component
of $E_{a}$, and let $\tL'$ be the $\HH$-connected component of
$(\angs{\psi}\otimes E)_{a'}$ with
$\psi_{E}(\tL)\subset\tL'$. Since
$\tL$ is open in $E_{a}$ and $\tL'$ is open in
$(\angs{\psi}\otimes E)_{a'}$, the restriction
$$ \psi_{E}|_{\tL}:\tL\lra\tL' $$ 
is also a local homeomorphism.

Take $e\in\tL$ and assume that $g'\otimes e'$ is a point of
$\tL'$ which lies in the same path-connected component of $\tL'$
as $1_{\psi_{0}(p(e))}\otimes e$.
Since $(\angs{\psi}\otimes E)_{a'}$ is the orbit space
of $\GG'^{a'}_{1}\times_{\GG'_{0}}E$,
Proposition \ref{Iprop11} implies that
there exists a $\GG$-path
$$ \sigma_{n}\cdot g_{n}\cdot\ldots\cdot g_{1}\cdot\sigma_{0} $$
from $(1_{\psi_{0}(p(e))},e)$ to $(g',e')$ in
$\GG'^{a'}_{1}\times_{\GG'_{0}}E$.
Since $\GG'^{a'}_{1}$ is discrete, it follows that
for any $i=0,1,\ldots,n$ we have
$$ \sigma_{i}(t)=(g'_{i},\rho_{i}(t))\;\;\;\;\;\;\;\; t\in [0,1]\;,$$
for an element $g'_{i}\in\GG'^{a'}_{1}$ and a path $\rho_{i}$ in $E$.
In particular, $\psi_{0}(p(\rho_{i}(t)))=\dom g'_{i}$ for any
$t\in [0,1]$, and since $\psi_{0}$ is a local homeomorphism,
$p\com\rho_{i}$ is a constant function as well,
for all $i=0,1,\ldots,n$.
Further, $g_{i}^{-1}\cdot\rho_{i}(0)=\rho_{i-1}(1)$ and
$$ g'_{i}=\psi(g_{i}\com\ldots\com g_{1})^{-1} $$
for any $i=1,2,\ldots,n$.
Now the paths $((g_{i}\com\ldots\com g_{1})^{-1}\cdot\oo)\com\rho_{i}$
concatenate in a path from $e$ to 
$(g_{n}\com\ldots\com g_{1})^{-1}\cdot e'$ in $\tL$, and
$$ \psi_{E}((g_{n}\com\ldots\com g_{1})^{-1}\cdot e')
   =g'\otimes e'\;.$$
This proves that any path-connected component of $\tL'$ which
intersects $\psi_{E}(\tL)$ lies in $\psi_{E}(\tL)$.
Since $\psi_{E}$ is $\HH$-equivariant, this implies that
$\psi_{E}(\tL)=\tL'$. Now since $\psi_{E}|_{\tL}:\tL\ra\tL'$ is an
open surjection, it induces a homeomorphism between the corresponding
leaves.
\eop

Let $\GG$ be an \'{e}tale groupoid, $\HH$ a topological groupoid and
$E$ a principal $\GG$-$\HH$-bibundle with locally path-connected fibers.
Theorem \ref{IItheo17} shows that the decomposition of $\HH_{0}$ into
the leaves of $E$ is invariant under the left composition of $E$ with
$\angs{\psi}$ if $\psi_{0}$ is a local homeomorphism.
In particular, the effect-functor
$$ \eee:\GG\lra\Gm(\GG_{0})\;,$$
described in Section \ref{secTopGro}, is the identity on objects.
Thus $\angs{\eee}\otimes E$ has the same leaves as $E$. However, the
holonomy groups may change.

To be more precise,
take $e\in E$, and let $L$ be the leaf of $E$ with $\xW(e)\in L$.
We define the {\em geometric holonomy homomorphism} of
$L$ with respect to the base point $e$ to be the composition
$$ \eee\cH_{e}=\eee\com\cH_{e}:\pi_{1}(\HH(L),\xW(e))\lra
   \Gm(\GG_{0})(p(e),p(e))\;.$$
As the holonomy homomorphism, the geometric holonomy homomorphism
depends up to an isomorphism only on the leaf $L$.
Thus we can define (up to an isomorphism) the {\em geometric holonomy
homomorphism} $\eee\cH_{L}:\pi_{1}(\HH(L))\ra\Gm(\GG_{0})$ of $L$.
The {\em geometric holonomy group} $\eee\cH(L)$ of
$L$ is the image of $\eee\cH_{L}$ in $\Gm(\GG_{0})$.

Now Theorem \ref{IItheo16}, applied to the
$\eee$-$\id_{\!\HH}$-equivariant map $\eee_{E}$ implies that the
holonomy group of $L$ regarded as
a leaf of $\angs{\eee}\otimes E$ is exactly the 
geometric holonomy group of $L$.

\begin{ex}  \label{IIex18}  \rm
(1) Let $\GG$ be an \'{e}tale groupoid and $X$ a topological space.
Let $\ccc$ be a Haefliger $\GG$-structure on $X$ represented
by a $\GG$-cocycle $c=(c_{ij})$ on an open cover
$\cU=(U_{i})_{i\in I}$ of $X$, and let $(\Sigma(c),p,\xW)$
be the corresponding principal $\GG$-bundle over $X$.

Recall that $\Sigma(c)$ is a quotient of the space
$\bar{\Sigma}(c)=\GG_{1}\times_{\GG_{0}}\coprod_{i\in I}U_{i}$.
The quotient projection $q:\bar{\Sigma}(c)\ra\Sigma(c)$
is a local diffeomorphism. Indeed, observe that
if $V\subset\GG_{1}$ is such that both $\dom|_{V}$ and $\cod|_{V}$
are injective, then both $q|_{V\times_{\GG_{0}}U_{i}}$ and
$\xW|_{q(V\times_{\GG_{0}}U_{i})}$ are injective, for any $i\in I$.
Therefore, it is clear that the principal $\GG$-bundle
$\Sigma(c)$ has locally path-connected fibres if and only if
the maps $c_{ii}$ have locally path-connected fibres.

Assume now that the maps
$c_{ii}$ have locally path-connected fibres.
We will describe the leaves and holonomy of $\Sigma(c)$ in
terms of $c$. Choose $x\in X$ and let
$L$ be the leaf of $\Sigma(c)$ with $x\in L$.
Using the fact that the open subsets
of the form $q(V\times_{\GG_{0}} U_{i})$ as above cover
$\Sigma(c)$, it is obvious that a point
$x'\in X$ belongs to $L$ if and only if there exists a path
$\sigma:[0,1]\ra X$ from $x$ to $x'$, with a partition
$0=t_{0}<t_{1}<\ldots <t_{n-1}<t_{n}=1$ such that
$$ \sigma[t_{k-1},t_{k}]\subset c_{i_{k}i_{k}}^{-1}(a_{k}) $$
for some $i_{1},\ldots,i_{n}\in I$ and
$a_{1},\ldots,a_{n}\in\GG_{0}$. Any path in $L$ is of
this form.

Now assume that we have a loop $\sigma$ in $L$ with the
base-point $x$. Choose $i_{0}\in I$ with $x\in U_{i_{0}}$,
and put $e=[1_{c_{i_{0}i_{0}}(x)},x,i_{0}]\in\Sigma(c)$.
Let $\tL$ be the connected component of $\Sigma(c)_{p(e)}$
with $e\in\tL$. We will compute $\cH_{e}([\sigma])$.

First, choose a partition
$0=t_{0}<t_{1}<\ldots <t_{n-1}<t_{n}=1$,
$i_{1},\ldots,i_{n}\in I$ and $a_{1},\ldots,a_{n}\in\GG_{0}$
such that
$$ \sigma[t_{k-1},t_{k}]\subset c_{i_{k}i_{k}}^{-1}(a_{k})\;.$$
Take $g_{k}=c_{i_{k-1},i_{k}}(\sigma(t_{k-1}))$ and
$\varrho_{k}=g_{1}\com g_{2}\com\ldots\com g_{k}$, for
$k=1,2,\ldots,n$.
Define $\sigma_{k}:[t_{k-1},t_{k}]\ra\bar{\Sigma}(c)$ by
$$ \sigma_{k}(t)=(\varrho_{k},\sigma(t),i_{k})\;\;\;\;\;\;\;\; 
   k=1,2,\ldots,n\;.$$
The paths $q\com\sigma_{k}$ amalgamate in a path
$\tilde{\sigma}:[0,1]\ra\tilde{L}$. Clearly
$\xW\com\tilde{\sigma}=\sigma$ and $\tilde{\sigma}(0)=e$.
Since
$$ \tilde{\sigma}(1)=[\varrho_{n},x,i_{n}]=
   [\varrho_{n}\com c_{i_{n},i_{0}}(x),x,i_{0}]\;,$$
this yields that
$\cH_{e}([\sigma])=c_{i_{0},i_{n}}(x)\com\varrho^{-1}_{n}$,
that is
$$ \cH_{e}([\sigma])=
   c_{i_{0},i_{n}}(\sigma(t_{n}))\com
   c_{i_{n},i_{n-1}}(\sigma(t_{n-1}))\com\ldots\com
   c_{i_{1},i_{0}}(\sigma(t_{0}))\;.$$

(2) Let $\GG$ be an \'{e}tale groupoid, $\HH$ a topological
groupoid and $\phi:\HH\ra\GG$ a continuous functor.
A fiber of $\angs{\phi}$ over $a\in\GG_{0}$ is
$\dom^{-1}(a)\times_{\GG_{0}}\HH_{0}\subset\angs{\phi}$.
Since $\dom^{-1}(a)$ is discrete, the fibers of
$\angs{\phi}$ are locally path-connected if and only if
the fibers of $\phi_{0}$ are locally path-connected.

Assume that the fibers of $\phi_{0}$ are locally
path-connected. Then a leaf of $\angs{\phi}$ is
a minimal union of connected components of the fibers
of $\phi_{0}$ which is $\HH$-invariant. Each connected
component of a fiber of $\phi_{0}$ is open in the
corresponding leaf with respect to the leaf topology,
and embedded in $\HH_{0}$.

The holonomy of an $\HH$-loop
$\ell=\sigma_{n}\cdot h_{n}\cdot\ldots\cdot h_{1}\cdot\sigma_{0}$
in a leaf $L$ of $\angs{\phi}$ with the base-point $x\in L$
is clearly
$$ \cH_{e}([\ell])=\phi(h_{n})\com\phi(h_{n-1})\com\ldots\com
   \phi(h_{1})\;,$$
where $e=(1_{\phi_{0}(x)},x)$.

(3) Let $G$ be a discrete group acting continuously on
a topological space $Y$. Let $\HH$ be a topological groupoid
with $\HH_{0}$ locally path-connected, and assume that
$|\HH|$ is connected. Let $(E,p,\xW)$ be a principal
$G(Y)$-$\HH$-bibundle.
Observe that there is a canonical continuous functor
between \'{e}tale groupoids
$$ \ppp:G(Y)\lra G\;.$$
Therefore $\angs{\ppp}\otimes E$ is a principal
$G$-$\HH$-bibundle. Note however that as a space,
$\angs{\ppp}\otimes E$ is homeomorphic to $E$ and
this homeomorphism is $\ppp$-$\id_{\!\HH}$-equivariant.
The underlying $G$-bundle over $\HH_{0}$ of the
bibundle $\angs{\ppp}\otimes E$ is principal, therefore
$\xW$ is a covering projection.

Now since $G_{0}$ is a one-point space and $|\HH|$ is
connected, it is clear that $E$ has only one leaf, i.e.
$\HH_{0}$. We have the holonomy homomorphism of this leaf
$$ \cH_{\HH_{0}}:\pi_{1}(\HH)\lra G\;.$$
As we know, with a base-point in $\HH_{0}$ fixed, this
homomorphism is determined up to a conjugation in $G$.
This is closely related to the notion of development.
In fact, if $\HH$ is just a space $X$ and if $\tilde{X}$
is an universal covering space of $X$, the development map
can be defined as $p\com\zeta:\tilde{X}\ra Y$, where
$\zeta:\tilde{X}\ra E$ is a covering projection.
\end{ex}

\chapter{Stability}
\label{chapSta}
\startchapterskip

A Haefliger $\GG$-structure on a topological space $X$ can be
seen as a generalized foliation on $X$. Moreover, the
Haefliger-Reeb-Ehresmann stability theorem \cite{Hae} generalizes
the Reeb stability theorem for foliations on manifolds
(Theorem \ref{Itheo1}).
In Chapter \ref{chapHilSkaMap} we showed that
in fact any transitive $\GG$-$\HH$-bibundle -- and hence any
Hilsum-Skandalis map between topological groupoids $\HH$ and
$\GG$ -- can be seen as a generalized foliation on $\HH$.

In this Chapter we prove a stability theorem for a
transitive $\GG$-$\HH$-bibundle, which in particular
generalizes the Haefliger-Reeb-Ehresmann stability theorem.
Furthermore, we show that the Reeb-Thurston stability theorem
(Theorem \ref{Itheo2})
can be extended to transitive $\GG$-$\HH$-bibundles in the case
where $\GG$ is an \'{e}tale $\eCe$-groupoid. In particular, we give a
version of the Reeb-Thurston stability theorem for
Haefliger $\GG$-structures.
In Chapter \ref{chapEquFol} we shall apply
these results to the $\eCe$-principal
$\GmqCe$-$G(M)$-bibundle associated to
a $G$-equivariant foliation on a $\eCe$-manifold $M$
(Example \ref{IIex8} (3)).

\section{Reeb Stability for a Transitive Bibundle}
\label{secReeSta}

To prove the Reeb stability theorem for a transitive bibundle,
we first deal with a special case of a transitive
$Y$-$X$-bibundle, i.e. a continuous map $p:X\ra Y$, where $X$ and
$Y$ are topological spaces. Recall that a leaf of a continuous map
$p:X\ra Y$ is a connected component of a fiber
of $p$ with the inherited topology from $X$, and the holonomy
group of a leaf of $p$ is trivial.

\begin{lem}  \label{IIIlem1}
Let $p:X\ra Y$ be a continuous map between topological spaces,
$K$ a compact subset of $X$ and $U$ an open subset of $X$ with
$U\subset K$. If $L$ is a closed Hausdorff leaf of $p$ such that
$L\cap K=L\cap U\neq\emptyset$, then $L\subset U$.
\end{lem}
\Proof
Since $L$ is closed, $L\cap K$ is closed in the compact $K$ and hence
compact. Since $L$ is Hausdorff, $L\cap K$ is closed in $L$. But
$L\cap K=L\cap U$ is also open in $L$ and non-empty, and therefore
$L\cap U=L$ by connectedness of $L$.
\eop

\begin{lem}  \label{IIIlem2}
Let $X$ be a locally compact space and $K$ a closed compact subset
of $X$. Then $X/K$ is locally compact and the restriction of the
quotient projection $q:X\ra X/K$ to the open subset $X\setminus K$
of $X$ is an open embedding.
\end{lem}
\Proof
If $U$ is an open neighbourhood of $\{K\}$ in $X/K$,
then $q^{-1}(U)$ is an open neighbourhood of $K$ in $X$. Since $X$
is locally compact and $K$ is compact,
there exists a compact neighbourhood $W\subset q^{-1}(U)$ of
$K$. Now $K\subset\Int(W)$, so $\Int(W)$ is saturated
and open. Therefore $q(\Int(W))$ is open in $X/K$,
$q(W)$ is a compact neighbourhood of $\{K\}$ in $X/K$ and
$q(W)\subset U$. The rest of the lemma is a trivial consequence
of the fact that $K$ is closed.
\eop

\begin{theo}[Reeb stability for a continuous function]
\label{IIItheo3}
Let $X$ be a locally compact space, $Y$ a locally Hausdorff space
and $p:X\ra Y$ a continuous map with Hausdorff leaves. Let $L$ be
a compact leaf of $p$ which is open in the fiber $p^{-1}(p(L))$.
Then for any open neighbourhood $V$ of $L$ there exists an open
neighbourhood $U\subset V$ of $L$ which is a union of compact
leaves of $p$.
\end{theo}
\Proof
Clearly we can assume without loss of generality that $Y$ is
Hausdorff. Since the points of $Y$ are closed, the fibers of $p$
are closed and therefore the leaves of $p$ are also closed.

(a) Assume first that $L=\{x\}$ for some point $x\in X$, and
denote $y=p(x)$. Since $L$ is open in $p^{-1}(y)$, there exists
an open neighbourhood $W$ of $x$ in $X$ such that
$W\cap p^{-1}(y)=\{x\}$. Since $X$ is locally compact there
is a compact neighbourhood $K$ of $x$ with $K\subset V\cap W$.
Put $R=K\setminus\Int(K)$.
Now $\Int(K)$ is open in $X$ and hence also in $K$, so $R$ is
closed in the compact $K$ and therefore compact. Hence $p(R)$
is compact in $Y$, and since
$Y$ is assumed to be Hausdorff, this yields that $p(R)$ is closed
in $Y$. Since $W\cap p^{-1}(y)=\{x\}$ and $x\in\Int(K)$,
we have $R\cap p^{-1}(y)=\emptyset$ and hence
$y\not\in p(R)$. Take
$$ U=\Int(K)\cap p^{-1}(Y\setminus p(R))\;.$$
The set $U$ is clearly an open neighbourhood of $x$ and
$U\subset K\subset V$. Let $L'$ be a leaf of $p$
such that $L'\cap U\neq\emptyset$. It follows that
$p(L')\cap p(R)=\emptyset$
and hence $L'\cap R=\emptyset$. This yields that
$L'\cap K=L'\cap U\neq\emptyset$, so
Lemma \ref{IIIlem1} implies that $L'\subset U$. Moreover,
$L'$ is a closed subset of the compact $K$ and therefore compact.
It follows that $U$ is a union of compact leaves of $p$.

(b) For the general case, denote by $q:X\ra X/L$ the quotient
projection. Now $p=p'\com q$ for a unique continuous function
$p':X/L\ra Y$. By
Lemma \ref{IIIlem2}, $X/L$ is locally compact and
$q|_{X\setminus L}$ is an open embedding. Now $\{L\}$ is a leaf
of $p'$. For any leaf $L'$ of $p'$,
$q^{-1}(L')$ is a leaf of $p$, and if $L'\neq\{L\}$
then $q^{-1}(L')$ is homeomorphic
to $L'$. The set $V$ is saturated and hence $q(V)$ is an open
neighbourhood of $\{L\}$ in $X/L$. Clearly all leaves of $p'$
are Hausdorff and $\{L\}$ is open in $p'^{-1}(p(L))$, so we can
use the part (a) to find an open neighbourhood $U'\subset q(V)$ of
$\{L\}$ which is a union of compact leaves of $p'$. Then we take
$$ U=q^{-1}(U') $$
and the proof is complete.
\eop

\begin{cor}  \label{IIIcor4}
Let $X$ be a locally compact space, $Y$ a locally Hausdorff
space and $p:X\ra Y$ a continuous map with locally connected
Hausdorff fibers. Let $L$ be a compact leaf of $p$. Then for any
open neighbourhood $V$ of $L$ there exists an open neighbourhood
$U\subset V$ of $L$ which is a union of compact leaves of $p$.
\end{cor}
\Proof
Since the fibers of $p$ are locally connected, the leaves are open
in the corresponding Hausdorff fibers, and hence also Hausdorff.
\eop

\begin{theo}[Reeb stability for a transitive bundle]
\label{IIItheo5}
Let $\GG$ be a topological groupoid and $(E,p,\xW)$ a
transitive $\GG$-bundle over a topological space $X$. 
Assume that $E$ is locally compact, $\GG_{0}$
locally Hausdorff and $(E,p,\xW)$ with locally connected
Hausdorff fibers. Let $L$ be a compact leaf of $(E,p,\xW)$
with compact holonomy group. 
Then for any open neighbourhood $V$ of $L$ in $X$
there exists an open neighbourhood $U\subset V$ of $L$ in $X$ 
which is a union of compact leaves of $(E,p,\xW)$.
\end{theo}
\Proof
Choose $a\in\GG_{0}$ and a connected component $\tL$ of the
fiber $E_{a}=p^{-1}(a)$ such that $L$ is associated to $\tL$.
Since $\cH(\tL)$ is compact and $L\cong\tL/\cH(\tL)$ is compact,
$\tL$ is also compact.
Applying Corollary \ref{IIIcor4} to the map $p:E\ra\GG_{0}$ we
can find an open neighbourhood
$$ \tilde{U}\subset \xW^{-1}(V) $$
of $\tL$ in $E$ which is a union of compact leaves of $p$, i.e.
of compact connected components of the fibers of $(E,p,\xW)$.
Since $\xW$ is open,
$$ U=\xW(\tilde{U})\subset V $$
is an open neighbourhood of $L$ in $X$ which is a union of
compact leaves of $(E,p,\xW)$.
\eop

\begin{theo}[Reeb stability for a transitive bibundle]
\label{IIItheo6}
Let $\GG$ and $\HH$ be topological gro\-upoids and $(E,p,\xW)$
a transitive $\GG$-$\HH$-bibundle. Assume that $E$ is locally
compact, $\GG_{0}$ locally Hausdorff and $(E,p,\xW)$ with locally
connected $\HH$-Hausdorff fibers. Let $L$ be an $\HH$-compact
leaf of $(E,p,\xW)$ with compact holonomy group. 
Then for any $\HH$-invariant open neighbourhood $V$ of $L$ in
$\HH_{0}$ there exists an $\HH$-invariant open neighbourhood
$U\subset V$ of $L$ in $\HH_{0}$ which is a union of
$\HH$-compact leaves of $(E,p,\xW)$.
\end{theo}
\Proof
Note that the associated $\GG$-bundle $E/\HH$ over $|\HH|$ is
locally compact with locally connected Hausdorff fibres.
By Proposition \ref{IIprop13}, $L/\HH$ is a compact leaf
of $(E/\HH,p/\HH,\xW/\HH)$ with compact holonomy group,
and $V/\HH$ is an open neighbourhood
of $L/\HH$ in $|\HH|$. The result now follows from
Theorem \ref{IIItheo5} and Proposition \ref{IIprop13}.
\eop
\Rem
In particular, Theorem \ref{IIItheo6} extends the Reeb stability
theorem to principal $\GG$-$\HH$-bibundles, i.e. to the
Hilsum-Skandalis maps from $\HH$ to $\GG$. 
In the special cases, it gives a version
of the Reeb stability theorem for
nested foliations (Example \ref{IIex8} (4))
and for foliations invariant under a group action
(Example \ref{IIex8} (3)), which we will discuss extensively
in Chapter \ref{chapEquFol}.

\begin{theo}[Reeb stability for a principal bundle]
\label{IIItheo7}
Let $\GG$ be an \'{e}tale groupoid with $\GG_{0}$ locally
Hausdorff, and let $(E,p,\xW)$ be a principal 
$\GG$-bundle over a locally compact space $X$
with locally path-connected Hausdorff fibers. 
Let $L$ be a compact leaf of $(E,p,\xW)$ with finite
holonomy (respectively geometric holonomy) group. 
Then for any open neighbourhood $V$ of $L$ in $X$
there exists an open neighbourhood $U\subset V$ of $L$ in
$X$ which is a union of compact leaves of $(E,p,\xW)$ with
finite holonomy (respectively geometric holonomy) groups.
\end{theo}
\Rem
Observe that if $X$ is Hausdorff,
the assumption that the fibers of $E$ are Hausdorff
is superfluous.
\vspace{4 mm}

\Proof
(a) Assume first that the holonomy of $L$ is finite.
Since $X$ is locally compact and $\xW$ is a local
homeomorphism, $E$ is also locally compact.
Choose $e\in E$ with $\xW(e)\in L$, and let $\tL$ be the
connected component of the fiber 
$E_{p(e)}$ with $e\in\tL$. Since $\cH(\tL)$ is
finite and $L\cong\tL/\cH(\tL)$ is compact, $\tL$ is also
compact. Applying Corollary \ref{IIIcor4} to the map $p$ we
can find an open neighbourhood
$$ \tilde{U}\subset \xW^{-1}(V) $$
of $\tL$ in $E$ which is a union of compact leaves of $p$,
i.e. of compact connected components of the fibers of $(E,p,\xW)$. 
Proposition \ref{IIprop15} implies that the holonomy group
of a compact connected component of a fiber of $(E,p,\xW)$ is
finite. Thus
$$ U=\xW(\tilde{U})\subset V $$
is an open neighbourhood of $L$ in $X$ which is a union of
compact leaves of $(E,p,\xW)$ with finite holonomy groups. 

(b) Assume that the geometric holonomy of $L$ is finite.
Recall from Section \ref{secTopGro} that we have the
effect-functor
$$ \eee:\GG\lra\Gm(\GG_{0})\;,$$
which is the identity on objects.
By Theorem \ref{IItheo16} and Theorem \ref{IItheo17},
the tensor product $\angs{\eee}\otimes E$
has locally path-connected fibers and exactly
the same leaves as $E$. Moreover, the 
geometric holonomy group of a leaf of $E$ is exactly
the holonomy group of the same leaf regarded as a leaf of
$\angs{\eee}\otimes E$.
Since $E$ has Hausdorff fibers, the leaves of $E$ are
Hausdorff, and therefore the fibers of 
$\angs{\eee}\otimes E$ are also Hausdorff. 
The theorem thus follows from the part (a).
\eop
\Rem
In the case where $\GG$ is effective,
Theorem \ref{IIItheo7} is precisely the Haefliger-Reeb-Ehresmann
theorem for Haefliger $\GG$-structures \cite{Hae}.
In particular, it generalizes the classical Reeb stability theorem
for foliations on manifolds, and
gives a version of the Reeb stability theorem
for foliations with modular singularities
(Example \ref{IIex8} (2)).

\section{Differential Categories}  \label{secDifCat}

In Section \ref{secReeThuSta} we will
generalize the results of Section \ref{secReeSta}
in a similar way Thurston generalized the Reeb stability
theorem \cite{Thu} (see Theorem \ref{Itheo2}).
The proof of the Thurston generalization can be essentially
reduced to a theorem concerning groups of germs of
diffeomorphisms of $\RRR^{q}$ with differential 1,
which was simplified in \cite{ReeSch,Scha}.
In Subsection \ref{subsecRepGro} we prove the
infinite-dimensional version of that theorem, in
order to prove in Section \ref{secReeThuSta} the
Reeb-Thurston stability for transitive $\GG$-$\HH$-bibundles,
where $\GG$ is an \'{e}tale $\eCe$-groupoid.
To facilitate the proof,
we first introduce differential categories and
investigate some of their properties. 
\vspace{4 mm}

\Not
In this section, $\FFF$ stands for the field $\RRR$
or $\CCC$, and all the vector spaces are over $\FFF$ unless
we specify the field explicitly.

\subsection{Definition of a Differential Category}
\label{subsecDef}

\begin{dfn}  \label{IIIdfn8}
A right linear category over $\FFF$ is a category $\gC$ which 
satisfies
\begin{enumerate}
\item $\gC(a,a')$ is a vector space over $\FFF$, for any
      $a,a'\in \Ob\gC$ (we denote by $0_{a',a}$ the neutral
      element of $\gC(a,a')$),
\item $\oo\com c:\gC(a',a'')\ra\gC(a,a'')$ is linear, for any
      $a,a',a''\in\Ob\gC$ and any $c\in\gC(a,a')$,
\item $\gC(a',a'')\com 0_{a',a}=0_{a'',a}$ for any
      $a,a',a''\in\Ob\gC$.
\end{enumerate}
\end{dfn}
\Rem A linear category is a right linear category in which the obvious
duals of the conditions 2. and 3. are satisfied as well.
\vspace{4 mm}

\Not
Let $\gC$ be a right linear category.
A morphism $A\in\gC(a',a'')$ is linear if 
$A\com\oo:\gC(a,a')\ra\gC(a,a'')$ is linear for any
$a\in\Ob\gC$.
For example, $0_{a'',a'}$ and the identity morphism
$1_{a'}\in\gC(a',a')$ are linear. The linear morphisms of $\gC$
clearly form a linear subcategory of $\gC$, which will be denoted
by $\cL\gC$.

\begin{dfn}  \label{IIIdfn9}
Let $\gC$ and $\gD$ be right linear categories.
A functor $F:\gC\ra\gD$ is linear if
$$ F|_{\gC(a,a')}:\gC(a,a')\ra\gD(Fa,Fa') $$
is linear, for any $a,a'\in \Ob\gC$.
\end{dfn}

\begin{dfn}  \label{IIIdfn10}
Let $\gC$ be a right linear category. A (semi)norm
on $\gC$ is a function
$$ \|\oo\| : \!\!\!\!\coprod_{a,a'\in\Ob\gC}\!\!\!\!
   \gC(a,a')\lra [0,\infty ) $$
such that $\|\oo\| \,|_{\gC(a,a')}$ is a (semi)norm on
$\gC(a,a')$ for any $a,a'\in\Ob\gC$.
A (semi)norm $\|\oo\|$ on $\gC$ is
submultiplicative if $\|1_{a}\|=1$ and
$$ \|c'\com c\|\leq\|c'\| \|c\| $$
for any $c\in\gC(a,a')$, $c'\in\gC(a',a'')$,
$a,a',a''\in\Ob\gC$.
A (semi)normed (right) linear category is a pair
$(\gC,\|\oo\|)$, where
$\gC$ is a (right) linear category and $\|\oo\|$ is a
submultiplicative (semi)norm on $\gC$.
A (right) Banach category
is a normed (right) linear category $(\gB,\|\oo\|)$
such that the normed space
$$ (\gB(a,a'),\|\oo\| \,|_{\gB(a,a')}) $$
is complete, for any $a,a'\in\Ob\gB$.
\end{dfn}

\Rem 
Banach categories need not be additive.

\begin{ex}  \label{IIIex11}  \rm
(1) A Banach algebra $\cA$ with identity is a
Banach category, with $\Ob\cA$ an one-point set.

(2) The category $\Ban=\Ban_{\FFF}$ of Banach spaces
over $\FFF$ and bounded linear operators is a Banach category
over $\FFF$.

(3) Let $\Ob\gC=\Ob\Ban$ and let $\gC(\cX,\cY)$ be the
vector space of continuous functions $f:\cX\ra\cY$, for any
pair of Banach spaces $\cX,\cY$.
Then $\gC$ is a right linear category. We can define a function
$\|\oo\|_{\infty}:\coprod_{\cX,\cY\in\Ob\gC}\gC(\cX,\cY)\ra
[0,\infty]$ by
$$ \|f\|_{\infty}=\sup_{x\in \cX}\|f(x)\| $$
for any $f\in\gC(\cX,\cY)$. This function satisfies the usual
properties of a norm, and for any $f\in\gC(\cX,\cY)$ and
$f'\in\gC(\cY,\cZ)$ we have
$\|f'\com f\|_{\infty}\leq\|f'\|_{\infty}$. 
\end{ex}

\begin{dfn}  \label{IIIdfn12}
Let $\gC$ be a right linear category.
A $\gB$-linearization of $\gC$ is a linear functor
$$ D :\gC\lra\gB\;,$$
where $\gB$ is a Banach category.
A $\gB$-linearization $D$ of $\gC$ is 
called a $\gB$-differential on $\gC$
if there exists a sequence
$(\|\oo\|_{n})_{n=1}^{\infty}$ of norms on $\gC$ such that
$$ \lim\sup_{\!\!\!\!\!\!\!\!\!\!\!\! n\ra\infty}
   \frac{\|c\com c'-c\com c''\|_{n}}{\|c'-c''\|_{n}}
   \leq\|Dc\| $$
for any $c\in\gC(a',a'')$ and any two distinct
$c',c''\in\gC(a,a')$, $a,a',a''\in\Ob\gC$.
A differential category is a triple
$(\gC,D,\gB)$, where $\gC$ is a right linear category and $D$
is a $\gB$-differential on $\gC$.
\end{dfn}
\Rem 
In Subsection \ref{subsecDifCatAss} we will show that 
differential categories naturally arise from the
affine manifolds.

\begin{lem}  \label{IIIlem13}
Let $\gC$ and $\gD$ be right linear categories,
$D$ a $\gB$-differential on $\gD$ and $F:\gC\ra\gD$ a 
faithful linear functor. Then $D\com F$ is a
$\gB$-differential on $\gC$.
\end{lem}
\Proof
Since $F$ is faithful,
$F|_{\gC(a,a')}:\gC(a,a')\ra\gD(Fa,Fa')$ is injective,
for any $a,a'\in\Ob\gC$. Since $D$ is a $\gB$-differential, 
there exists a family of norms 
$(\|\oo\|_{n})_{n=1}^{\infty}$ on $\gD$
as in Definition \ref{IIIdfn12}. This family
induces a family of norms $(\|F\oo\|_{n})_{n=1}^{\infty}$ on
$\gC$. Now for any $c\in\gC(a',a'')$ and any two distinct
$c',c''\in\gC(a,a')$ we have
$$ \lim\sup_{\!\!\!\!\!\!\!\!\!\!\!\! n\ra\infty}
   \frac{\|F(c\com c'-c\com c'')\|_{n}}{\|F(c'-c'')\|_{n}}=
   \lim\sup_{\!\!\!\!\!\!\!\!\!\!\!\! n\ra\infty}
   \frac{\|Fc\com Fc'-Fc\com Fc''\|_{n}}{\|Fc'-Fc''\|_{n}}
   \leq\|D Fc\|\;,$$
and hence $D\com F$ is a $\gB$-differential.
\eop

\begin{ex}  \label{IIIex14}  \rm
Define a category $\gC$ with
$\Ob\gC=\Ob\Ban$ such that for any $\cX,\cY\in\Ob\Ban$,
$\gC(\cX,\cY)$ is the vector space
of functions $f:\cX\ra\cY$ which are continuously
differentiable on some open neighbourhood of $0\in \cX$,
preserve bounded sets and satisfy $f(0)=0$. Clearly
$\gC$ is a right linear category. Note that a morphism $f$ of
$\gC$ is linear if and only if it is a bounded linear operator,
so $\cL\gC=\Ban\subset\gC$.

If $f\in\gC(\cX,\cY)$ and $x\in \cX$ such that $f$ is
differentiable at $x$, we denote
by $df_{x}\in\Ban(\cX,\cY)$ the differential of $f$ at $x$.
We have a natural $\Ban$-linearization
$D$ on $\gC$ given as the identity on objects and by
$$ D(f)=df_{0} $$
for any $f\in\gC(\cX,\cY)$. Observe that $D^{2}=D$ and
$D(\gC)=\Ban$, i.e. $D$ is a projector on $\Ban$.

Now define for any $r\in [0,\infty)$ a seminorm on
$\gC$ by
$$ \|f\|_{r}=\sup_{\|x\|\leq r}\|f(x)\| $$
for any $f\in\gC(\cX,\cY)$. This is well-defined since
$f$ preserves bounded sets.
Since $f$ is continuous on an open neighbourhood of $0$,
we have
$$ \lim_{r\ra 0}\|f\|_{r}=\|f\|_{0}=\|f(0)\|=0 \;.$$
Note that if $f'\in\gC(\cY,\cZ)$, then
$$ \|f'\com f\|_{r}\leq \|f'\|_{\|f\|_{r}} $$
for all $r\geq 0$. Moreover, if $A\in\Ban(\cX,\cY)$, then
$\|A\|_{r}=r\,\|A\|$.
It follows that $\|\oo\|_{r}$ is actually a norm on $\Ban$,
for any $r>0$. We can define also $\|f\|=\|D(f)\|$ for any 
$f\in\gC(\cX,\cY)$, obtaining a
submultiplicative seminorm on $\gC$ which extends the standard
norm on $\Ban$. From the definition of differential
$df_{0}$ one can easily see that
$$ \|f\|=\|D(f)\|=\|df_{0}\|=
   \lim_{r\ra 0}\frac{\|f\|_{r}}{r}\;. $$
It follows that $D(f)$ is the unique linear morphism in
$\gC(\cX,\cY)$ which satisfies 
$$ \lim_{r\ra 0}\frac{\|f-D(f)\|_{r}}{r}=0\;.$$
Moreover, if $f',f''\in\gC(\cX,\cY)$ and $f\in\gC(\cY,\cZ)$,
the mean value theorem \cite{Lan} gives us for any $x\in \cX$
small enough
$$ \|f(f'(x))-f(f''(x))\|\leq\|f'(x)-f''(x)\|\,
   \sup_{x'\in S}\|df_{x'}\|\;, $$
where $S$ is the segment between $f'(x)$ and $f''(x)$.
It follows that for small $r$
$$ \|f\com f'-f\com f''\|_{r}\leq\|f'-f''\|_{r}
   \sup_{x'\in K}\|df_{x'}\|\;,$$
where now
$K=\{x'\in \cX\,|\,\|x'\|\leq\max
(\,\|f'\|_{r},\|f''\|_{r})\,\}$.
If $\germ_{0}f'\neq\germ_{0}f''$,
then $\|f'-f''\|_{r}>0$ for any $r>0$, and
by the continuity of the differential $df_{x'}$ around $0$ we get
\begin{eqnarray}  \label{presentability}
 \lim\sup_{\!\!\!\!\!\!\!\!\!\!\!\! r\ra 0}
 \frac{\|f\com f'-f\com f''\|_{r}}{\|f'-f''\|_{r}}\leq\|D(f)\|\;.
\end{eqnarray}

Define now a category $\cBan=\cBan_{\FFF}$
with $\Ob\cBan=\Ob\Ban$ such that
for any $\cX,\cY\in\Ob\Ban$, $\cBan(\cX,\cY)$ is the vector space
of germs at $0$ of continuously differentiable functions
$f:U\ra \cY$, defined on an open neighbourhood 
$U\subset \cX$ of $0$ with $f(0)=0$.
This is also a right linear category,
and the natural functor
$$ \germ_{0}:\gC\lra\cBan $$
is linear, the identity on objects and full on morphisms. 
Note also that we can identify $\cL\cBan=\Ban$.
Moreover, $D$ factors through $\germ_{0}$, so we obtain a 
$\Ban$-linearization of $\cBan$, which we will denote by $\cD$. 

For any $\cX,\cY\in\Ob\Ban$, we can choose a linear injection
$$ \alpha=\alpha_{\cX,\cY}:\cBan(\cX,\cY)\lra\gC(\cX,\cY) $$ 
such that $\germ_{0}(\alpha(\gm))=\gm$
for any $\gm\in\cBan(\cX,\cY)$. The map
$\alpha$ may not be a functor, but if
$\gm\in\cBan(\cX,\cY)$ and $\gm'\in\cBan(\cY,\cZ)$,
the functions $\alpha(\gm'\com \gm)$ and $\alpha(\gm')\com\alpha(\gm)$
determine the same germ at $0$,
and hence they coincide on some neighbourhood of $0$.

Using $\alpha$, the seminorm $\|\oo\|_{r}$ induces a
seminorm on $\cBan$, which will be denoted again by
$\|\oo\|_{r}$, for all $r\geq 0$.
But $\|\oo\|_{r}$ is in fact a norm on $\cBan$ for any $r>0$.
Moreover, the inequality (\ref{presentability}) gives us 
\begin{eqnarray*} 
   \lim\sup_{\!\!\!\!\!\!\!\!\!\!\!\! r\ra 0}
   \frac{\|\gm\com \gm'-\gm\com \gm''\|_{r}}{\|\gm'-\gm''\|_{r}} &=&
   \lim\sup_{\!\!\!\!\!\!\!\!\!\!\!\! r\ra 0}
   \frac{\|\alpha(\gm)\com\alpha(\gm')-
   \alpha(\gm)\com\alpha(\gm'')\|_{r}}
   {\|\alpha(\gm')-\alpha(\gm'')\|_{r}} \\
   &\leq& \|D(\alpha(\gm))\|=\|\cD(\gm)\|
\end{eqnarray*}
for any $\gm\in\cBan(\cY,\cZ)$ and any two distinct
$\gm',\gm''\in\cBan(\cX,\cY)$.
It follows that $\cD$ is a $\Ban$-differential, 
so $(\cBan,\cD,\Ban)$ is a differential category.
\end{ex}

\subsection{Representations of Groups in Differential Categories}
\label{subsecRepGro}

If $\GG$ is a groupoid and $\gC$ a category, we call a functor
$$ R:\GG\lra\gC $$
a {\em representation} of $\GG$ in $\gC$.
Such a representation $R$ is called {\em trivial} if the
restriction $R|_{\GG(a,a')}$
is a constant map, for any $a,a'\in\GG_{0}$.
\vspace{4 mm}

\Rem
If $\cX$ is a vector space, the underlying abelian group $\cX$
can be seen as a category with only one object,
so we can speak about $\cX$-representations of a groupoid $\GG$.
If $G$ is a group and $\gC$ a category, a representation
of $G$ in $\gC$ is thus just a homomorphism of groups
$R:G\ra\gC^{-1}(a,a)$, where $\gC^{-1}(a,a)$ is the group
of isomorphisms in $\gC(a,a)$, for some
$a\in\Ob\gC$. For example, if $\cX$ is a Banach space,
then $\Ban(\cX,\cX)$ is a Banach category, and
the representations of $G$ in $\Ban(\cX,\cX)$ are
just the representations of $G$ in the classical sense.

\begin{ex}  \label{IIIex15}  \rm
Let $\cE$ be a Banach space over $\RRR$ and $M$ a $\eCe$-manifold
modeled on $\cE$. For any $x\in M$,
choose a chart $\varphi_{x}:U_{x}\ra\cE$ on $M$ with
$x\in U_{x}$ and $\varphi_{x}(x)=0$.
Define a representation $\rho=\rho_{M}$ of the groupoid
$\Gm_{\eCe}(M)$ in the category $\cBan$
(Example \ref{IIIex14}) by
$$ \rho(\germ_{x}f)=\germ_{0}
   [\varphi_{f(x)}\com f\com\varphi_{x}^{-1}]
   \in\cBan(\cE,\cE)\;,$$
for any $\eCe$-diffeomorphism $f$ on $M$ defined on an
open neighbourhood of $x$. The representation $\rho$ is
clearly well-defined and faithful. Another choice of charts
gives a naturally isomorphic representation.
\end{ex}

\Not
Let $G$ be a finitely generated group, $B$ a finite basis
(i.e. a finite set of generators) for $G$
with $1\in B=B^{-1}$, and let $\cX$ be a normed space. 
Further, let $P$ be a subset of $G$ with $B\subset P$ and
let $\varepsilon\geq 0$. A
$(P,\varepsilon)$-{\em approximate representation} in $\cX$
is a function 
$$ \eta:P\lra\cX $$
such that for any $g,g'\in P$ with $g'g\in P$ we have
$$ \|\eta(g'g)-\eta(g')-\eta(g)\|\leq\varepsilon\;.$$
We say that $\eta$ is {\em normed} if 
$$ \max_{g\in B}\|\eta(g)\|=1\;.$$
Note that a $(G,0)$-approximate representation in $\cX$ is
exactly a representation of $G$ in $\cX$.
\vspace{4 mm} \\
Recall now the following result of Thurston \cite[Lemma 1]{Thu}:

\begin{lem}  \label{IIIlem16}
Let $G$ be a finitely generated group and let $B$ be a finite
basis for $G$ with $1\in B=B^{-1}$.
Then non-trivial representations of $G$ in $\FFF$ exist if
and only if normed $(B^{l},\varepsilon)$-approximate
representations in $\FFF$ exist for any $\varepsilon >0$
and any integer $l\geq 1$.
\end{lem}

\begin{theo}  \label{IIItheo17}
Let $G$ be a finitely generated group, $(\gC,D,\gB)$ a
differential category over $\FFF$ and $R$ a representation
of $G$ in $\gC$. Then either
\begin{enumerate}
\item [(i)]   $R$ is trivial, or
\item [(ii)]  $D\com R$ is non-trivial, or
\item [(iii)] there exists a non-trivial representation of $G$
              in $\FFF$.
\end{enumerate}
\end{theo}
\Proof
Assume that $R$ is non-trivial and $D\com R$ is trivial.
Let $a$ be the object of $\gC$ such that $R:G\ra\gC^{-1}(a,a)$. 
Choose a finite basis $B$ for $G$ with $1\in B=B^{-1}$.
Since $D$ is a $\gB$-differential, there exists a sequence of
norms
$$ (\|\oo\|_{n})_{n=1}^{\infty} $$
on $\gC$  as in Definition \ref{IIIdfn12}.
By Lemma \ref{IIIlem16} it is enough to show that for any
$\varepsilon >0$ and any integer $l\geq 1$,
there exists a normed $(B^{l},\varepsilon)$-approximate
representation in $\FFF$. Take therefore $\varepsilon >0$
and $l\geq 1$, and let $\delta>0$ be so small that
$((l-1)\delta +l)\delta\leq\varepsilon$.

If $g,g'\in G$ and $R(g)\neq 1_{a}$,
we have $D(R(g')-1_{a})=0_{D a,D a}$, and hence
$$ \lim_{n\ra\infty}\frac{\|(R(g')-1_{a})\com
   R(g)-(R(g')-1_{a})\|_{n}}{\|R(g)-1_{a}\|_{n}}=0\;.$$
Since $B^{l}$ is finite, we can choose $n\geq 1$ such that
$$ \|(R(g')-1_{a})\com R(g)-(R(g')-1_{a})\|_{n}\leq
   \|R(g)-1_{a}\|_{n} \delta $$
for any $g,g'\in B^{l}$. 

Denote $M=\max_{g\in B}\|R(g)-1_{a}\|_{n}$.
Since $R$ is non-trivial, we have $M>0$. Define now
$\eta:B^{l}\lra\gC(a,a)$ by
$$ \eta(g')=\frac{1}{M}(R(g')-1_{a})
   \;\;\;\;\;\;\;\;\;\; g'\in B^{l}\;.$$

First we shall prove by induction that for any integer $k$,
$1\leq k\leq l$,
\begin{eqnarray}  \label{est1}
   \|\eta(g')\|_{n}\leq (k-1)\delta +k
\end{eqnarray}
for any $g'\in B^{k}$. 
This is clearly true for $k=1$. Assume now that the inequality
(\ref{est1}) holds for some $k<l$. Take any 
$g''\in B^{k+1}$ and write $r=g'g$ for some $g'\in B^{k}$ and
$g\in B$. Then we have
\begin{eqnarray*}
\|\eta(g'')\|_{n} &=& \frac{1}{M}\|R(g'g)-1_{a}\|_{n}       \\
&=& \frac{1}{M}\|R(g'g)-R(g)-R(g')+1_{a}
    +R(g')-1_{a}+R(g)-1_{a}\|_{n}                           \\
&\leq&
    \frac{1}{M}\|(R(g')-1_{a})\com R(g)-(R(g')-1_{a})\|_{n} \\
& & +\frac{1}{M}\|R(g')-1_{a}\|_{n}
    +\frac{1}{M}\|R(g)-1_{a}\|_{n}                          \\
&\leq&
    \frac{1}{M}\|R(g)-1_{a}\|_{n}\delta+
    \frac{1}{M}\|R(g')-1_{a}\|_{n}+  
    \frac{1}{M}\|R(g)-1_{a}\|_{n}                           \\ 
&=& \|\eta(g)\|_{n}\delta+
    \|\eta(g')\|_{n} + \|\eta(g)\|_{n}                      \\
&\leq&
    \delta +(k-1)\delta+k+1 = k\delta +(k+1)\;.
\end{eqnarray*}
This proves the inequality (\ref{est1}).
In particular we have for any $g'\in B^{l}$
$$ \|\eta(g')\|_{n}\leq (l-1)\delta +l\;.$$

Now take any $g,g'\in B^{l}$ such that $g'g\in B^{l}$.
The inequality above yields
\begin{eqnarray*}
\lefteqn{\|\eta(g'g)-\eta(g')-\eta(g)\|_{n}} \hspace{21 mm}\\
&=&
   \frac{1}{M}\|R(g'g)-1_{a}-R(g')+1_{a}-R(g)+1_{a}\|_{n}  \\
&=&
   \frac{1}{M}\|(R(g')-1_{a})\com R(g)-(R(g')-1_{a})\|_{n} \\
&\leq&
   \frac{1}{M}\|R(g)-1_{a}\|_{n}\delta
   =\|\eta(g)\|_{n}\delta                                  \\
&\leq&
   ((l-1)\delta+l)\delta\leq\varepsilon \;.
\end{eqnarray*}
Thus $\eta$ is a $(B^{l},\varepsilon)$-approximate
representation in $(\gC(a,a),\|\oo\|_{n})$.
It is clearly normed. Choose now a bounded linear functional
$$ f:(\gC(a,a),\|\oo\|_{n})\lra\FFF $$
of norm $1$ such that $f(\eta(g))=1$ for an element $g\in B$
with $\|R(g)-1_{a}\|_{n}=M$. Then we take
$\bar{\eta}=f\com\eta$. It follows that
$\bar{\eta}$ is a normed $(B^{l},\varepsilon)$-approximate 
representation in $\FFF$.
\eop

\subsection{Differential Categories Associated to Affine Manifolds}
\label{subsecDifCatAss}

In this Subsection we show that the germs of
$\eCe$-maps of an affine manifold \cite{Ehr,Kui}
provide a natural example of a differential category. 

Let $\cE$ be a fixed Banach space over $\FFF$. 
Recall that if $\cX$ and $\cY$ are Banach spaces and $U$
an open subset of $\cX$, a mapping
$f:U\ra \cY$ is {\em locally affine} if for any
$x\in U$ there exists a bounded linear operator
$A\in\Ban(\cX,\cY)$ such that
$$ f(x+\xi)=f(x)+A\xi $$
for any $\xi\in \cX$ small enough. In particular, the
map $f$ is analytic
(real analytic when $\FFF=\RRR$) and $A$ is the
differential $df_{x}$ of $f$ at $x$.

\begin{dfn}  \label{IIIdfn18}
An affine $\cE$-modeled atlas on a set $M$ is an
$\cE$-modeled atlas
$$ (\varphi_{i}:U_{i}\lra\cE)_{i\in I} $$
on $M$ such that the change of coordinates maps
$$ \varphi_{ij}=\varphi_{i}\com
   \varphi_{j}^{-1}|_{\varphi_{j}(U_{i}\cap \,U_{j})}:
   \varphi_{j}(U_{i}\cap \,U_{j})\lra\varphi_{i}
   (U_{i}\cap \,U_{j}) $$
is locally affine, for any $i,j\in I$.
An affine manifold modeled on $\cE$ is a set $M$ together
with a maximal affine $\cE$-modeled atlas on $M$.
\end{dfn}
\Rem
In particular, any affine manifold is an analytic 
(real analytic when $\FFF=\RRR$) manifold.
If $f:M\ra M'$ is a map between affine
manifolds, we say that
$f$ is {\em affine} 
if for any $x\in M$ there exist affine charts $\varphi:U\ra\cE$ for
$M$ and $\varphi':U'\ra\cE'$ for $M'$ with $x\in U$ and
$f(U)\subset U'$ such that
$$ \varphi'\com f\com\varphi^{-1}|_{\varphi(U)}:
   \varphi(U)\lra\varphi'(U') $$
is locally affine. One sees immediately that this condition is
in fact independent on the choice of charts.

\begin{ex}  \label{IIIex19}  \rm
It is easy to see that $S^{1}$, $T^{2}$  and the
M\"{o}bius strip can be equipped with
an affine structure. Of course, any open subset of $\cE$
is an affine manifold modeled on $\cE$.
A disjoint union or a product
of affine manifolds is again an affine manifold.
\end{ex}
\Not
Let $M$ be a $\eCe$-manifold modeled on $\cE$.
If $x,x'\in M$, we denote by
$\Delta(M)(x,x')$ the set of germs at $x$ of
$\eCe$-maps $f:U\ra M$, defined on
an open neighbourhood $U$ of $x$, with $f(x)=x'$. Denote
$$ \Delta(M)=\!\!\!\coprod_{x,x'\in M}\!\!\!\Delta(M)(x,x')\;.$$
Clearly $\Delta(M)$ has a structure of a category, given by
the composition of germs. There is a natural
$\eCe$-structure on
$\Delta(M)$ such that $\dom$ is a local
$\eCe$-diffeomorphism.
 
\begin{dfn}  \label{IIIdfn20}
Let $M$ be an affine manifold modeled on $\cE$,
and $x,x'\in M$. Define a structure of 
a vector space on the set $\Delta(M)(x,x')$ by
\begin{enumerate}
\item [(a)] $\germ_{x}f+\germ_{x}f'=\germ_{x}\left[\varphi^{-1}
            \!\!\com (\varphi\com f
            +\varphi\com f'-\varphi(x'))\right]\;,$
\item [(b)] $\lambda\,\germ_{x}f=\germ_{x}\left[\varphi^{-1}
            \!\!\com (\lambda\,\varphi\com f+(1-\lambda)\,
            \varphi(x'))\right]\;,$
\end{enumerate}
for any $\germ_{x}f,\,\germ_{x}f'\in\Delta(M)(x,x')$,
any $\lambda\in\FFF$ and
any affine chart $\varphi:U\ra\cE$ of $M$ with $x'\in U$.
\end{dfn}
\Rem
Note that with a fixed chart, this is well defined.
We have to check that it is also independent on the
choice of a chart. Let $\varphi':U'\ra\cE$ be another
affine chart of $M$ with $x'\in U'$. Since $M$ is affine, 
$\varphi\com\varphi'^{-1}$ is of the form
$A+\varphi(x')-A\varphi'(x')$ on a neighbourhood of
$\varphi'(x')$ for some bounded linear operator
$A\in\Ban(\cE,\cE)$. Therefore we have
on a small neighbourhood of $x$
\begin{eqnarray*}
\lefteqn{\varphi^{-1}\!\!\!\com (\varphi\com f+\varphi
         \com f'-\varphi(x'))} \hspace{6 mm}                   \\
&=&
   \varphi^{-1}(\varphi\varphi'^{-1}\varphi' f
   +\varphi\varphi'^{-1}\varphi' f'-\varphi(x'))               \\
&=&
   \varphi^{-1}(A\varphi' f+\varphi(x')-A\varphi'(x')+A\varphi' 
   f'+\varphi(x')-A\varphi'(x')-\varphi(x'))                   \\
&=&
   \varphi^{-1}\varphi\varphi'^{-1}(\varphi' f+
   \varphi' f'-\varphi'(x'))=
   \varphi'^{-1}\!\!\!\com (\varphi'\com f+\varphi'\com f'
   -\varphi'(x')) \;,
\end{eqnarray*}
so the definition of $\germ_{x}f+\germ_{x}f'$ is indeed
independent on the choice of a chart. In a similar way we can
check this also for (b).

With a fixed chart, it is trivial to see that this definition
gives a linear structure on $\Delta(M)(x,x')$. The neutral
element in $\Delta(M)(x,x')$ is the germ of the constant map.

\begin{obs}  \label{IIIobs21}
Let $M$ be an affine manifold modeled on $\cE$. 
Then $\Delta(M)$ is a right linear category.
\end{obs}
\Not
Let $M$ be an affine manifold, and $x\in M$. Recall that the
tangent space $T_{x}M$ can be described as the space of
equivalence classes of triples
$(U,\varphi,e)$, where $\varphi:U\ra\cE$ is an affine
chart of $M$, $x\in U$ and $e\in \cE$, and 
where $(U,\varphi,e)$ is equivalent to $(U',\varphi',e')$ if
$$ d(\varphi'\com\varphi^{-1})_{\varphi(x)}e=e'\;.$$
A chosen chart $\varphi:U\ra\cE$ determines an isomorphism
$\varphi_{\ast\, x}:T_{x}M\ra\cE$
and hence a norm on $T_{x}M$.
Another chart gives an equivalent norm,
so $T_{x}M$ is a Banachable space \cite{Lan}.

Now choose one of the equivalent norms $\|\oo\|_{x}$
on $T_{x}M$, for any $x\in M$. If $M$ is paracompact
(and hence separated),
it admits continuous partitions of
unity and in this case we can choose these norms such that
$$ \|\oo\|=\coprod_{x\in M}\|\oo\|_{x}:TM\lra [0,\infty) $$
is continuous.

Any element $\germ_{x}f\in\Delta(M)(x,x')$ has
the derivative
$$ D(\germ_{x}f)=f_{\ast\, x}\in\Ban(T_{x}M,T_{x'}M)\;.$$
With $D(x)=T_{x}M$ for any $x\in M$,
$D$ is a $\Ban$-linearization of $\Delta(M)$,
$$ D:\Delta(M)\lra\Ban\;.$$
Note that another choice of norms on the Banachable tangent
spaces gives a naturally isomorphic linearization.

\begin{theo}  \label{IIItheo22}
Let $M$ be an affine manifold modeled on $\cE$. 
Then
$$ (\Delta(M),D,\Ban) $$
is a differential category, for any choice of norms
on the Banachable tangent spaces of $M$.
\end{theo}
\Proof
Choose a norm on the Banachable tangent space $T_{x}M$,
for any $x\in M$. We have to prove that the corresponding
linearization $D$ is a $\Ban$-differential.

Let $\germ_{x}f\in\Delta(M)(x,x')$, where $f:W\ra M$ is a
$\eCe$-map defined on an open neighbourhood $W$ of $x$.
Let $\varphi:U\ra\cE$ and $\vartheta:V\ra\cE$ be affine
charts of $M$ with
$x\in U$ and $x'\in V$ such that
$\varphi(x)=0$ and $\vartheta(x')=0$.
We can assume that $W\subset U$ and
$f(W)\subset V$. Now we can define
$s(f):\varphi_{\ast\, x}^{-1}(\varphi(W))\ra T_{x}M$ by
$$ s(f)=\vartheta_{\ast\, x'}^{-1}\com\vartheta\com f
   \com\varphi^{-1}\com
   \varphi_{\ast\, x}|_{\varphi_{\ast\, x}^{-1}
   (\varphi(W))}\;.$$
The map $s(f)$ is clearly of class $\eCe$, and 
$0\in\varphi_{\ast\, x}^{-1}(\varphi(W))$ with $s(f)(0)=0$.
Define $F(\germ_{x}f)\in\cBan(T_{x}M,T_{x'}M)$ by
$$ F(\germ_{x}f)=\germ_{0}s(f)\;.$$
Note first that this is well-defined on germs.
Secondly, this definition is 
independent on the choice of charts. To prove this, let 
$\varphi':U'\ra\cE$ and $\vartheta':V'\ra\cE$ be another charts with
$x\in U'$, $x'\in V'$, $\varphi'(x)=0$ and $\vartheta'(x')=0$.
We can shrink $W$ further such that
$W\subset U\cap U'$ and $f(W)\subset V\cap V'$.
Now on a small neighbourhood
of $0\in T_{x}M$ we have
$$ \vartheta_{\ast\, x'}^{-1}\com\vartheta\com 
   f\com\varphi^{-1}\com\varphi_{\ast\, x}=
   \vartheta_{\ast\, x'}^{-1}\com\vartheta\com
   \vartheta'^{-1}\com\vartheta'\com f\com
   \varphi'^{-1}\com\varphi'\com\varphi^{-1}\com
   \varphi_{\ast\, x}\;. $$
But since $M$ is affine, we have $\varphi'\com\varphi^{-1}=d(\varphi'\com\varphi^{-1})_{0}=
\varphi'_{\ast\, x}\com\varphi_{\ast\, x}^{-1}$, and also
$\vartheta'\com\vartheta^{-1}=d(\vartheta'\com
\vartheta^{-1})_{0}=\vartheta'_{\ast\, x'}\com
\vartheta_{\ast\, x'}^{-1}$.
Therefore
$$ \vartheta_{\ast\, x'}^{-1}\com\vartheta\com 
   f\com\varphi^{-1}\com\varphi_{\ast\, x}=
   \vartheta'^{-1}_{\ast\, x'}\com\vartheta'\com 
   f\com\varphi'^{-1}\com\varphi'_{\ast\, x} $$
on a small neighbourhood of $0$.
This means exactly that the definition of $F$ is independent
on the choice of charts. Extending $F$ on objects of
$\Delta(M)$ by $F(x)=T_{x}M$ for any $x\in M$, we get a functor
$$ F:\Delta(M)\lra\cBan\;.$$
It is easy to check that $F$ is linear and faithful.
In example \ref{IIIex14}
we defined the differential $\cD:\cBan\ra\Ban$.
But for any $\germ_{x}f\in\Delta(M)(x,x')$ we have
\begin{eqnarray*}
\cD(F(\germ_{x}f)) &=&
    \cD(\germ_{0}(\vartheta_{\ast\, x'}^{-1}\com\vartheta
    \com f\com\varphi^{-1}\com\varphi_{\ast\, x}))           \\
&=& d(\vartheta_{\ast\, x'}^{-1}\com\vartheta\com f
    \com\varphi^{-1}\com\varphi_{\ast\, x})_{0}              \\ 
&=& \vartheta_{\ast\, x'}^{-1}\com\vartheta_{\ast\, x'}
    \com f_{\ast\, x}\com\varphi_{\ast\, x}^{-1}\com
    \varphi_{\ast\, x}=f_{\ast\, x}=D(\germ_{x}f)\;,
\end{eqnarray*}
for any pair of charts $\varphi$ and $\vartheta$
as before. Therefore we have
$$ \cD\com F=D\;.$$
Now Lemma \ref{IIIlem13} implies that $D$ is a
$\Ban$-differential.
\eop

\section{Reeb-Thurston Stability for Bibundles}
\label{secReeThuSta}

In this Section we use Theorem \ref{IIItheo17} to generalize the
Reeb-Thurston stability to transitive $\GG$-$\HH$-bibundles.
First, we define the notion of linear holonomy of
a leaf of a transitive bibundle.

For this, we assume that $\GG$ is an \'{e}tale $\eCe$-groupoid
modeled on a Banach space $\cE$ over $\RRR$, i.e.
both $\GG_{1}$ and $\GG_{0}$ are modeled on $\cE$. The effect-functor
$\eee$ now maps $\GG$ in the \'{e}tale $\eCe$-groupoid
$\Gm_{\eCe}(\GG_{0})$,
$$ \eee:\GG\lra\Gm_{\eCe}(\GG_{0})\;,$$
and it is of class $\eCe$.
In particular, for each $g\in\GG_{1}$ we
can compute the
differential $d\eee(g)$ of the germ $\eee(g)$,
which is a continuous linear map between the Banachable
tangent spaces of $\GG_{0}$
$$ d\eee(g)=\eee(g)_{\ast\, \dom g}:
   T_{\dom g}\GG_{0}\lra T_{\cod g}\GG_{0}\;.$$
Thus $\GG$ acts continuously on the tangent bundle 
$T(\GG_{0})$ of $\GG_{0}$ with respect to the projection
$T(\GG_{0})\ra\GG_{0}$ by
$$ g\cdot\xi=d\eee(g)(\xi)\;,$$
for any $g\in\GG_{1}$ and any $\xi\in T_{\dom g}\GG_{0}$.

Choose a norm on the Banachable tangent space
$T_{a}\GG_{0}$, for any $a\in\GG_{0}$. Then
$$ d\eee:\GG\lra\Ban $$
is a representation of $\GG$ in $\Ban$.

Let $\HH$ be a topological groupoid and $(E,p,\xW)$ a
transitive $\GG$-$\HH$-bibundle. If $L$ is a leaf of
$E$, we define the {\em linear holonomy group} $d\cH(L)$
of $L$ to be the image of $\cH(L)$ with $d\eee$.
It is defined uniquely up to a conjugation.

Assume now that $(E,p,\xW)$ is principal with locally
path-connected fibers. Take $e\in E$ and let $L$ be the
leaf of $E$ with $\xW(e)\in L$. The {\em linear holonomy
homomorphism} of $L$ with respect to the base point $e$
is the composition
$$ d\cH_{e}=d\eee\com\cH_{e}:\pi_{1}(\HH(L),\xW(e))\lra
   \Ban(T_{p(e)}\GG_{0},T_{p(e)}\GG_{0})\;.$$
In the same way as the holonomy homomorphism, the linear
holonomy homomorphism depends up to an isomorphism only on
the leaf $L$. Therefore we can define
(up to an isomorphism) the {\em linear holonomy homomorphism}
$d\cH_{L}:\pi_{1}(\HH(L))\ra\Ban$ of $L$. Clearly, the linear
holonomy group $d\cH(L)$ is the image of $d\cH_{L}$.

If $\GG$ is finite-dimensional, we can define (by analogy
with foliations on manifolds) that $E$ is transversely
orientable around $L$ if all the elements in
the linear holonomy group preserve orientation, i.e.
they have positive determinant.

\begin{prop}  \label{IIIprop23}
Let $\GG$ be an effective $\eCe$-gro\-up\-oid modeled on a
Banach space $\cE$,
$\HH$ a topological gro\-up\-oid and $(E,p,\xW)$ a
transitive $\GG$-$\HH$-bibundle.
Let $L$ be a leaf of $(E,p,\xW)$ such that
$\Ker d\eee|_{\cH(L)}$
is finitely generated and
$$ \Hom(\Ker d\eee|_{\cH(L)},\RRR)=\{0\}\;.$$
Then $\Ker d\eee|_{\cH(L)}$ is trivial.
\end{prop}
\Proof
A choice of charts as in example \ref{IIIex15}
gives a representation $\rho$ of $\Gm_{\eCe}(\GG_{0})$
in $\cBan$, and the composition
$$ \rho\com\eee:\GG\lra\cBan $$
is a faithful representation of $\GG$ in $\cBan$.
(Example \ref{IIIex14}). The composition
$\cD\com\rho\com\eee$ is isomorphic to
$d\eee$. The restriction $R$ of the representation
$\rho\com\eee$ to the group
$\Ker d\eee|_{\cH(L)}$ is now a
representation of $\Ker d\eee|_{\cH(L)}$, and
the composition $\cD\com R$ is trivial.
Now Theorem \ref{IIItheo17} yields that $R$ is trivial.
But $\rho\com\eee$ is faithful, therefore
the group $\Ker d\eee|_{\cH(L)}$ is trivial.
\eop

\begin{cor}  \label{IIIcor24}
Let $\GG$ be an effective $\eCe$-groupoid modeled on a
Banach space $\cE$, $\HH$ a topological groupoid and 
$(E,p,\xW)$ a locally compact transitive $\GG$-$\HH$-bibundle
with locally connected $\HH$-Hausdorff fibers.
Let $L$ be an $\HH$-compact leaf of $(E,p,\xW)$
such that 
\begin{enumerate}
\item  $\Ker d\eee|_{\cH(L)}$ is finitely generated,
\item  $\Hom(\Ker d\eee|_{\cH(L)},\RRR)=\{0\}$, and
\item  the linear holonomy group of $L$ is finite.
\end{enumerate}
Then the holonomy group of $L$ is finite, and
for any $\HH$-invariant open neighbourhood $V$ of $L$
in $\HH_{0}$ there exists an $\HH$-invariant open
neighbourhood $U\subset V$ of $L$ in $\HH_{0}$
which is a union of $\HH$-compact leaves of $(E,p,\xW)$.
\end{cor}
\Proof
Proposition \ref{IIIprop23} implies that the holonomy of
$L$ is also finite, and hence the conditions of Theorem
\ref{IIItheo6} are satisfied.
\eop

\begin{theo} \label{IIItheo25}
Let $\GG$ be an effective $\eCe$-groupoid modeled on a
Banach space $\cE$, 
$\HH$ a topological groupoid with $\HH_{0}$ locally compact
and $(E,p,\xW)$ a principal $\GG$-$\HH$-bibundle
with locally path-connected $\HH$-Hausdorff fibers.
Let $L$ be an $\HH$-compact leaf of $(E,p,\xW)$
such that 
\begin{enumerate}
\item  $\Ker d\cH_{L}$ is finitely generated,
\item  $\Hom(\Ker d\cH_{L},\RRR)=\{0\}$, and
\item  the linear holonomy group of $L$ is finite.
\end{enumerate}
Then the holonomy group of $L$ is finite, and
for any $\HH$-invariant open neighbourhood $V$ of $L$ in
$\HH_{0}$ there exists an $\HH$-invariant open neighbourhood 
$U\subset V$ of $L$ in $\HH_{0}$
which is a union of $\HH$-compact leaves of $(E,p,\xW)$.
\end{theo}
\Proof
Since $\HH_{0}$ is locally compact and $\xW$ a local
homeomorphism, the space $E$ is also locally compact.
Note that for the surjective 
homomorphism $\cH_{L}:\pi_{1}(\HH(L))\ra\cH(L)$ we have
$$ \cH_{L}(\Ker d\cH_{L})=\Ker d\eee|_{\cH(L)}\;.$$
This implies that $\Ker d\eee|_{\cH(L)}$
in finitely generated and
$$ \Hom(\Ker d\eee|_{\cH(L)},\RRR)=\{0\}\;.$$
Proposition \ref{IIIprop23} now yields that
the holonomy of $L$ is also finite, 
and hence we can use Theorem \ref{IIItheo6}.
\eop
\Rem 
Theorem \ref{IIItheo25} generalizes the Reeb-Thurston
stability theorem for foliations on manifolds to the
Hilsum-Skandalis maps. 

\begin{cor} \label{IIIcor26}
Let $X$ be a locally compact space, $\GG$ an
\'{e}tale $\eCe$-groupoid modeled on a Banach space
$\cE$, and $(E,p,\xW)$ a principal $\GG$-bundle over $X$
with locally path-connected Hausdorff fibers.
Let $L$ be a compact leaf of $(E,p,\xW)$ such that
\begin{enumerate}
\item  $\Ker d\cH_{L}$ is finitely generated,
\item  $\Hom(\Ker d\cH_{L},\RRR)=\{0\}$, and
\item  the linear holonomy group of $L$ is finite.
\end{enumerate}
Then the geometric holonomy group of $L$ is finite, and
for any open neighbourhood $V$ of $L$ there exists an
open neighbourhood $U\subset V$ of $L$ which is a union
of compact leaves of $(E,p,\xW)$ with finite geometric holonomy
groups.
\end{cor}
\Proof
By Theorem \ref{IItheo16} and Theorem \ref{IItheo17},
we can assume without loss of generality that $\GG$ is
effective. The corollary now follows from Theorem
\ref{IIItheo25} and Theorem \ref{IIItheo7}.
\eop
\Rem
This is the Reeb-Thurston stability theorem for Haefliger
structures on topological spaces.

\begin{prop}  \label{IIIprop27}
Let $\GG$ be an \'{e}tale $\eCe$-groupoid of
dimension one, $\HH$ a topological groupoid and 
$(E,p,\xW)$ a principal $\GG$-$\HH$-bibundle
with locally path-connected fibers.
Let $L$ be a leaf of $(E,p,\xW)$ such that 
$$ \Hom(\pi_{1}(\HH(L)),\RRR)=\{0\}\;.$$
Then the linear holonomy group of $L$ has at most two
elements.
\end{prop}
\Proof 
Since $\GG$ is of dimension one, the linear holonomy
homomorphism $d\cH_{L}$ of $L$ maps $\pi_{1}(\HH(L))$ into
a group isomorphic to the multiplicative group
$\RRR^{\ast}=\RRR\setminus\{0\}$. Hence we can form the
homomorphism
$$ \log |d\cH_{L}|:\pi_{1}(\HH(L))\lra\RRR\;,$$
which must be trivial by assumption. Therefore the image of
$d\cH_{L}$ has at most two elements.
\eop

\begin{cor}  \label{IIIcor28}
Let $X$ be a locally compact space, $\GG$ an
\'{e}tale $\eCe$-groupoid of dimension one,
and $(E,p,\xW)$ a principal $\GG$-bundle over $X$
with locally path-connected Hausdorff fibers.
Let $L$ be a compact leaf of $(E,p,\xW)$ such that
\begin{enumerate}
\item  $E$ is transversely orientable around $L$,
\item  $\pi_{1}(H(L))$ is finitely generated, and
\item  $\Hom(\pi_{1}(H(L)),\RRR)=\{0\}$.
\end{enumerate}
Then the geometric holonomy group of $L$ is trivial, and
for any open neighbourhood $V$ of $L$ there exists an
open neighbourhood $U\subset V$ of $L$ which is a union
of compact leaves of $(E,p,\xW)$. Moreover, if
$L'$ is a leaf of $(E,p,\xW)$ with $L'\subset U$, the
geometric holonomy group of $L'$ is either trivial or
isomorphic to $\ZZZ/2\ZZZ$.
\end{cor}
\Proof
By Proposition \ref{IIIprop27}, the linear holonomy group
of $L$ has at most two elements, but because of the transversal
orientability around $L$ it is in fact trivial. Now apply
Corollary \ref{IIIcor26}. Finally, observe that if the
geometric holonomy group of a leaf of $E$ is finite then
it has at most two elements. Thus the assumption (1) implies
that the geometric holonomy group of $L$ is trivial.
\eop
\Rem
Note that if $X$ is Hausdorff, the assumption
that the fibers of $E$ are Hausdorff in
Corollary \ref{IIIcor26} and in Corollary \ref{IIIcor28},
is superfluous.

\chapter{Equivariant Foliations}
\label{chapEquFol}
\startchapterskip

As an application of the results of Chapters \ref{chapHilSkaMap}
and \ref{chapSta}, we study in this Chapter the $\eCe$-principal
$\GmqCe$-$G(M)$-bibundle associated to
a $G$-equivariant foliation on a $\eCe$-manifold $M$
(Example \ref{IIex8} (3)).
We give a geometric interpretation for the holonomy
of this bibundle, and show that the stability theorems take
a much simpler form in this case. Finally, we give some
examples of concrete equivariant foliations.

\section{Equivariant Holonomy}  \label{secEquHol}

Let $M$ be a $\eCe$-manifold of dimension $n$ and
$\nu:M\times G\ra M$ a right $\eCe$-action of a discrete group
$G$ on $M$. Let $\cF$ be a foliation of codimension $q$ on $M$
which is invariant under the action of $G$,
see Example \ref{IIex8} (3). For any $g\in G$
denote by $\hg$ the $\eCe$-diffeomorphism $\nu(\oo,g):M\ra M$.
As usual, we will write $\nu(x,g)=x\cdot g$.

Let $L$ be a leaf of $\cF$. The {\em isotropy group} $G_{L}$
of $L$ is the subgroup of $G$ 
$$ G_{L}=\{\,g\in G\,|\,L\cdot g\subset L\,\}\;.$$
Note that $L\cdot g\subset L$ implies $L\cdot g=L$.
The action of $G$ clearly restricts to a $\eCe$-action of
$G_{L}$ on $L$, where $L$ has the leaf topology.
The immersion of $L$ into $M$ induces a continuous injection of
the orbit space $L/G_{L}$ into the orbit space $M/G$.
As a set, we may identify $L/G_{L}$ with the corresponding
subset of $M/G$. However, the topology of $L/G_{L}$ is in general
finer than the one inherited from the space $M/G$.

Let $L'$ be another leaf of $\cF$. If there exists $g\in G$ such
that $L\cdot g=L'$, then clearly $G_{L'}=g^{-1}G_{L}\,g$, and
the $\eCe$-diffeomorphism
$$ \hg|_{L}:L\lra L' $$
induces the identity homeomorphisms between the orbit spaces
$L/G_{L}$ and $L'/G_{L'}$. On the other hand, if
$L\cdot G\cap L'=\emptyset$ then $L/G_{L}\cap L'/G_{L'}=\emptyset$.
Thus $\cF$ induces a partition of the orbit space $M/G$.

If the action of $G$ on $M$ is properly discontinuous, the orbit
space $M/G$ is a $\eCe$-manifold and the induced partition is the
induced $\eCe$-foliation $\cF/G$ on $M/G$. In general, the induced
partition of $M/G$ may be seen as a generalized foliation on $M/G$.

Let $L$ be a leaf of $\cF$ and $x_{0}\in L$.
We will denote by $\pi^{G}_{1}(L,x_{0})$ the group of
pairs
$$ \pi^{G}_{1}(L,x_{0})=\{\,(g,\vsig)\,|\,g\in G,\,
   \vsig\in\bpi_{1}(L)(x_{0},x_{0}\cdot g)\,\} $$
with multiplication given by
$$ (g',\vsig')(g,\vsig)=(g'g,\hg_{\ast}(\vsig')\,\vsig)\;.$$
We will call this group the {\em equivariant fundamental
group} of $L$. Recall from Proposition \ref{Iprop12} that the
equivariant fundamental group of $L$ is isomorphic to the
fundamental group $\pi_{1}(G_{L}(L),x_{0})$ of the \'{e}tale
$\eCe$-groupoid $G_{L}(L)$ associated to the action of $G_{L}$
on $L$, and that there is a short exact sequence
\begin{equation}  \label{IVshoexaseq}
\CD
  1 \cdr{}{}        \pi_{1}(L,x_{0})
  \cdr{\inc_{L}}{}  \pi_{1}^{G}(L,x_{0})
  \cdr{\pr_{L}}{}   G_{L}
  \cdr{}{} 1
\endCD
\end{equation}
of homomorphisms of groups, where $\inc_{L}(\vsig)=(1,\vsig)$
and $\pr_{L}(g,\vsig)=g$.

Note that if $x'_{0}$ is another point of $L$, then a path in $L$
between $x_{0}$ and $x'_{0}$ induces  an isomorphism
$\pi_{1}^{G}(L,x_{0})\cong\pi_{1}^{G}(L,x'_{0})$. More precisely,
if $\tau$ is a homotopy class (with fixed end-points) of a path
from $x_{0}$ to $x'_{0}$, then the isomorphism
$$ \tau^{\diamond}:\pi_{1}^{G}(L,x'_{0})\lra
   \pi_{1}^{G}(L,x_{0}) $$
is given by
$$ \tau^{\diamond}(g,\vsig)=
   (g,\hg_{\ast}(\tau)^{-1}\!\vsig\tau)\;.$$
Therefore
we will write $\pi_{1}^{G}(L)$ for the isomorphism class of the
equivariant fundamental group of $L$.
Note that if $G_{L}$ acts properly discontinuously on $L$,
then the covering projection $L\ra L/G_{L}$ induces an isomorphism
$\pi_{1}^{G}(L)\cong \pi_{1}(L/G_{L})$.

\begin{dfn}  \label{IVdfn1}
Let $M$ be a $\eCe$-manifold of dimension $n$ equipped
with a right $\eCe$-action of a discrete group $G$, and let $\cF$
be a $G$-equivariant foliation of codimension $q$ on $M$.
Let $L$ be a leaf of $\cF$ and let $T:\RRR^{q}\ra M$
be a transversal section of $\cF$ with $T(0)\in L$.
The equivariant holonomy homomorphism of $L$ with
respect to the transversal section $T$ is the homomorphism
$$ \Hol_{T}^{G}:\pi_{1}^{G}(L,T(0))\lra\Dif{q} $$
given by
$$ \Hol_{T}^{G}(g,\vsig)=\Hol_{\hg\com T,T}(\vsig) $$
for any $(g,\vsig)\in\pi_{1}^{G}(L,T(0))$.
\end{dfn}
\Rem
First observe that the equivariant holonomy homomorphism
is well-defined. In particular, the composition $\hg\com T$ is
a transversal section with $(\hg\com T)(0)=T(0)\cdot g$.
Note that the restriction of $\Hol^{G}_{T}$ to the fundamental
group $\pi_{1}(L,T(0))$ is the usual holonomy homomorphism 
$\Hol_{T}$ of $L$ with respect to $T$, i.e.
$$ \Hol_{T}^{G}\com \inc_{L}=\Hol_{T}\;.$$
Let $T'$ be another transversal section with $T'(0)\in L$, and
choose a homotopy class (with fixed end-points) $\tau$ of a path
in $L$ from $T(0)$ to $T'(0)$. Now we have
\begin{eqnarray*}
\Hol^{G}_{T'}(g,\vsig) &=&   \Hol_{\hg\com T',T'}(\vsig)     \\
&=&
   \Hol_{\hg\com T',\hg\com T}(\hg_{\ast}(\tau))\com
   \Hol_{\hg\com T,T}(\hg_{\ast}(\tau)^{-1}\vsig\tau)\com
   \Hol_{T,T'}(\tau^{-1})                                    \\
&=&
   \Hol_{T',T}(\tau)\com
   \Hol^{G}_{T}(\tau^{\diamond}(g,\vsig))\com
   \Hol_{T',T}(\tau)^{-1}
\end{eqnarray*}
for any $(g,\vsig)\in\pi_{1}^{G}(L,T'(0))$. In other words,
the homomorphisms $\Hol_{T'}^{G}$ and
$\Hol_{T}^{G}\com\tau^{\diamond}$ differ by the conjugation
by $\Hol_{T',T}(\tau)$ in $\Dif{q}$, i.e. we have a canonical
isomorphism $\Hol_{T'}^{G}\cong\Hol_{T}^{G}\com\tau^{\diamond}$.
Therefore we will write 
$$ \Hol_{L}^{G}:\pi_{1}^{G}(L)\lra\Dif{q} $$
for the isomorphism class of the equivariant holonomy homomorphism
of $L$. The group $\Hol_{L}^{G}(\pi_{1}^{G}(L))$ is
called the {\em equivariant holonomy group} of $L$.

Intuitively, for any $g\in G_{L}$ the germ
$$ \Hol_{T}^{G}(g,\vsig) $$
is the transversal part of the diffeomorphism $\hg$ around $L$,
but determined up to the choice of a homotopy class $\vsig$ of a path
from $T(0)$ to $T(0)\cdot g$ in $L$.
In fact, it is only the holonomy class
of $\vsig$ which matters.

The geometric picture is particularly clear
if the group $G_{L}$ has a fixed point in $L$, so we can choose
$T$ so that $T(0)$ is that fixed point.
In this case there is a splitting
$$ s:G_{L}\lra\pi_{1}^{G}(L,T(0)) $$
of the sequence (\ref{IVshoexaseq}) given by
$s(g)=(g,1)$, and $\pi_{1}^{G}(L,T(0))$ is a semi-direct
product of $G_{L}$ and $\pi_{1}(L,T(0))$.
In particular, the transversal part of $\hg$ may be seen
as the germ $\Hol_{T}^{G}(s(g))$.

In general, since $\pi_{1}(L,T(0))$ is a 
normal subgroup of $\pi_{1}^{G}(L,T(0))$, 
it follows that the holonomy group $\Hol_{T}(\pi_{1}(L,T(0)))$
of $L$ is a normal subgroup of the equivariant holonomy group
$\Hol^{G}_{T}(\pi_{1}^{G}(L,T(0)))$ of $L$. The quotient
$$ \Hol^{G}_{T}(\pi_{1}^{G}(L,T(0)))/
   \Hol_{T}(\pi_{1}(L,T(0))) $$
is thus the contribution of the action of $G$ to the equivariant
holonomy of $L$.

We say that the foliation $\cF$ is {\em equivariantly
transversely orientable} around a leaf $L$ of $\cF$ if all
the germs in $\Hol_{L}^{G}(\pi_{1}^{G}(L))$ preserve the
orientation of $\RRR^{q}$, i.e. if the equivariant holonomy
homomorphism $\Hol_{L}^{G}$ actually maps the equivariant
holonomy group of $L$ into the group $\Difp{q}$.
In particular, if $\cF$ is equivariantly transversely orientable
around $L$, then $\cF$ is transversely orientable around $L$.
Moreover, note that if $\cF$ is transversely orientable around $L$
and $g\in G_{L}$, then all the germs in
$\Hol^{G}_{L}(\pr_{L}^{-1}(g))$ either preserve or invert
the orientation of $\RRR^{q}$.
\vspace{4 mm}

Let $L$ be a leaf of $\cF$. Define the group 
$$ G_{L}^{\sharp}=\pr_{L}(\Ker\Hol_{L}^{G})\subset G_{L}\;.$$
Observe that this group does not depend on the choice of a transversal 
section through the leaf $L$. An element $g\in G_{L}$ belongs to 
$G_{L}^{\sharp}$ if and only if there exists a homotopy class
$\vsig$ of a path in $L$ from a point $x$ to $x\cdot g$ such that
$$ \Hol_{L}^{G}(g,\vsig)=1\;.$$
With $x$ fixed, the class $\vsig$ is determined uniquely up to 
holonomy. According to the previous remark, $G_{L}^{\sharp}$ may
be regarded as the subgroup of those elements of $G_{L}$ which
act transversely trivially around $L$.

The important fact about the group $G_{L}^{\sharp}$ is that it
canonically acts on the holonomy covering space $\hcL$
of the leaf $L$, i.e. on a covering space of $L$ which
corresponds to the kernel of the holonomy homomorphism of $L$.
Explicitly, if
$$ \zeta:\hcL\lra L $$
is the holonomy covering projection, we define
for any $g\in G_{L}^{\sharp}$ and $e\in\hcL$
$$ e\cdot g=\tilde{\sigma}(1)\;,$$
where $\sigma$ is a path in $L$ such that $\sigma(0)=\zeta(e)$,
$\sigma(1)=\sigma(0)\cdot g$ and
$$ \Hol_{L}^{G}(g,[\sigma])=1\;,$$
and where $\tilde{\sigma}$ is the unique lift of $\sigma$ in
$\hcL$ with $\tilde{\sigma}(0)=e$.
Here $[\sigma]$ denotes the homotopy class of the path $\sigma$.
Since such a path is uniquely determined up to holonomy,
the definition makes sense. It gives a right $\eCe$-action
of $G_{L}^{\sharp}$ on $\hcL$ such that $\zeta$ is
equivariant. Hence $\zeta$ induces a continuous map
$$ \bar{\zeta}:\hcL/G_{L}^{\sharp}\lra L/G_{L}\;.$$
We say that the action of $G$ is {\em leafwise separated} if
the orbit space $\hcL/G_{L}^{\sharp}$ is Hausdorff, for any leaf
$L$ of $\cF$.

\begin{ex}  \label{IVex2}  \rm
Let $M$ and $N$ be finite-dimensional $\eCe$-manifolds
(not necessarily Hausdorff), and let
$p:M\ra N$ be a $\eCe$-submersion.
The submersion $p$ induces a foliation $\cF$ on $M$. The leaves
of $\cF$ are the connected components of the fibers of $p$.

Assume that $M$ is equipped with a $\eCe$-action of a discrete group
$G$ such that
$$ p(x\cdot g)=p(x) $$
for any $x\in M$ and $g\in G$.
The foliation $\cF$ is clearly invariant under this action.
Such an equivariant foliation $\cF$ will be referred to as
{\em simple} equivariant foliation.
If $L$ is a leaf of $\cF$, it is clear that both the holonomy
group and also the equivariant holonomy group of $L$ are trivial.
Thus $\hcL\cong L$ and $G_{L}^{\sharp}=G_{L}$. In particular,
the action of $G$ is leafwise separated if and only if the orbit
spaces of the leaves of $\cF$ are Hausdorff.

Denote by $r:M\ra M/G$ the quotient projection. The submersion
$p$ induces a continuous map $p/G:M/G\ra N$ such that
$$ p=p/G\com r\;.$$
Let $y\in p(M)$. The restriction
$$ r|_{p^{-1}(y)}:p^{-1}(y)\lra(p/G)^{-1}(y) $$
is an open surjection. Since $p$ is a submersion, the fiber
$p^{-1}(y)$ is locally connected, and hence so is the fiber
$(p/G)^{-1}(y)$. This implies that if $L$ is a leaf of $\cF$
then $r(L)$ is a connected component of the
corresponding fiber of $p/G$. Moreover, the restriction
$r|_{L}:L\ra r(L)$ is open, hence
$$ r(L)\cong L/G_{L}\;.$$
In other words, the orbit spaces of the leaves of $\cF$ are
isomorphic with the connected components of the fibers of $p/G$.
\end{ex}

Now let $(\varphi_{i}:U_{i}\ra\RRR^{n-q}\times\RRR^{q})_{i\in I}$
be the maximal atlas for $\cF$.
If $i\in I$ and $g\in G$, denote by $i\cdot g$ the element of
$I$ such that 
$$ \varphi_{i}\com\hg^{-1}|_{U_{i}\cdot g}=\varphi_{i\cdot g}\;,$$
as in Example \ref{IIex8} (3). Also denote
$s_{i}=\pr_{2}\com\varphi_{i}$, and let $c=(c_{ij})$ be the
$\GmqCe$-cocycle on $(U_{i})$ corresponding
to the family of submersions $(s_{i})$ (see Section \ref{secFol}).
Let $(\Sigma(c),p,\xW)$ be the principal $\GmqCe$-bundle over
$M$ associated to this cocycle. Write $E=\Sigma(c)$.
An element of $E$ is of the form
$[\gm,x,i]$, where $\gm\in\GmqCe$, $i\in I$
and $x\in U_{i}$ with $\dom\gm=s_{i}(x)$. The map $p:E\ra\RRR^{q}$
is a submersion. In particular, the fibers of $p$, i.e. the
fibers if $E$, are locally path-connected.

Recall from Example \ref{IIex8} (3) that $G$ acts on $E$ by
$$ [\gm,x,i]\cdot g=[\gm,x\cdot g, i\cdot g]\;.$$
This action induces an action of the \'{e}tale groupoid $G(M)$
on $E$, and $E$ becomes a $\eCe$-principal
$\GmqCe$-$G(M)$-bibundle.
We shall now investigate the correspondence between
the equivariant holonomy groups of leaves of $\cF$ and the
holonomy groups of leaves of the $\GmqCe$-$G(M)$-bibundle
$E$.

First note that since $p$ is a submersion, it defines a foliation
$\cFsh$ on $E$. Clearly we have
$$ p(e\cdot g)=p(e) $$
for any $e\in E$ and $g\in G$, i.e. $\cFsh$ is a simple
$G$-equivariant foliation on $E$ (Example \ref{IVex2}).
Since the underlying principal $\GmqCe$-bundle
of $E$ is just the one which represents the foliation $\cF$,
the map $\xW$ maps the leaves of $\cFsh$ onto the leaves of $\cF$.
Moreover, the map $\xW$ restricted to a leaf of $\cFsh$ is
the holonomy covering projection onto the image leaf of $\cF$.
Thus we proved:

\begin{prop}  \label{IVprop3}
Let $M$ be a $\eCe$-manifold of dimension $n$ equipped with
a $\eCe$-action of a discrete group $G$, and let $\cF$ be a
$G$-equivariant foliation of dimension $q$
on $M$. There exists a $\eCe$-manifold
$M^{\sharp}$ equipped with a $\eCe$-action of $G$, a simple
$G$-equivariant foliation $\cFsh$ on $M^{\sharp}$ given by a
submersion $p:M^{\sharp}\ra\RRR^{q}$, and a
$G$-equivariant local $\eCe$-diffeomorphism $\xW:M^{\sharp}\ra M$
which maps the leaves of $\cFsh$ onto the leaves of $\cF$
as the holonomy covering projection.
\end{prop}

Let $y\in\RRR^{q}$. A $G(M)$-connected component of the fiber
$E_{y}$ is of course a minimal $G$-invariant union of
leaves of $\cFsh$. Since the map $\xW$ is $G$-equivariant, a
leaf of the $\GmqCe$-$G(M)$-bibundle $E$ is a minimal
$G$-invariant union of leaves of $\cF$, equipped with
the action of $G$. Each leaf of $\cF$ is open and embedded
in the corresponding leaf of $E$. Therefore a leaf of $E$ is
exactly the disjoint union of those leaves of $\cF$ which
have the same orbit space in $M/G$. Observe that the leaves
of the associated transitive $\GmqCe$-bundle over $M/G=|G(M)|$
are exactly the orbit spaces of leaves of $\cF$ in $M/G$.

Choose $x_{0}\in M$ and $i\in I$
such that $x_{0}\in U_{i}$ and $\varphi_{i}(x_{0})=0$.
In particular, we have $s_{i}(x_{0})=0$.
We can assume without loss of generality that $\varphi_{i}$ is
surjective.
Let $T:\RRR^{q}\ra M$ be the transversal section of $\cF$ given by
$$ T(y)=\varphi_{i}^{-1}(0,y)\;.$$
Choose $e\in E$ such that $\xW(e)=x_{0}$. Hence
$$ e=[\gm,x_{0},i] $$
for a unique $\gm\in\GmqCe$. Let $\hcL$ be the leaf
of $\cFsh$ with $e\in\hcL$. Hence $\hcL\cdot G$ is the corresponding
$G(M)$-connected component of the fiber $E_{p(e)}$.
Let $L$ be the leaf of $\cF$ with $x_{0}\in L$. The restriction
$$ \zeta=\xW|_{\hcL}:\hcL\lra L $$
is the holonomy covering projection.
Denote by $L\cdot G$ the leaf of $E$ corresponding to $\hcL\cdot G$.
Observe that
$$ \pi_{1}(G(M)(L\cdot G),x_{0})\cong\pi_{1}(G_{L}(L),x_{0})
   \cong\pi_{1}^{G}(L,x_{0})\;.$$
Denote this isomorphism by
$\Phi:\pi_{1}^{G}(L,x_{0})\ra\pi_{1}(G(M)(L\cdot G),x_{0})$.
We will show that for any $(g,\vsig)\in\pi_{1}^{G}(L,x_{0})$
we have
$$ \cH_{e}(\Phi(g,\vsig))=
   \gm\com\Hol^{G}_{T}(g,\vsig)\com\gm^{-1}\;.$$
In other words, the equivariant holonomy homomorphism of $L$ and the
holonomy homomorphism of the corresponding leaf of $E$ differ by
the conjugation by $\gm$.

Let $(g,\vsig)\in\pi_{1}^{G}(L,x_{0})$, and let $\sigma$ be a path
in $L$ from $x_{0}$ to $x_{0}\cdot g$ which represents $\vsig$.
Let $\tilde{\sigma}$ be the lift of $\sigma$ along $\zeta$ such that
$\tilde{\sigma}(0)=e$. Observe that
$$ (\hg\com T)(y)=\varphi_{i\cdot g}^{-1}(0,y)\;.$$
This and the definition of holonomy imply that
\begin{eqnarray*}
\tilde{\sigma}(1) &=&
[\gm\com\Hol_{\hg\com T,T}(\vsig)^{-1},x_{0}\cdot g,i\cdot g]   \\
&=&
[\gm\com\Hol^{G}_{T}(g,\vsig)^{-1},x_{0},i]\cdot g\;.
\end{eqnarray*}
Therefore
\begin{equation}  \label{IVeq1}
  (\gm\com\Hol^{G}_{T}(g,\vsig)\com\gm^{-1})\cdot
  \tilde{\sigma}(1)\cdot g^{-1}=e\;.
\end{equation}
On the other hand, the element $\Phi(g,\vsig)$ is the
homotopy class of the $G(M)$-loop
$$ \varpi(x_{0})\cdot (x_{0},g)\cdot\sigma\in\Omega(G(M)
   (L\cdot G),x_{0})\;,$$
where $\varpi(x_{0})$ denotes the constant path with the image
point $x_{0}$. The lift of this $G(M)$-loop is clearly
the $G(M)$-path in $E$
$$ \varpi(\tilde{\sigma}(1)\cdot g^{-1})
   \cdot(x_{0},g)\cdot\tilde{\sigma}\;,$$
hence $\tilde{\sigma}(1)\cdot g^{-1}$ is the end-point of this lift.
Together with the equation (\ref{IVeq1}) this implies
$$ \cH_{e}([\varpi(x_{0})\cdot (x_{0},g)\cdot\sigma])
   =\gm\com\Hol^{G}_{T}(g,\vsig)\com\gm^{-1}\;.$$
Since another choice of $e$ and $T$ gives isomorphic holonomy
homomorphisms, we can conclude that
$$ \cH_{L}\com\Phi=\Hol^{G}_{L}\;.$$

Observe that the equation (\ref{IVeq1}) implies that
$$ G^{\sharp}_{L}=G_{\hcL}\;. $$
Indeed, if $g\in G^{\sharp}_{L}$, there exists a path $\sigma$
in $L$ from $x_{0}$ to $x_{0}\cdot g$ such that
$\Hol^{G}_{T}(g,[\sigma])=1$. Therefore by
the equation (\ref{IVeq1}) we have $e\cdot g\in\hcL$, hence 
$g\in G_{\hcL}$. Conversely, if $g\in G_{\hcL}$,
there exists a path $\tilde{\sigma}$ in $\hcL$ from $e$ to
$e\cdot g$, and thus the equation (\ref{IVeq1}) gives
$$ (\gm\com\Hol^{G}_{T}(g,[\zeta\com\tilde{\sigma}])\com\gm^{-1})
   \cdot e=e\;.$$
Since $E$ is principal, this implies that
$\Hol^{G}_{T}(g,[\zeta\com\tilde{\sigma}])=1$, and hence
$g\in G^{\sharp}_{L}$.

The equation (\ref{IVeq1}) also shows that the 
action of $G_{\hcL}$ on $\hcL$ is precisely the canonical
action of $G^{\sharp}_{L}=G_{\hcL}$ on the holonomy covering
space $\hcL$ of $L$.

\section{Equivariant Stability Theorems}  \label{secEquStaThe}

As an application of the stability theorems of Chapter
\ref{chapSta}, we shall give in this Section the equivariant
version of the Reeb-Thurston stability theorem.
We start with a theorem for the special case of a
simple equivariant foliation.

\begin{theo}  \label{IVtheo4}
Let $M$ be a $\eCe$-manifold of dimension $n$ equipped with a
$\eCe$-action of a discrete group $G$. Let $\cF$ be a simple
$G$-equivariant foliation on $M$ given by a submersion
$p:M\ra N$, where $N$ is a $\eCe$-manifold. Assume that
the orbit space of any leaf of $\cF$ is Hausdorff.
If $L$ is a leaf of $\cF$ with compact orbit space,
then for any $G$-invariant open neighbourhood of $L$
there exists a smaller $G$-invariant open neighbourhood of $L$ 
which is a union of leaves of $\cF$ with compact orbit spaces.
\end{theo}
\Proof
We can assume without loss of generality that $N=\RRR^{q}$.
Let $L$ be a leaf with compact orbit space and $V$ an open
$G$-invariant neighbourhood of $L$. Denote by $r:M\ra M/G$ the
quotient projection and by $p/G:M/G\ra\RRR^{q}$ the map induced by
$p$ as in Example \ref{IVex2}.

We know from Example \ref{IVex2} that the fibers of $p/G$ are
locally path-connected, and that the connected components of the
fibers of $p/G$ are exactly the orbit spaces of leaves of $\cF$.
In particular, the connected components of the fibers of $p/G$
are Hausdorff.

Now $r(V)$ is an open neighbourhood of the compact connected
component $L/G_{L}$ of $p/G$. Hence we can apply Corollary \ref{IIIcor4}
to find an open neighbourhood $W\subset r(V)$ of $L/G_{L}$ in $M/G$
which is a union of compact leaves of $p/G$. Now we take
$$ U=r^{-1}(W)\;.$$
Since $V$ is $G$-invariant, we have $U\subset V$. Clearly $U$
is an open $G$-invariant neighbourhood of $L$ which is a union
of leaves of $\cF$. If $L'$ is a leaf of $\cF$ with $L'\subset U$,
then
$$ r(L')=L'/G_{L'}\subset W\;,$$
therefore the orbit space of $L'$ is compact.
\eop 
\Rem
The condition that the action is leafwise separated is a
necessary one.
Indeed, take for example $G=1$ and define $M$ to be the
(non-Hausdorff) $\eCe$-manifold
$$ M=(S^{1}\times\RRR)\amalg (\RRR\times (0,\infty))/\sim\;,$$
where $((\phi,x),1)\sim ((e^{i\phi},x),0)$ if $0<x<e^{\phi}$.
Here  $S^{1}=\{\,z\in\CCC\,|\,|z|=1\,\}$. The submersion
$p:M\ra\RRR$ given by the second projection on both summands
in the coproduct has connected fibres. The fiber over $0$ is
diffeomorphic to $S^{1}$, but a fiber over any $x>0$ in not compact.
In other words, the classical Reeb stability theorem does not
hold for non-Hausdorff manifolds.

If we assume that $M$ is Hausdorff, this does not really simplify
the situation. For example, take $M=\RRR^{2}\setminus\{0\}$
and let $p$ be the second
projection. Define an action of $\ZZZ$ on $M$ along the fibers of $p$
generated by the $\eCe$-diffeomorphism
$$ \hat{1}(x,y)=(x+\frac{\sin(xy)}{2y},y)\;.$$
The orbit spaces of the leaves in $p^{-1}(0)$ are diffeomorphic with
$S^{1}$, but all the other leaves have non-compact non-Hausdorff
orbit spaces.

\begin{theo}[Equivariant Reeb stability]  \label{IVtheo5}
Let $\cF$ be a foliation on a finite-dimensional
$\eCe$-manifold $M$, invariant under a leafwise
separated $\eCe$-action of a discrete group $G$ on $M$.
Let $L$ be a leaf of $\cF$ such that
\begin{enumerate}
\item [(i)]  the orbit space of $L$ is compact, and
\item [(ii)] the equivariant holonomy group of $L$ is finite.
\end{enumerate}
Then for any $G$-saturated open neighbourhood of $L$ there exists a
smaller $G$-saturated open neighbourhood of $L$ which is a union of
leaves of $\cF$ with compact orbit spaces.
\end{theo}
\Proof
Consider the $\eCe$-principal $\GmqCe$-$G(M)$-bibundle $(E,p,\xW)$
associated to the equivariant foliation $\cF$, as in Section
\ref{secEquHol}. Since $p$ is a submersion, the fibers of $E$ are
locally path-connected. Let $\cFsh$ be the foliation
on $E$ given by the submersion $p$.

Let $\hcLa$ be a leaf of $\cFsh$. We know that $\xW(\hcLa)$ is
the underlying set of a leaf $L'$ of $\cF$, and the
restriction
$$ \zeta'=\xW|_{\hcLa}:\hcLa\ra L' $$
is the holonomy covering projection. Since
$$ G^{\sharp}_{L'}=G_{\hcLa}\;,$$
and the action of $G_{\hcLa}$ on $\hcLa$ is precisely the canonical
action of $G^{\sharp}_{L'}=G_{\hcLa}$ on the holonomy covering space
$\hcLa$ of $L'$, we have
$$ \hcLa/G^{\sharp}_{L'}=\hcLa/G_{\hcLa}\;.$$
By the assumption that the action of $G$ on $M$ is leafwise separated
this implies that the action of $G$ on $E$ is leafwise separated with
respect to $\cFsh$. Moreover, since
$$ |G(M)(\hcLa\cdot G)|=|G_{\hcLa}(\hcLa)|=\hcLa/G_{\hcLa}\;,$$
the leaves of the bibundle $E$ are $G(M)$-Hausdorff.

Next, we saw in Section \ref{secEquHol} that the leaves of
the bibundle $E$ are just the minimal $G$-invariant (disjoint) unions
of leaves of $\cF$. In particular, if $V$ is an open $G$-invariant
neighbourhood of $L$, it is also a $G(M)$-invariant neighbourhood
of the corresponding leaf of $E$. Moreover, a leaf of $\cF$ has
compact orbit space if and only if the corresponding leaf of $E$ is
$G(M)$-compact. In particular, the leaf $L\cdot G$ of $E$ which
correspond to the leaf $L$ is $G(M)$-compact. The
holonomy group of $L\cdot G$ is finite since it is isomorphic
to the equivariant holonomy group of $L$.
Thus we can apply Theorem \ref{IIItheo6}.
\eop

Next theorem gives the equivariant version of the Thurston
generalization \cite{Thu} of the Reeb stability theorem for
foliations of codimension one:

\begin{theo}[Equivariant Reeb-Thurston stability]  \label{IVtheo6}
Let $\cF$ be a foliation of codimension one on a
finite-dimensional $\eCe$-manifold $M$, invariant under a leafwise
separated $\eCe$-action of a discrete group $G$ on $M$.
Let $L$ be a leaf of $\cF$ with finitely generated equivariant
fundamental group such that
\begin{enumerate}
\item [(i)]   the foliation $\cF$ is equivariantly
              transversely orientable around $L$,
\item [(ii)]  the orbit space of $L$ is compact, and
\item [(iii)] the first equivariant cohomology group of $L$
              with coefficients in $\RRR$ is trivial.
\end{enumerate}
Then the equivariant holonomy group of $L$ is trivial, and 
for any $G$-saturated open neighbourhood of $L$ there exists a
smaller $G$-saturated open neighbourhood of $L$ which is a union
of leaves of $\cF$ with compact orbit spaces.
\end{theo}
\Rem
By definition, the first equivariant cohomology group 
$H_{G_{L}}^{1}(L;\RRR)$ is the first ordinary cohomology group of
the space $L\times_{G_{L}}\mbf{E} G_{L}$, where $\mbf{E} G_{L}$
is the universal $G_{L}$-bundle \cite{AtiBot}.
Since
$$ \pi^{G}_{1}(L)\cong\pi_{1}(L\times_{G_{L}}\mbf{E} G_{L})\;,$$
the Hurewicz formula implies that the condition (iii) is equivalent
to the condition
$$ \Hom(\pi_{1}^{G}(L),\RRR)=\{0\}\;.$$
In particular, the exact sequence (\ref{IVshoexaseq}) yields that
if the first ordinary cohomology group 
$H^{1}(L;\RRR)$ is trivial and $\Hom(G_{L},\RRR)=\{0\}$ then the
condition (iii) is satisfied. Next note that if $\pi_{1}(L)$ and
$G_{L}$ are finitely generated, then the equivariant fundamental
group $\pi_{1}^{G}(L)$ is also finitely generated.
\vspace{4 mm}

\Proof
As in the proof of Theorem \ref{IVtheo5}, we consider the
$\eCe$-principal $\GmqCe$-$G(M)$-bibundle $(E,p,\xW)$
associated to the equivariant foliation $\cF$ (see Section
\ref{secEquHol}). With the description of leaves and holonomy
of this bundle as in Section \ref{secEquHol}, the theorem follows
from Corollary \ref{IIIcor28}.
\eop

\begin{ex}  \label{IVex11}  \rm
(1)
Let $M=\RRR^{2}\setminus\{0\}$ and let $\cF$ be given by
the submersion $p=\pr_{2}:M\ra\RRR$. Take $G=\ZZZ$ and define
a properly discontinuous $\eCe$-action of $G$ on $M$ by
$$ (x,y)\cdot g=(2^{g}x,2^{g}y) $$
for any $(x,y)\in M$ and $g\in G$ (see Example \ref{Iex3} (5)).
The foliation $\cF$ is invariant under this action.

All the leaves of $\cF$ have trivial fundamental group and trivial
holonomy. Also, any leaf $L$ of $\cF$ with $0\not\in p(L)$ has trivial
isotropy group, therefore trivial equivariant fundamental group and
trivial equivariant holonomy group. Now let $L$ be the leaf 
$$ \{\,(x,0)\,|\,x>0\,\} $$
of $\cF$. The isotropy group of $L$ is clearly all $G$, and also the
equivariant fundamental group and the equivariant holonomy group of
$L$ are isomorphic to $\ZZZ$.

This example shows that we can not expect a stability theorem without
assuming that the equivariant holonomy group of a leaf with compact
orbit space is finite.

(2)
Let $M=\RRR\times S^{1}$, where $S^{1}=\{\,z\in\CCC\,|\,|z|=1\,\}$,
and let $\cF$ be the foliation on $M$ given by the submersion
$p=\pr_{2}:M\ra S^{1}$. Let $G=\ZZZ\oplus\ZZZ$, let 
$\alpha\in\RRR\setminus\QQQ$, and define a
$\eCe$-action of $G$ on $M$ by
$$ (x,z)\cdot (k,k')=(x+k,{\rm e}^{2\pi i \alpha k'}\!z) $$
for any $(x,z)\in M$ and $(k,k')\in G$. The foliation $\cF$ is
invariant under this action. If $L$ is a leaf of $\cF$, we have
$$ G_{L}\cong\pi_{1}^{G}(L)\cong\ZZZ\;,$$
but the equivariant holonomy group
of $L$ is trivial. The group $G_{L}$ acts properly discontinuously
on $L$ and the orbit space of $L$ is diffeomorphic to $S^{1}$,
so the action is leafwise separated. Thus the conditions of Theorem
\ref{IVtheo5} are satisfied.

(3)
Let $n\geq 2$, and denote by $\cM_{n}$ the $\eCe$-manifold
of $n\times n$-matrices over $\RRR$.
Consider the determinant map $\det:\cM_{n}\ra\RRR$.
This map is a submersion on the open submanifold of $\cM_{n}$
$$ M=\{\,A\in\cM_{n}\,|\,\dim\Ker A\leq 1\,\}\;,$$
with connected fibers. Hence the fibers of the determinant map
are the leaves of a foliation $\cF$ on $M$.
We will denote by $L_{0}$ the
fiber $\det^{-1}(0)$. We have $M=GL_{n}(\RRR)\cup L_{0}$,
where $GL_{n}(\RRR)$ is open, but not closed in $M$.

Now take $G$ to be the group $SL_{n}(\RRR)$ of matrices in
$\cM_{n}$ with determinant $1$.
There is a $\eCe$-action of the group $G$ on $M$ given by
$$ A\cdot g=g^{-1}A \;\;\;\;\;\;\;\; A\in M,\,g\in G\;.$$
This action is along the fibers of the determinant map. In particular,
$$ M/G\cong GL_{n}(\RRR)/G\cup L_{0}/G\cong\RRR^{\ast}\cup L_{0}/G\;,$$
where $\RRR^{\ast}=\RRR\setminus\{0\}$. For any $A\in L_{0}$, the kernel
of $A$ is a line in $\RRR^{n}$, hence we have a $\eCe$-map
$\Ker:L_{0}\ra\RRR P^{n-1}$. This map is clearly a surjective submersion.
On the other hand, we have
$$ \Ker (A\cdot g)=\Ker A $$
for any $A\in L_{0}$ and $g\in G$, therefore $\Ker$ factors as
$$ f:L_{0}/G\lra\RRR P^{n-1}\;.$$
Since for any $A,B\in L_{0}$ with $\Ker A=\Ker B$ one can
find a matrix $g\in G$ such that $g A=B$, the map
$f$ is a homeomorphism. In particular, $L_{0}/G$ is compact.

The orbit spaces of the leaves of $\cF$ which lie in $GL_{n}(\RRR)$
are the points in $\RRR^{\ast}$, and the orbit space of $L_{0}$ is
homeomorphic to $\RRR P^{n-1}$. The space $L_{0}/G$ is not
open in $M/G$. This example satisfies the conditions of
Theorem \ref{IVtheo4}.

(4)
Let $N$ be a connected $\eCe$-manifold of dimension $n$.
A family
$(f_{\lambda})_{\lambda\in\Lambda}$ of
$\eCe$-diffeomorphisms of $N$ has a {\em compact filling}
if there exists a compact set $K\subset N$ such that
$$ \bigcup\{(f_{\lambda_{1}}^{n_{1}}\com
   \ldots\com f_{\lambda_{k}}^{n_{k}})(K)\;|
   \,k\geq 1,\;n_{i}\in\ZZZ,\;\lambda_{i}\in\Lambda\}=N\;.$$
Observe that $(f_{\lambda})_{\lambda\in\Lambda}$
has a compact filling if and only if the quotient space
$N/G_{(f_{\lambda})}$ is compact, where $G_{(f_{\lambda})}$
is the group of diffeomorphisms of $N$ generated by
$(f_{\lambda})_{\lambda\in\Lambda}$.
The family $(f_{\lambda})_{\lambda\in\Lambda}$
is called {\em separated} if the space
$N/G_{(f_{\lambda})}$ is Hausdorff.

A $\eCe$-perturbation of a family
$(f_{\lambda})_{\lambda\in\Lambda}$ of
$\eCe$-diffeomorphisms of $N$ is a family of $\eCe$-maps
$$ (P_{\lambda}:(-\varepsilon,\varepsilon)
   \times N\lra N)_{\lambda\in\Lambda} $$
such that
$P^{t}_{\lambda}=P_{\lambda}(t,\oo)$
is a $\eCe$-diffeomorphism and
$P^{0}_{\lambda}=f_{\lambda}$, for any
$\lambda\in\Lambda$ and $t\in (-\varepsilon,\varepsilon)$.
The perturbation $(P_{\lambda})_{\lambda\in\Lambda}$ is
called {\em separated} if the family
$(P^{t}_{\lambda})_{\lambda\in\Lambda}$ is separated, for
any $t\in (-\varepsilon,\varepsilon)$.

Assume that $(P_{\lambda})_{\lambda\in\Lambda}$ is a separated
$\eCe$-perturbation of a family
$(f_{\lambda})_{\lambda\in\Lambda}$
of $\eCe$-diffeo\-morphisms of $N$ with a compact filling.
We will show that Theorem \ref{IVtheo4} implies that there exists
$0<\delta<\varepsilon$ such that for any $|\,t\,|\leq\delta$, the
family
$$ (P^{t}_{\lambda})_{\lambda\in\Lambda} $$
has a compact filling.

To see this, let $M=(-\varepsilon,\varepsilon)\times N$ and let $\cF$
be the foliation of codimension one on $M$ given by the first projection
$pr_{1}$. For any $\lambda$, we have a $\eCe$-diffeomorphism
$\alpha_{\lambda}$ of $M$, given by
$\alpha_{\lambda}(t,x)=(t,P_{\lambda}(t,x))$,
which preserves $\cF$. Let $G$ be the group of diffeomorphisms
of $M$ generated by $(\alpha_{\lambda})_{\lambda\in\Lambda}$.
This action of $G$ on $M$ is leafwise separated since the
perturbation is separated.
Observe that a leaf $L_{t}=pr_{1}^{-1}(t)$ has compact orbit space
for the action of $G$ precisely if the family 
$(P^{t}_{\lambda})_{\lambda\in\Lambda}$
has a compact filling.
\end{ex}

\chapter{Invariants of Hilsum-Skandalis Maps}
\label{chapInvHilSkaMap}
\startSchapterskip

In this Chapter we study some algebraic invariants
of topological gro\-up\-oids. First we use the Hilsum-Skandalis
maps to introduce (higher) homotopy groups and singular homology
groups of topological groupoids.
These invariants can be described particularly
easy if the groupoids are \'{e}tale. However, our main
example is the groupoid of leaves of a
foliation, which is not \'{e}tale but Morita equivalent to an
\'{e}tale groupoid. It is therefore our primary objective
to establish that the higher homotopy groups and the singular
homology groups are invariant under Morita equivalence.

In Section \ref{secBib} we proved that the Morita equivalent
groupoids are exactly the isomorphic groupoids in the category
of Hilsum-Skandalis maps. Thus the invariance under Morita
equivalence is guaranteed if we define the higher homotopy
groups and singular homology groups as functors on this category.
Moreover, the definition of these algebraic invariants is
formally the same as the classical definitions for topological
spaces, provided that we replace the category of topological
spaces with that of Hilsum-Skandalis maps.

If $\GG$ is an object-separated finite-dimensional
\'{e}tale $\eCr$-groupoid, the complex $\eCrc$-functions on
$\GG_{1}$ with the convolution product
form the associative algebra $\eCrc(\GG)$ with local units, see
\cite{BryNis,Con78,Con82,FacSka,HilS,Ren,RenMW}.
Moreover, there is a natural
norm on $\eCrc(\GG)$ so that the completion of $\eCrc(\GG)$ is
a $C^{\ast}$-algebra. In the last Section we associate to
a Hilsum-Skandalis $\eCr$-map $E$ between object-separated
finite-dimensional \'{e}tale $\eCr$-groupoids $\HH$
and $\GG$ an isomorphism class of a $\eCrc(\GG)$-$\eCrc(\HH)$-bimodule 
$\eCrc(E)$. We show that this gives a functor from the category of 
Hilsum-Skandalis $\eCr$-maps between separated finite-dimensional
\'{e}tale $\eCr$-groupoids to the category of isomorphism classes
of locally unital bimodules over algebras with local units.
This implies, for example, that the algebras associated
to $\eCr$-Morita equivalent separated finite-dimensional 
\'{e}tale $\eCr$-groupoids have isomorphic
cyclic homology groups \cite{Con85,Con86,Qui88}.

\section{Homotopy Groups of Topological Groupoids}
\label{secHomGro}

In this Section we present a possible definition
of higher homotopy groups
of topological groupoids. However, before we give the definition,
we have to study how the principal bibundles can be amalgamated.

\subsection{Amalgamations of Principal Bibundles}
\label{subsecAmaPriBib}

Let $\GG$ and $\HH$ be topological groupoids, and let $(E,p,\xW)$
be a principal $\GG$-$\HH$-bibundle. If $U$ is an $\HH$-invariant
subset of $\HH_{0}$, then we have the full open subgroupoid
$\HH|_{U}$ of $\HH$ with $(\HH|_{U})_{0}=U$, which is naturally
isomorphic to $\HH(U)$. Denote by $\iii:\HH|_{U}\ra\HH$ the
embedding. Define
$$ E|_{U}=\xW^{-1}(U)\;.$$
Since all the structure maps of $E$ restrict to $E|_{U}$, it is
clear that $E|_{U}$ is a principal $\GG$-$\HH|_{U}$-bibundle,
naturally isomorphic to $E\otimes\angs{\iii}$. The bibundle
$E|_{U}$ is called the {\em restriction} of $E$ to $\HH|_{U}$.
With this, we get a notion of restriction of a Hilsum-Skandalis
map. Observe that if $p:X\ra Y$ is a continuous
map and $U\subset X$, then
$\angs{p|_{U}}=\angs{p}|_{U}$.

Assume now that $(U_{i})_{i\in I}$ is a cover of $\HH_{0}$,
consisting of $\HH$-invariant subsets of $\HH_{0}$,
which is either
\begin{enumerate}
\item [(i)]  open, or
\item [(ii)] closed and locally finite.
\end{enumerate}
Further assume that we have a family $(E_{i})_{i\in I}$, where each
$E_{i}=(E_{i},p_{i},\xW_{i})$ is a principal
$\GG$-$\HH|_{U_{i}}$-bibundle.  
We would like to amalgamate these bibundles into a single
principal $\GG$-$\HH$-bibundle, just like this can be done for
continuous maps on topological spaces. Of course, we should assume
that the bibundles match on the intersections $U_{i}\cap U_{j}$.
In our context, this means that we should give a family of
equivariant homeomorphisms
$$ (\alpha_{ij}:E_{j}|_{U_{i}\cap U_{j}}\lra E_{i}|_{U_{i}\cap
   U_{j}})_{i,j\in I} $$
(more precisely, the map $\alpha_{ij}$ is
$\GG$-$\HH|_{U_{i}\cap U_{j}}$-equivariant) which satisfy
the cocycle condition
$$ \alpha_{ij}(\alpha_{jk}(e))=\alpha_{ik}(e) $$
for any $e\in E_{k}|_{U_{i}\cap U_{j}\cap U_{k}}$. In particular,
$\alpha_{ii}=\id_{E_{i}}$.

With this, we define a principal $\GG$-$\HH$-bibundle
$(E,p,\xW)$ as follows: The space $E$ is the quotient of 
$$ \coprod_{i\in I}E_{i}=
   \{\,(e,i)\,|\,i\in I,\,e\in E_{i}\,\} $$
by identifying $(e,j)$ with $(\alpha_{ij}(e),i)$, for any $i,j\in I$
and $e\in E_{j}|_{U_{i}\cap U_{j}}$. Denote by 
$[e,i]\in E$ the equivalence class of
an element $(e,i)\in\coprod_{i\in I}E_{i}$.
The quotient projection $\coprod_{i\in I}E_{i}\ra E$
is injective on the closed-open subsets $E_{i}\times\{i\}$
of $\coprod_{i\in I}E_{i}$. In fact, observe that the canonical map
$$ t_{i}:E_{i}\ra E\;,$$
given by $t_{i}(e)=[e,i]$, is an open embedding in case (i) and
closed embedding in case (ii). In particular,
$(t_{i}(E_{i}))_{i\in I}$
is either an open or a closed locally finite
cover of $E$. 

Now we define the structure maps of $E$ by
$p([e,i])=p_{i}(e)$ and $\xW([e,i])=\xW_{i}(e)$. Obviously this is
well-defined, and since $p\com t_{i}=p_{i}$ and $\xW\com t_{i}=\xW_{i}$,
the maps $p$ and $\xW$ are continuous. Next, observe that $\xW$ is an
open surjection, since $\xW|_{t_{i}(E_{i})}:t_{i}(E_{i})\ra U_{i}$ is
an open surjection for any $i\in I$. Similarly we define the actions
of $\GG$ and $\HH$ by amalgamating the actions on
$t_{i}(E_{i})\cong E_{i}$, i.e.
$g\cdot [e,i]=[g\cdot e,i]$ and $[e,i]\cdot h=[e\cdot h,i]$. Again,
these two actions are continuous since they are just the actions
on $E_{i}$ when restricted to $t_{i}(E_{i})\cong E_{i}$.
It is easy to check that $E$ is indeed a
principal $\GG$-$\HH$-bibundle.

The principal $\GG$-$\HH$-bibundle $E$ is
called the {\em amalgamation} of $(E_{i})$ with respect to
$(\alpha_{ij})$. We will denote this bibundle by
$$ \Sigma(E_{i},\alpha_{ij})\;.$$
Observe that for each $i\in I$ we have the natural
$\GG$-$\HH|_{U_{i}}$-equivariant homeomorphism
$$ t_{i}:E_{i}\lra\Sigma(E_{i},\alpha_{ij})|_{U_{i}}\;.$$
Moreover, for any $e\in E_{j}|_{U_{i}\cap U_{j}}$ we have
$$ t_{j}(e)=t_{i}(\alpha_{ij}(e))\;.$$
\vspace{4 mm}

\begin{ex}  \label{Vex1}  \rm
(1) Let $\GG$ be an \'{e}tale groupoid, $X$ a topological space with
an open cover $\cU=(U_{i})_{i\in I}$
and $c=(c_{ij})$ a $\GG$-cocycle on $\cU$.
For any $i\in I$ we have the map
$$ s_{i}=\dom\com c_{ii}:U_{i}\ra\GG_{0}\;,$$
which gives a principal $\GG$-bundle
$$ \angs{s_{i}}=\GG_{1}\times_{\GG_{0}}U_{i} $$
over $U_{i}$. Moreover, for any $i,j\in I$ we have the
equivariant map
$\alpha_{ij}:\angs{s_{j}}|_{U_{i}\cap U_{j}}\ra
\angs{s_{i}}|_{U_{i}\cap U_{j}}$ given by
$$ \alpha_{ij}(g,x)=(g\com c_{ji}(x),x)\;.$$
Since $c$ is a cocycle, it follows that 
$\alpha_{ij}(\alpha_{jk}(g,x))=\alpha_{ik}(g,x)$, thus we can
construct the principal $\GG$-bundle $\Sigma(\angs{s_{i}},\alpha_{ij})$.
But observe that there is a natural equivariant homeomorphism
$$ \Sigma(\angs{s_{i}},\alpha_{ij})\cong\Sigma(c)\;,$$
where $\Sigma(c)$ is the principal $\GG$-bundle associated to
the cocycle $c$. Thus the construction of $\Sigma(c)$
in Section \ref{secHaeStr} is in fact a special case of amalgamation
of principal bibundles.

(2) Let $\GG$ be an \'{e}tale groupoid and $X$ a topological space.
If $\cW=(W_{i})_{i\in I}$ is a closed locally finite cover of $X$,
we can define a $\GG$-cocycle on $\cW$ in exactly the same way as
a $\GG$-cocycle on an open cover, 
i.e. as a family of continuous maps
$$ c=(c_{ij}:W_{i}\cap W_{j}\lra\GG_{1}) $$
which satisfy the usual conditions, see Section \ref{secHaeStr}.
We say that $c$ is constant (with the value $a\in\GG_{0}$)
if $c_{ij}(x)=1_{a}$ for any $x\in W_{i}\cap W_{j}$.
We can also define the families which intertwine $\GG$-cocycles on
$\cW$,
and hence obtain the set $H^{1}(\cW,\GG)$ of cohomology classes of
$\GG$-cocycles on $\cW$. Finally, we can take the filtered colimit
$$ H^{1}_{cl}(X,\GG)=\lim_{\ra_{\,\cW}}H^{1}(\cW,\GG) $$
with respect to the partially ordered set of closed locally finite
covers of $X$. We call the elements of $H^{1}_{cl}(X,\GG)$ the
{\em closed Haefliger $\GG$-structures} on $X$.

If $c$ is a $\GG$-cocycle on a closed locally finite cover of $X$,
we take $s_{i}$ and $\alpha_{ij}$ to be defined
in the same way as in (1), and we obtain a
principal $\GG$-bundle
$$ \Sigma(c)=\Sigma(s_{i},\alpha_{ij})\;.$$
It is easy to verify that if $c'$ is another $\GG$-cocycle on a closed
locally finite cover of $X$ which represents the same closed Haefliger
$\GG$-structure on $X$ as $c$, then $\Sigma(c)$ and $\Sigma(c')$ are 
isomorphic.

In fact, if $X$ is paracompact, $\Sigma$ induces a bijective
correspondence between the closed Haefliger $\GG$-structures on $X$ and
the isomorphism classes of principal $\GG$-bundles over $X$.
In particular we have $H^{1}(X,\GG)\cong H^{1}_{cl}(X,\GG)$.
The proof is analogous to the proof of Proposition \ref{Iprop7}.
If $X$ is compact, a closed Haefliger $\GG$-structure on $X$ 
can be represented by a $\GG$-cocycle defined on a closed finite
cover of $X$.
\end{ex}

\subsection{Definition of Homotopy Groups}
\label{subsecDefHomGro}

A {\em marked topological groupoid} is a pair $(\GG,A)$, where
$\GG$ is a topological groupoid and $A$ is a subset of $\GG_{0}$.
If $(\HH,B)$ is another marked topological groupoid, a {\em marked
principal $(\GG,A)$-$(\HH,B)$-bibundle} is a pair
$(E,s)$, where $E=(E,p,\xW)$ is a principal $\GG$-$\HH$-bibundle and
$s:B\ra E$ is a continuous section of $\xW$ such that
$$ p(s(B))\subset A\;.$$
The section $s$ is called the {\em trivialization} of $E$
over $B$.

If $(E,s)$ and $(E',s')$ are two marked principal
$(\GG,A)$-$(\HH,B)$-bibundles, a continuous map
$\alpha:E\ra E'$ is called $(\GG,A)$-$(\HH,B)$-equivariant if
it is $\GG$-$\HH$-equivariant (and hence a homeomorphism) and
$$ \alpha\com s=s'\;.$$
If such a map exists, the marked bibundles $(E,s)$ and $(E',s')$
are called isomorphic. An isomorphism class of marked principal
$(\GG,A)$-$(\HH,B)$-bibundles is called a {\em marked
Hilsum-Skandalis map} from $(\HH,B)$ to $(\GG,A)$.
As for the principal bibundles, we will denote an isomorphism class
of a marked principal bibundle $(E,s)$ again by $(E,s)$.

Let $(\KK,C)$ be another marked topological groupoid,
let $((E,p,\xW),s)$ be a marked principal
$(\GG,A)$-$(\HH,B)$-bibundle, and let $((E',p',\xW'),s')$ be a marked
principal $(\HH,B)$-$(\KK,C)$-bibundle. Then we define
$$ s\otimes s':C\lra E\otimes E' $$
by $(s\otimes s')(c)=s(p'(s'(c)))\otimes s'(c)$, for any $c\in C$.
It is clear that $(E\otimes E',s\otimes s')$ is a marked
principal $(\GG,A)$-$(\KK,C)$-bibundle, which we will call
the {\em tensor product} of marked principal bibundles
$(E,s)$ and $(E',s')$, and denote by
$$ (E,s)\otimes (E',s')\;.$$
With this tensor product, the marked Hilsum-Skandalis maps
between marked topological groupoids form a category, which
will be denoted by $\MGpd$.
For example, the identity in $\MGpd((\GG,A),(\GG,A))$ is 
the isomorphism class of the marked principal
$(\GG,A)$-$(\GG,A)$-bibundle
$$ ((\GG_{1},\cod,\dom),\uni|_{A})\;.$$
The full subcategory of $\MGpd$ with objects all marked \'{e}tale
groupoids will be denoted by $\MGpde$.

A {\em pointed topological groupoid} is a marked
topological groupoid $(\GG,A)$ such that $A$ is a one-point set,
i.e. $A=\{a\}$. We will write in this case $(\GG,a)=(\GG,\{a\})$.
If $(\HH,b)$ is another pointed topological groupoid and
$((E,p,\xW),s)$ a marked principal $(\GG,a)$-$(\HH,b)$-bibundle, the section
$s$ is determined by the point $e=s(b)$ such that $p(e)=a$.
We shall thus denote $(E,s)=(E,e)$ and call $(E,e)$ a {\em pointed
principal $(\GG,a)$-$(\HH,b)$-bibundle}. An isomorphism class of
pointed principal bibundles is called a {\em pointed Hilsum-Skandalis map}.
The full subcategory of $\MGpd$ with objects all pointed topological
groupoids will be denoted by
$$ \PGpd\;,$$
and the full subcategory of
$\MGpde$ with objects all pointed \'{e}tale topological groupoids will
be denoted by $\PGpde$. Observe that the category of pointed topological
spaces $\PTop$ is a full subcategory of $\PGpde$.

\begin{ex}  \label{Vex2}  \rm
(1) A marked topological space is a pair $(X,B)$ with $B\subset X$, i.e.
a topological pair. Also, a marked Hilsum-Skandalis map
between two marked topological spaces is clearly just a continuous map
between topological pairs. Hence the category of topological
pairs is a full subcategory of $\MGpd_{e}$. If $(X,B)$ is a
marked topological space and $(\GG,A)$ a marked topological
groupoid, a marked principal $(\GG,A)$-$(X,B)$-bibundle will be
referred to as marked principal $(\GG,A)$-bundle over $(X,B)$.

(2) Let $\GG$ be a topological groupoid. Then
$\aangs{\GG}=(\GG,\GG_{0})$ is
a marked topological groupoid. Moreover, if
$\phi:\HH\ra\GG$ is a continuous functor between topological groupoids,
the principal $\GG$-$\HH$-bibundle $\angs{\phi}$ is canonically
marked by $(\uni\com\phi_{0},\id_{\HH_{0}}):\HH_{0}\ra\angs{\phi}$,
so
$$ \aangs{\phi}=(\angs{\phi},(\uni\com\phi_{0},\id_{\HH_{0}})) $$
is a marked principal $(\GG,\GG_{0})$-$(\HH,\HH_{0})$-bibundle.
This gives a functor
$$ \aangs{\oo}:\Gpd\lra\MGpd\;.$$
This functor is full and faithful. Indeed,
Proposition \ref{IIprop5} implies that $\aangs{\oo}$ is
full. On the other hand, if $\phi,\psi:\HH\ra\GG$ are continuous
functors such that there is a $(\GG,\GG_{0})$-$(\HH,\HH_{0})$-equivariant 
homeomorphism
$$ \alpha:\aangs{\phi}\lra\aangs{\psi}\;,$$
this clearly implies $\phi_{0}=\psi_{0}$, and for any $h\in\HH_{1}$
we have
\begin{eqnarray*}
(\psi(h),\dom h) &=& \alpha(1_{\phi_{0}(\cod h)},\cod h)\cdot h=
                     \alpha(\phi(h),\dom h)                       \\
&=& \phi(h)\cdot\alpha(1_{\phi_{0}(\dom h)},\dom h)=
    (\phi(h),\dom h)\;,
\end{eqnarray*}
hence $\phi=\psi$. Thus $\aangs{\oo}$ is faithful.

Observe that if $(E,s)$ and $(E',s')$ are isomorphic marked principal
$(\GG,\GG_{0})$-$(\HH,\HH_{0})$-bibundles, then there exists a unique
isomorphism between them.

(3) Observe that $\cGpd$ is also a full subcategory of $\MGpd$ if
we identify a topological groupoid $\GG$ with the marked
topological groupoid $(\GG,\emptyset)$, and a principal
bibundle $E$ with the marked principal bibundle $(E,\emptyset)$.

(4) Let $(E,s)$ be a marked principal $(\GG,A)$-$(\HH,B)$-bibundle,
and let $U$ be an $\HH$-invariant subset of $\HH_{0}$.
Define $(\HH,B)|_{U}=(\HH|_{U},B\cap U)$, and let
$$ (E,s)|_{U}=(E|_{U},s|_{B\cap U})\;.$$
This is clearly a marked principal $(\GG,A)$-$(\HH,B)|_{U}$-bibundle,
which we call the
{\em restriction} of $(E,s)$ on $(\HH,B)|_{U}$. This gives
a notion of restriction of marked Hilsum-Skandalis maps.
\end{ex}

Let $(\GG,a)$ be a pointed topological groupoid, and let $n\geq 0$.
Denote by $\III^{n}=[0,1]^{n}$ the $n$-cube and by $\rob\III^{n}$ its
boundary. In particular, $\III^{0}$ is a one-point space and
$\rob\III^{0}=\emptyset$. Now put
$$ \Theta_{n}(\GG,a)=\MGpd((\III^{n},\rob\III^{n}),(\GG,a))\;.$$
The elements of $\Theta_{n}(\GG,a)$ will be called the
{\em Hilsum-Skandalis $n$-loops} in $(\GG,a)$.

For any $t\in\III$,
denote by $l_{t}:\III^{n}\ra\III^{n}\times\III$ the continuous
map given by $l_{t}(x)=(x,t)$. Note that
$l_{t}\in\MGpd((\III^{n},\rob\III^{n}),
(\III^{n}\times\III,\rob\III^{n}\times\III))$. A marked
Hilsum-Skandalis map
$$ \gH\in\MGpd((\III^{n}\times\III,\rob\III^{n}\times\III),
   (\GG,a)) $$
is called the {\em homotopy} between the Hilsum-Skandalis $n$-loops
$\gH\otimes l_{0}$ and $\gH\otimes l_{1}$ in $(\GG,a)$.
We will denote $\gH_{t}=\gH\otimes l_{t}$, for any $t\in\III$.
Two Hilsum-Skandalis $n$-loops in $(\GG,a)$ are {\em homotopic}
if there exists a homotopy between them.

\begin{lem}  \label{Vlem3}
Let $(\GG,a)$ be a pointed topological groupoid.
Then the homotopy is an equivalence relation on
$\Theta_{n}(\GG,a)$.
\end{lem}
\Proof
Assume that $\gH$ and $\gH'$
are homotopies such that $\gH_{1}=\gH'_{0}$. We shall define
a homotopy between $\gH_{0}$ and $\gH'_{1}$ as follows:
Represent $\gH$ respectively $\gH'$ by marked principal bibundles
$(E,s)$ respectively $(E',s')$. Define
$v:\III^{n}\times[0,1/2]\lra\III^{n}\times\III$
by
$$ v(x,t)=(x,2t)\;,$$
and $v':\III^{n}\times[1/2,1]\lra\III^{n}\times\III$ by
$$ v'(x,t)=(x,2t-1)\;.$$
Now observe that the composition $(E,s)\otimes v$
is a marked principal
$(\GG,a)$-bundle over
$(\III^{n}\times[0,1/2],\rob\III^{n}\times[0,1/2])$,
and $(E',s')\otimes v'$ is a marked principal $\GG$-bundle over
$(\III^{n}\times[1/2,1],\rob\III^{n}\times[1/2,1])$.
Moreover, since $\gH_{1}=\gH'_{0}$ there exists an
equivariant homeomorphism of marked principal $(\GG,a)$-bundles
over $(\III^{n}\times\{1/2\},\rob\III^{n}\times\{1/2\})$
$$ \alpha:(E,s)\otimes v|_{\III^{n}\times\{1/2\}}\lra
   (E',s')\otimes v'|_{\III^{n}\times\{1/2\}}\;.$$
The bundles $E\otimes v$ and $E'\otimes v'$ thus amalgamate
with respect to $\alpha$
into a single principal $\GG$-bundle $E''$ over $\III^{n}\times\III$,
as we described in Subsection \ref{subsecAmaPriBib}. Moreover,
since $\alpha$ is a homeomorphism of marked bundles, it is easy
to see that using $s$ and $s'$ one can define
$s'':\rob\III^{n}\times\III\ra E''$, such that
$(E'',s'')$ is a marked principal $(\GG,a)$-bibundle over
$(\III^{n}\times\III,\rob\III^{n}\times\III)$ which represents
a homotopy between
$\gH_{0}$ and $\gH'_{1}$. This proves the transitivity, and
the rest is trivial.
\eop

Assume now that $n\geq 1$, and
let $(E,s)$ and $(E',s')$ be marked principal $(\GG,a)$-bundles
over $(\III^{n},\rob\III^{n})$. We shall now define a marked principal $(\GG,a)$-bundle
$$ (E\ast E',s\ast s') $$
over $(\III^{n},\rob\III^{n})$,
called the {\em concatenation} of $(E,s)$ and $(E',s')$, as follows:
Define $u:[1/2,1]\times\III^{n-1}\ra\III^{n}$ by $u(t,x)=(2t-1,x)$
and $u':[0,1/2]\times\III^{n-1}\ra\III^{n}$ by $u'(t,x)=(2t,x)$.
Observe that $(E,s)\otimes u$ is a marked principal $(\GG,a)$-bundle
over $([1/2,1]\times\III^{n-1},\rob([1/2,1]\times\III^{n-1}))$, and
$(E',s')\otimes u'$ is a marked principal $(\GG,a)$-bundle
over $([0,1/2]\times\III^{n-1},\rob([0,1/2]\times\III^{n-1}))$.
Moreover, there is exactly one equivariant homeomorphism of
marked principal bibundles
$$ \alpha:(E,s)\otimes u|_{\{1/2\}\times\III^{n-1}}\lra
   (E',s')\otimes u'|_{\{1/2\}\times\III^{n-1}}\;,$$
since $\{1/2\}\times\III^{n-1}=\rob([0,1/2]\times\III^{n-1})\cap
\rob([1/2,1]\times\III^{n-1})$, see Example \ref{Vex2} (2).
Thus we can amalgamate the bundles $(E,s)\otimes u$ and
$(E',s')\otimes u'$ with respect to $\alpha$. We take $E\ast E'$
to be the amalgamation.
Using $s$ and $s'$ we can define a section
$s\ast s':\rob\III^{n}\ra E\ast E'$ such that
$(E\ast E',s\ast s')$ is a marked principal $(\GG,a)$-bibundle
over $(\III^{n},\rob\III^{n})$.

The concatenation of principal $(\GG,a)$-bundles over
$(\III^{n},\rob\III^{n})$ gives a notion of concatenation of
Hilsum-Skandalis $n$-loops in $(\GG,a)$. Indeed, if
$(E,s)\cong (E_{1},s_{1})$ and $(E',s')\cong (E'_{1},s'_{1})$,
then these two isomorphisms amalgamate into an isomorphism
$$ (E\ast E',s\ast s')\cong (E_{1}\ast E'_{1},s_{1}\ast
   s'_{1})\;.$$
Moreover, the concatenation on $\Theta_{n}(\GG,a)$ induces a
concatenation between the homotopy classes of Hilsum-Skandalis
$n$-loops in $(\GG,a)$, as one can check by amalgamating the
homotopies. Finally, it is easy to see that with the operation
of concatenation, the set of homotopy classes of Hilsum-Skandalis
$n$-loops in $(\GG,a)$ becomes a group, called the
{\em $n$-th homotopy group} of $(\GG,a)$. We will denote this
group by 
$$ \Pi_{n}(\GG,a)\;.$$
For example, the unit element in $\Pi_{n}(\GG,a)$ is the homotopy
class of the marked principal $(\GG,a)$-bundle
$(\angs{\varpi(a)},s)$, where $\varpi(a):\III^{n}\ra\GG_{0}\subset\GG$
is the constant map with the image point $a$, and
$s(x)=(1_{a},x)$ for any $x\in\rob\III^{n}$.
Observe that we can also define $\Pi_{0}(\GG,a)$ as the set of homotopy
classes of $0$-loops in $(\GG,a)$.

Now let $(E,e)\in\PGpd((\HH,b),(\GG,a))$. We can define
for any $n\geq 0$ a function
$$ (E,e)_{\#}=\MGpd((\III^{n},\rob\III^{n}),(E,e)):
   \Theta_{n}(\HH,b)\lra\Theta_{n}(\GG,a)\;,$$
i.e. $(E,e)_{\#}$ is given by the tensor product with $(E,e)$.
The function $(E,e)_{\#}$ preserves homotopy, because one can
compose a homotopy with $(E,e)$. Therefore there is
an induced map
$$ (E,e)_{\ast}=\Pi_{n}(E,e):\Pi_{n}(\HH,b)\lra\Pi_{n}(\GG,a)\;.$$
It is straightforward to check that this is a homomorphism of
groups if $n\geq 1$.
We get a functor $\Pi_{0}$ from $\PGpd$ to the category
of sets, and a functor $\Pi_{n}$ from $\PGpd$ to the
category of groups, for any $n\geq 1$. One can prove that
$\Pi_{n}(\GG,a)$ is in fact abelian if $n\geq 2$, using exactly the
same argument as in the proof of this fact for the higher homotopy
groups of topological spaces.

\begin{ex}  \label{Vex5}  \rm
(1) Let $(X,b)$ be a pointed topological space. In this case the
definition of homotopy groups of $(X,b)$ as a pointed
groupoid is exactly the classical definition of homotopy groups
for the pointed space $(X,b)$, hence
$$ \Pi_{n}(X,b)=\pi_{n}(X,b)\;\;\;\;\;\;\;\; n\geq 0\;.$$

(2) Let $G$ be a discrete group. As an \'{e}tale groupoid,
$G$ is canonically pointed. A principal $G$-bundle over $\III^{n}$ is
a covering space, hence isomorphic to the product
$G\times\III^{n}$. A marked principal $G$-bundle over
$(\III^{n},\rob\III^{n})$ is thus isomorphic to
$$ (G\times\III^{n},s)\;,$$
where $s:\rob\III^{n}\ra G\times\III^{n}$ is a section of
the second projection. In other words, $s$ may be represented
by a locally constant function from $\rob\III^{n}$ to $G$.
Now if $n\geq 2$ such a function has to be constant, therefore
$$ \Pi_{n}(G)=1\;\;\;\;\;\;\;\; n\geq 2\;.$$
Also it is clear that $\Pi_{0}(G)=1$. Finally, for $n=1$, any
marked principal $G$-bundle over $(\III,\rob\III)$ is
isomorphic to $(G\times\III,s)$ with
$s(0)=(1,0)$, hence uniquely characterized by $g\in G$ such that
$s(1)=(g,1)$. Moreover, this element characterize the homotopy
class of the bundle as well. It is then straightforward that
$$ \Pi_{1}(G)=G\;.$$
\end{ex}

\subsection{Elementary Properties of Homotopy Groups}
\label{subsecElePro}

\begin{prop}  \label{Vprop4}
Let $\GG$ and $\HH$ be Morita equivalent topological groupoids.
Then for any $b\in\HH_{0}$ there exists $a\in\GG_{0}$ such that
$$ \Pi_{n}(\HH,b)\cong\Pi_{n}(\GG,a)\;\;\;\;\;\;\;\; n\geq 0\;.$$
In particular, if $\phi:\HH\ra\GG$ is an essential equivalence,
then
$$ (\angs{\phi},(1_{\phi_{0}(b)},b))_{\ast}:\Pi_{n}(\HH,b)\lra
   \Pi_{n}(\GG,\phi_{0}(b)) $$
is an isomorphism.
\end{prop}
\Proof
Choose $b\in\HH_{0}$. By Corollary \ref{IIcor7}
there exists a principal $\GG$-$\HH$-bibundle $(E,p,\xW)$
which represents an isomorphism between $\HH$ and $\GG$ in
$\cGpd$. Let $(E',p',\xW')$ be a bundle which represents the inverse
of $E$. In particular, there exists an equivariant homeomorphism
$\alpha:\HH_{1}\ra E'\otimes E$. Let $\alpha(1_{b})=e'\otimes e$.
Put $a=p(e)=\xW'(e')$. Hence $(E,e)$ is a pointed principal
$(\GG,a)$-$(\HH,b)$-bibundle and $(E',e')$ is a pointed principal
$(\HH,b)$-$(\GG,a)$-bibundle.
We claim that $(E,e)$ and $(E',e')$
are inverse to each other in the category $\PGpd$.

To see this, observe first that $(E',e')\otimes (E,e)\cong
(\HH_{1},a)$. Next, we know that there exists an equivariant
homeomorphism $\beta:E\otimes E'\ra\GG_{1}$, hence
$(E,e)\otimes (E',e')\cong (\GG_{1},g)$ for some $g\in\GG_{1}$
with $\dom g=\cod g=a$. But since $(\GG_{1},g)$ is clearly invertible
in $\PGpd$ with the inverse given by $(\GG_{1},g^{-1})$, this yields
that $(\GG_{1},g)\cong (\GG_{1},1_{a})$ and
that $(E,e)$ is the inverse of $(E',e')$.
The first part of the proposition now follows from the functoriality
of $\Pi_{n}$ on the category $\PGpd$.
For the second part, we can take $E=\angs{\phi}$ and
$e=(1_{\phi_{0}(b)},b)$.
\eop

In Section \ref{secFunGro} we defined the fundamental group
of a topological groupoid $\HH$ in terms of $\HH$-loops in
$\HH_{0}$. We shall now prove that the fundamental
group $\pi_{1}(\HH,b)$ of $\HH$ with a base-point $b$ is isomorphic
to the first homotopy group $\Pi_{1}(\HH,b)$ of the pointed groupoid
$(\HH,b)$ in case that $\HH$ is a suitable \'{e}tale groupoid.

\begin{theo}  \label{Vtheo6}
Let $\HH$ be an \'{e}tale groupoid such that $\HH_{0}$ is
locally path-connected, and let $b\in\HH_{0}$. Then there is
an isomorphism
$$ \Phi:\pi_{1}(\HH,b)\lra\Pi_{1}(\HH,b)\;.$$
\end{theo}
\Rem
The isomorphism $\Phi$ is natural in the sense that if
$\GG$ is another \'{e}tale groupoid and
$\phi:\HH\ra\GG$ is a continuous functor, then
$$ \Phi\com\phi_{\ast}=(\angs{\phi},
   (1_{\phi_{0}(b)},b)_{\ast}\com\Phi\;.$$
\Proof
Define a function $\bar{\Phi}:\Omega(\HH,b)\lra\Theta_{1}(\HH,b)$
as follows: Let
$$ \sigma_{n}\cdot h_{n}\cdot\ldots
   \cdot h_{1}\cdot\sigma_{0} $$
be an element of $\Omega(\HH,b)$. Put
$t_{i}^{n}=\frac{i}{n+1}$ and
let $J_{i}^{n}:[t_{i}^{n},t_{i+1}^{n}]\ra\III$ be the affine map with
$J_{i}^{n}(t_{i}^{n})=0$ and $J_{i}^{n}(t_{i+1}^{n})=1$,
for $0\leq i\leq n$.
Let $\tau_{i}=\sigma_{i}\com J_{i}^{n}$. Now for any
$1\leq i\leq n$ we define the equivariant homeomorphism
of principal $\HH$-bundles over the point $\{t_{i}^{n}\}$
$$ \alpha_{i}:\angs{\tau_{i-1}}|_{\{t_{i}^{n}\}}\lra
   \angs{\tau_{i}}|_{\{t_{i}^{n}\}} $$
by $\alpha_{i}(h,t_{i}^{n})=(h\com h_{i}^{-1},t_{i}^{n})$. With respect to
$(\alpha_{i})$ we can now amalgamate the bundles $\angs{\tau_{i}}$ into
a principal $\HH$-bundle $E$ over $\III$. In other words, we view
the $\HH$-loop $\sigma_{n}\cdot h_{n}\cdot\ldots
\cdot h_{1}\cdot\sigma_{0}$ as a closed $\HH$-cocycle on $\III$, and
we take $E$ to be the bundle associated to this cocycle
(see Example \ref{Vex1} (2)).
We can define $s:\rob\III\ra E$ by
$s(0)=(1_{b},0)\in\angs{\tau_{0}}\subset E$ and
$s(1)=(1_{b},1)\in\angs{\tau_{n}}\subset E$. Thus $(E,s)$ is a
marked principal $(\HH,b)$-bundle over $(\III,\rob\III)$,
and we take $\bar{\Phi}(\sigma_{n}\cdot h_{n}\cdot\ldots
\cdot h_{1}\cdot\sigma_{0})$ to be the isomorphic class of $(E,s)$.

Now it is obvious that the equivalent $\HH$-loops in $\Omega(\HH,b)$
gives homotopic Hilsum-Skandalis $1$-loops in $(\HH,b)$.
Moreover, a deformation between two $\HH$-loops 
$\sigma_{n}\cdot h_{n}\cdot\ldots\cdot h_{1}\cdot\sigma_{0}$ and
$\sigma'_{n}\cdot h'_{n}\cdot\ldots\cdot h'_{1}\cdot\sigma'_{0}$
in $\Omega(\HH,b)$ can be seen as a closed
$\HH$-cocycle on the finite cover $(\III\times[t_{i}^{n},t_{i+1}^{n}])$
of $\III^{2}$. The associated
principal $\HH$-bundle over $\III^{2}$ (which is naturally marked)
provides a homotopy between
$\bar{\Phi}(\sigma_{n}\cdot h_{n}\cdot\ldots
\cdot h_{1}\cdot\sigma_{0})$ and 
$\bar{\Phi}(\sigma'_{n}\cdot h'_{n}\cdot\ldots
\cdot h'_{1}\cdot\sigma'_{0})$.
Therefore $\bar{\Phi}$ induces a map
$$ \Phi:\pi_{1}(\HH,b)\lra\Pi_{1}(\HH,b) $$
which is clearly a homomorphism of groups. We shall prove that
it is an isomorphism.

Let $(E,s)$ be a Hilsum-Skandalis $n$-loop in $(\HH,b)$. Since
$E=(E,p,\xW)$ is principal and $\HH$ \'{e}tale, the map $\xW$ is a local
homeomorphism. Thus we can choose $n\geq 1$ and sections
$s_{i}:[t_{i}^{n},t_{i+1}^{n}]\ra E$ of $\xW$, for $0\leq i\leq n$.
We can choose this sections so that $s_{0}(0)=s(0)$ and
$s_{n}(1)=s(1)$. By composing these sections with $p$,
we get paths in $\HH_{0}$
$$ \sigma_{i}=p\com s_{i}\com (J^{n}_{i})^{-1}\;.$$
Moreover, for any $1\leq i\leq n$ there is a uniquely determined
$h_{i}$ such that
$s_{i}(t^{n}_{i})=h_{i}\cdot s_{i-1}(t^{n}_{i})$.
Now $\sigma_{n}\cdot h_{n}\cdot\ldots\cdot h_{1}\cdot\sigma_{0}\in
\Omega(\HH,b)$, and clearly
$$ \bar{\Phi}(\sigma_{n}\cdot h_{n}\cdot\ldots\cdot h_{1}
   \cdot\sigma_{0})=(E,s) $$
in $\Theta_{1}(\HH,b)$. This proves that $\Phi$ is surjective.

To prove that $\Phi$ is also injective, assume that
$\sigma_{n}\cdot h_{n}\cdot\ldots
\cdot h_{1}\cdot\sigma_{0}\in\Omega(\HH,b)$ is such that
$\bar{\Phi}(\sigma_{n}\cdot h_{n}\cdot\ldots
\cdot h_{1}\cdot\sigma_{0})$ represents the unit in $\Pi_{1}(\HH,b)$.
The unit in $\Pi_{1}(\HH,b)$ may be represented by the
marked principal $(\HH,a)$-bundle over $(\III,\rob\III)$
$$ (\dom^{-1}(b)\times\III,u)\;,$$
where $u(0)=(1_{b},0)$ and $u(1)=(1_{b},1)$.
Therefore there is a homotopy $\gH$ with
$\gH_{0}=\bar{\Phi}(\sigma_{n}\cdot h_{n}\cdot\ldots
\cdot h_{1}\cdot\sigma_{0})$ and $\gH_{1}=(\dom^{-1}(b)\times\III,u)$.
If $(E',s')$ is a marked principal $(\HH,b)$-bundle over
$(\III\times\III,\rob\III\times\III)$ which represents $\gH$, we can
find $k\geq 1$ and a multiple $m$ of $n$ big enough so that $E$ can be
represented by an $\HH$-cocycle $c=(c_{(i,j)(i',j')})$ defined on the
finite closed cover of $\III\times\III$
$$ ([t^{m}_{i},t^{m}_{i+1}]\times
   [t^{k}_{j},t^{k}_{j+1}])_{0\leq i\leq m,\,0\leq j\leq k}\;.$$
Moreover, we choose $c$ such that its
restriction on $\III\times\{0\}$ (viewed as an $\HH$-loop) is
equivalent to $\sigma_{n}\cdot h_{n}\cdot\ldots
\cdot h_{1}\cdot\sigma_{0}$,
and its restriction on the rest of the boundary of $\III^{2}$
is constant with the value $b$.
Hence $c$ provides a chain of deformations and equivalences between
$\sigma_{n}\cdot h_{n}\cdot\ldots\cdot h_{1}\cdot\sigma_{0}$ 
and the $\HH$-loop
$$ \varpi(b)\cdot 1_{b}\cdot\ldots\cdot 1_{b}\cdot\varpi(b)\;,$$
which obviously represents the unit in $\pi_{1}(\HH,b)$.
Here $\varpi(b)$ denotes the constant path in $\HH_{0}$ with the image
point $b$. This completes the proof.
\eop
\Rem
Thus an $\HH$-loop in $\Omega(\HH,b)$ can be viewed as a
$\GG$-cocycle on a closed finite cover of $\III$ which represents
a Hilsum-Skandalis $1$-loop in $(\HH,b)$, while equivalences
and deformations between the elements of $\Omega(\HH,b)$ identify
the cocycles which represents the same homotopy class of
Hilsum-Skandalis $1$-loops in $(\HH,b)$. In other words, we have
a description of $\Pi_{1}(\HH,b)$ by $\HH$-cocycles on closed finite
covers. One can easily generalize this theorem by representing
a Hilsum-Skandalis $n$-loop in $(\HH,b)$ by an $\HH$-cocycle
on a closed finite cover of $\III^{n}$ which is constant
(with value $b$) when restricted on $\rob\III^{n}$,
and by representing a homotopy by an $\HH$-cocycle on a closed finite
cover of $\III^{n}\times\III$ which is constant (with value $b$)
when restricted on $\rob\III^{n}\times\III$.

Theorem \ref{Vtheo6} implies that $\Pi_{1}(\HH,b)$
coincides with the fundamental group of $\HH$ as described by
W. T. van Est \cite{Est} if $\HH$ is effective and $\HH_{0}$
simply-connected, as well as with the fundamental group
of the classifying space and of the
classifying topos of $\HH$ \cite{Moe4,Moe1}.
In fact, the results in \cite{Moe4} indicate that
the same is true for all the higher
homotopy groups as well.

\begin{cor}  \label{Vcor7}
Let $\phi:\HH\ra\GG$ be an essential equivalence between \'{e}tale
groupoids, and let $b\in\HH_{0}$. Then
$$ \phi_{\ast}:\pi_{1}(\HH,b)\lra\pi_{1}(\GG,\phi_{0}(b)) $$
is an isomorphism.
\end{cor}

\section{Singular Homology of Topological Groupoids}
\label{secSinHom}

In this Section we introduce a homology theory
of topological groupoids which we call the singular homology.
This homology is the classical singular
homology when restricted to the topological spaces.
Moreover, the general definition is literally the
same as the definition of the singular homology of topological spaces
if one replace the category of topological spaces with that
of Hilsum-Skandalis maps.
Since the singular homology is a functor on the category $\cGpd$,
it is invariant under Morita equivalence.

Let $\GG$ be an \'{e}tale groupoid.
Then there is a chain homomorphism
$$ \omega:S(\GG_{1})\lra S(\GG_{0}) $$
between the singular chain complexes of the topological spaces
$\GG_{1}$ and $\GG_{0}$, given by
$$ \omega=\cod_{\sh}-\dom_{\sh}\;.$$
We shall prove that the singular homology groups of $\GG$ are
exactly the homology groups of the cokernel of $\omega$.
Finally, we shall prove that the effect-functor induces isomorphisms
between the singular homology groups of $\GG$ and $\Eff(\GG)$.

\subsection{Definition of Singular Homology}  \label{secDefSinHom}

For each $n\geq 0$ denote by $\ssx{n}$ the standard $n$-simplex,
i.e. the convex hull
of the standard basis in $\RRR^{n+1}$. For $n\geq 1$ and
$0\leq i\leq n$, denote by $d_{i}=d_{i}^{n}:\ssx{n-1}\ra\ssx{n}$
the affine map given on the basis by
$$ d_{i}(e_{k})=\left\{ \begin{array}{ccl}
                        e_{k}   & ; & k<i \\
                        e_{k+1} & ; & k\geq i
                        \end{array}
                \right. \;\;\;\;\;\;\;\;\;\; 0\leq k\leq n-1\;.$$

A {\em singular $n$-simplex} in a topological groupoid $\GG$ is a
Hilsum-Skandalis map $u\in\cGpd(\ssx{n},\GG)$.
For any $n\geq 1$ and $0\leq i\leq n$ we have the map
$$ \rob^{i}=\rob^{i}_{n}=\cGpd(d_{i},\GG):
   \cGpd(\ssx{n},\GG)\lra\cGpd(\ssx{n-1},\GG)\;.$$
It maps a singular $n$-simplex $u$ to its $i$-th face
$\rob^{i}u=d_{i}^{\ast}u=u\otimes d_{i}$.

The {\em group of $n$-chains} $\gS_{n}(\GG)$ in $\GG$ is the free
abelian group generated by the set of singular $n$-simplexes in
$\GG$,
$$ \gS_{n}(\GG)=\ZZZ\cGpd(\ssx{n},\GG)\;.$$
We define the boundary homomorphism
$$ \rob=\rob_{n}:\gS_{n}(\GG)\lra\gS_{n-1}(\GG) $$
on the generators by
$$ \rob u=\sum_{i=0}^{n}(-1)^{i}\rob^{i}u\;.$$
Since the face maps satisfy the equation
$d_{i}\com d_{j}=d_{j+1}\com d_{i}$
if $i\leq j$, we have $\rob^{2}=0$. Hence
$\gS(\GG)=(\gS_{n}(\GG),\rob)$
is a chain complex of abelian groups.

Let $E\in\cGpd(\HH,\GG)$. The maps
$$ \cGpd(\ssx{n},E):\cGpd(\ssx{n},\HH)\lra\cGpd(\ssx{n},\GG) $$
extend to a chain homomorphism
$$ E_{\sh}=\gS(E):\gS(\HH)\lra\gS(\GG)\;,$$
and this gives a functor $\gS$ from $\cGpd$ to the category
of (non-negative) chain complexes of abelian groups.

Let $\GG$ be a topological groupoid, $A$ an abelian group and
$n\geq 0$. Define the {\em $n$-th singular homology group}
$H_{n}(\GG;A)$ of the topological groupoid $\GG$ with
coefficients in $A$ to be the $n$-th homology group
of the chain complex $\gS(\GG)\otimes A$,
$$ H_{n}(\GG;A)=H_{n}(\gS(\GG)\otimes A) \;\;\;\;\;\; n\geq 0\;.$$
Further, we define the {\em $n$-th singular cohomology group} 
$H^{n}(\GG;A)$ of the topological groupoid $\GG$ with
coefficients in $A$ to be the $n$-th cohomology group
of the chain cocomplex $\Hom(\gS(\GG),A)=(\Hom(\gS_{n}(\GG),A),\delta)$,
$$ H^{n}(\GG;A)=H^{n}(\Hom(\gS(\GG),A)) \;\;\;\;\;\; n\geq 0\;.$$

Let $E\in\cGpd(\HH,\GG)$. Then denote
$$ E_{\ast}=H_{n}(E;A)=H_{n}(E_{\sh}\otimes A):
   H_{n}(\HH;A)\lra H_{n}(\GG;A)\;,$$
and
$$ E^{\ast}=H^{n}(E;A)=H^{n}(\Hom(E_{\sh},A)):
   H^{n}(\GG;A)\lra H^{n}(\HH;A)\;.$$
If $\phi:\HH\ra\GG$ is a continuous functor, we denote
$\phi_{\ast}=\angs{\phi}_{\ast}$ and
$\phi^{\ast}=\angs{\phi}^{\ast}$.
Further, if $f:A\ra B$ is a homomorphism of abelian groups,
write
$$ H_{n}(\GG;f)=H_{n}(\gS(\GG)\otimes f):
   H_{n}(\GG;A)\lra H_{n}(\GG;B)\;,$$
and
$$ H^{n}(\GG;f)=H^{n}(\Hom(\gS(\GG),f)):
   H^{n}(\GG;A)\lra H^{n}(\GG;B)\;.$$
This gives the sequences of functors
$$ H_{n}:\cGpd\times\Ab\lra\Ab\;\;\;\;\;\;\;\; n\geq 0 $$
and
$$ H^{n}:\cGpd^{op}\times\Ab\lra\Ab\;\;\;\;\;\;\;\; n\geq 0\;.$$

\begin{ex}  \label{Vex10}  \rm
(1) Let $X$ be a topological space. Since the category
$\Top$ is a full subcategory of $\cGpd$, the
chain complex $\gS(X)$ is exactly the singular chain complex 
$S(X)$ of the topological space
$X$. The singular homology and cohomology groups of the \'{e}tale
groupoid $X$ are thus just the singular homology and cohomology
groups of the topological space $X$.

(2) Let $G$ be a discrete group. In paricular, $G$ is an \'{e}tale
groupoid with $\GG_{0}$ an one-point space.
Any principal $G$-bundle over $\ssx{n}$ is isomorphic to the
trivial bundle $G\times\ssx{n}$ since $\ssx{n}$ is simply
connected. Therefore $\cGpd(\ssx{n},G)$ has exactly one
element, and
$$ H_{n}(G;A)=H^{n}(G;A)
              =\left\{ \begin{array}{ccl}
                      A  & ; & n=0 \\
                      0  & ; & n\geq 1
                        \end{array}
                \right.$$
for any abelian group $A$.
\end{ex}

\begin{theo}  \label{Vtheo11}
If $\GG$ and $\HH$ are Morita equivalent topological groupoids,
then
$$ H_{n}(\GG;A)\cong H_{n}(\HH;A) $$
and
$$ H^{n}(\GG;A)\cong H^{n}(\HH;A)\;,$$
for any abelian grop $A$ and any $n\geq 0$.
\end{theo}
\Proof
By Corollary \ref{IIcor7} the groupoids $\GG$ and $\HH$ are
isomorphic in $\cGpd$. The theorem thus follows from the
functoriality of $H_{n}$ and $H^{n}$.
\eop
 
\begin{prop} \label{Vprop12}
Let $\GG$ be a topological groupoid.
A short exact sequence of abelian groups
$$\CD
  0 \cdr{}{} A \cdr{}{} B \cdr{}{} C \cdr{}{} 0
  \endCD$$
induces the long exact sequences
$$ \ldots \lra                    H_{n}(\GG;A) 
       \lra                       H_{n}(\GG;B)
       \lra                       H_{n}(\GG;C) 
       \lra                       H_{n-1}(\GG;A)
       \lra \ldots $$
and
$$ \ldots \lra                    H^{n}(\GG;A) 
       \lra                       H^{n}(\GG;B)
       \lra                       H^{n}(\GG;C)
       \lra                       H^{n+1}(\GG;A)
       \lra \ldots $$
\end{prop}
\Proof
This follows from the fact that $\gS(\GG)$ is a free complex.
\eop

\begin{ex}  \label{Vex12a}  \rm
Let $X$ be a topological space equipped with a properly
discontinuous action of a discrete group $G$. In particular,
$$ p:X\lra X/G $$
is a covering projection. The groupoid $G(X)$ is \'{e}tale
and Morita equivalent to the space $X/G$.
Theorem \ref{Vtheo11} thus implies that
$$ H_{n}(G(X);A)\cong H_{n}(X/G;A)\;\;\;\;\;\;\mbox{and}
   \;\;\;\;\;\; H^{n}(G(X);A)\cong H^{n}(X/G;A)\;,$$
for any abelian group $A$ and $n\geq 0$.
\end{ex}

\subsection{Singular Homology of \'{E}tale Groupoids}

In this Subsection we give a simpler description
of the singular homology groups in the case where the
groupoid is \'{e}tale.

Let $\GG$ be an \'{e}tale groupoid. The domain and the codomain
map $\dom,\cod:\GG_{1}\ra\GG_{0}$ of $\GG$ induce the chain
homomorphisms
$$ \dom_{\sh},\cod_{\sh}:S(\GG_{1})\lra S(\GG_{0}) $$
between the singular chain complexes of the topological spaces
$\GG_{0}$ and $\GG_{1}$. Put
$$ \omega=\cod_{\sh}-\dom_{\sh}\;.$$
We denote by $S(\GG)=(S_{n}(\GG),\rob)$ the cokernel
of $\omega$,
$$\CD
  S(\GG_{1}) \cdr{\omega}{}      S(\GG_{0})
             \cdr{\epsilon}{}    S(\GG)      \cdr{}{}  0
  \endCD\;.$$
If $\phi:\HH\ra\GG$ is a continuous functor between \'{e}tale groupoids,
it induces a chain homomorphism
$S(\phi):S(\HH)\ra S(\GG)$, as
in the following diagram:
$$\CD
  S(\HH_{1}) \cdr{\omega}{}    S(\HH_{0})
             \cdr{\epsilon}{}    S(\HH)       \cdr{}{}  0         \\
  \cdd{S(\phi_{1})}{}  \cdd{}{S(\phi_{0})}  \cdd{}{S(\phi)}  \cd. \\
  S(\GG_{1}) \cdr{\omega}{} S(\GG_{0})
             \cdr{\epsilon}{}    S(\GG)       \cdr{}{} 0
  \endCD$$
This gives a functor $S$ from the category $\Gpde$ of \'{e}tale groupoids
and continuous functors to the category of (non-negative)
chain complexes of abelian groups. It extends the usual functor $S$
defined on the full subcategory $\Top$ of topological spaces, i.e.
the notation is justified.
Note that $S(\uuu)=\epsilon$, where $\uuu:\GG_{0}\ra\GG$ is the
canonical functor.

Let $\GG$ be an \'{e}tale groupoid.
We define the {\em balanced chain complex} $BS(\GG)$ of $\GG$ to be
the kernel of the chain homomorphism $\epsilon$, so we have a short
exact sequence
\begin{equation}  \label{Vshexseq}
\CD
  0 \cdr{}{} BS(\GG) \cdr{\iota}{} S(\GG_{0})
    \cdr{\epsilon}{} S(\GG) \cdr{}{} 0
  \endCD\;.
\end{equation}
The chain complex $BS(\GG)$ is free.
For an abelian group $A$, the homology of the chain complex 
$BS(\GG)\otimes A$ will be called the
{\em balanced homology} of $\GG$ with coefficients in $A$,
and denoted by
$$ BH_{n}(\GG;A)=H_{n}(BS(\GG)\otimes A)\;\;\;\;\;\;\;\;n\geq 0\;.$$
Analogously we define the {\em balanced cohomology} of $\GG$ with
coefficients in $A$ as $BH^{n}(\GG;A)=H^{n}(\Hom(BS(\GG),A))$.

Let us describe $BS_{n}(\GG)$ explicitly. We say that two
singular $n$-simplexes $u,v:\ssx{n}\ra\GG_{0}$ are {\em similar} or
{\em $\GG$-similar} if there exists a singular $n$-simplex
$f:\ssx{n}\ra\GG_{1}$
such that $\dom\com f=u$ and $\cod\com f=v$.
This is clearly an equivalence relation, which will be denoted by
$\sml$ or $\sml_{\GG}$.
Since $BS(\GG)$ is also the image of $\omega$, it is
clear that $BS_{n}(\GG)$ is generated by the elements
$$ u-v\in S_{n}(\GG_{0})\;,$$
where $u,v\in\Top(\ssx{n},\GG_{0})$ are similar.
Furthermore, we have a natural isomorphism
$$ S_{n}(\GG)\cong\ZZZ\{\Top(\ssx{n},\GG_{0})/\!\!\sml\}\;.$$
In particular, the chain complex $S(\GG)$ is free and hence
the short exact sequence (\ref{Vshexseq}) splits.

Observe that the singular $n$-simplexes in $\GG_{0}$
can be identified with
the continuous functors in $\Gpd(\ssx{n},\GG)$.
Let $u$ be a singular $n$-simplex in $\GG_{0}$.
Then $\angs{\uuu\com u}$ is a principal $\GG$-bundle over $\ssx{n}$.
Let $v$ be a singular $n$-simplex in $\GG_{0}$ similar
to $u$, and let $f$ be a singular $n$-simplex in $\GG_{1}$ such
that $\dom\com f=u$ and $\cod\com f=v$. Then there is a
$\GG$-equivariant map $\alpha:\angs{\uuu\com v}\ra\angs{\uuu\com u}$
given by
$$ \alpha(g,x)=(g\cdot f(x),x)\;,$$
hence $\angs{\uuu\com v}$ and $\angs{\uuu\com u}$ are isomorphic.
Conversely, if $\alpha:\angs{\uuu\com v}\lra\angs{\uuu\com u}$
is a $\GG$-equivariant map, there is a unique singular
$n$-simplex $f$ in $\GG_{1}$ such that
$$ f(x)\cdot (1_{u(x)},x)=\alpha(1_{v(x)},x)\;,$$
therefore $u\sml v$. In other words, we proved that
there is natural bijective correspondence
$$ \Top(\ssx{n},\GG_{0})/\!\!\sml\;\cong
   \langle\Gpd(\ssx{n},\GG)\rangle\subset\cGpd(\ssx{n},\GG)\;.$$
In this view, we will consider $S(\GG)$ as a subcomplex of
$\gS(\GG)$. A singular $n$-simplex $u=(E,p,\xW)$ in $\GG$ is in
$S(\GG)$ if and only if there exists a global section of $\xW$.
Moreover, if $\phi:\HH\ra\GG$ is a continuous functor,
we have $\gS(\angs{\phi})|_{S(\HH)}=S(\phi)$.

\begin{theo}  \label{Vtheo13}
Let $\GG$ be an \'{e}tale groupoid. The chain inclusion
$S(\GG)\hookrightarrow\gS(\GG)$ is a chain equivalence.
In particular,
$H_{n}(\GG;A)\cong H_{n}(S(\GG)\otimes A)$ and
$H^{n}(\GG;A)\cong H^{n}(\Hom(S(\GG),A))$, for any abelian group
$A$ and any $n\geq 0$.
\end{theo}
\Proof
First note that we can apply the classical construction of the
barycentric subdivision \cite{Spa} to obtain a chain map
$$ \gsd:\gS(\GG)\lra\gS(\GG) $$
and a chain homotopy
$$ D=D_{n}:\gS_{n}(\GG)\lra\gS_{n+1}(\GG) $$
such that $\rob\com D+ D\com\rob=1-\gsd$. The map $\gsd$ sends a
singular $n$-simplex $u$ into a sum of the restrictions of $u$
on the $n$-simplexes of the barycentric subdivision of $\ssx{n}$,
with the appropriate signs. Observe that both $\gsd$ and $D$ are
natural in $\GG$ and restrict to the subcomplex $S(\GG)$.

Let $u$ be a singular $n$-simplex in $\GG$, represented by
a principal $\GG$-bundle $(E,p,\xW)$ over $\ssx{n}$.
Since $\xW$ is a local homeomorphism and $\ssx{n}$ is compact,
there exists a finite open cover $\cU=(U_{i})$ of $\ssx{n}$ with
sections $s_{i}:U_{i}\ra E$ of $\xW$. Now we
can barycentrically subdivide $\ssx{n}$ sufficiently many times
such that each simplex of the subdivision lies in an element of
$\cU$. The restriction of $E$ on such a simplex is therefore in
$S_{n}(\GG)$. This implies that for any $\xi\in\gS_{n}(\GG)$
there exists $j\geq 0$ such that
$$ \gsd^{j}\xi\in S_{n}(\GG)\;.$$
The standard argument (see for example \cite{Spa}) now
implies that the inclusion $S(\GG)\hookrightarrow\gS(\GG)$
is a chain equivalence.
\eop

\begin{prop}  \label{Vcor14}
Let $\GG$ be an \'{e}tale groupoid.
There are long exact sequences
$$ \ldots \lra                          BH_{n}(\GG;A) 
          \stackrel{\iota_{\ast}}{\lra} H_{n}(\GG_{0};A)
          \stackrel{\uuu_{\ast}}{\lra}  H_{n}(\GG;A) 
          \lra                          BH_{n-1}(\GG;A)
          \lra \ldots $$
and
$$ \ldots \lra                          H^{n}(\GG;A)
          \stackrel{\uuu^{\ast}}{\lra}  H^{n}(\GG_{0};A) 
          \stackrel{\iota^{\ast}}{\lra} BH^{n}(\GG;A)
          \lra                          H^{n+1}(\GG;A) 
          \lra \ldots $$
\end{prop}
\Proof
This follows from Theorem \ref{Vtheo13} and the short
exact sequence (\ref{Vshexseq}).
\eop

\begin{prop}[Mayer-Vietoris sequence]  \label{Vprop14a}
Let $\GG$ be an \'{e}tale gro\-up\-oid and $A$ an abelian group.
Let $U$ and $V$ be open $\GG$-invariant subsets of $\GG_{0}$
such that $U\cup V=\GG_{0}$. Then there are long exact
sequences
\begin{eqnarray*}
   \ldots \lra  H_{n}(\GG|_{U\cap V};A) 
          \lra  H_{n}(\GG|_{U};A)\oplus H_{n}(\GG|_{V};A)
          \lra  H_{n}(\GG;A)                              \\
   \lefteqn{\lra  H_{n-1}(\GG|_{U\cap V};A)
            \lra  \ldots} \hspace{111 mm}
\end{eqnarray*}
and
\begin{eqnarray*}
   \ldots \lra  H^{n}(\GG;A) 
          \lra  H^{n}(\GG|_{U};A)\oplus H^{n}(\GG|_{V};A)
          \lra  H^{n}(\GG|_{U\cap V};A)                   \\
   \lefteqn{\lra  H^{n+1}(\GG;A)
            \lra  \ldots} \hspace{112 mm} 
\end{eqnarray*}
\end{prop}
\Proof
We will adapt the standard proof to our case.
Write $S^{\{U,V\}}(\GG_{0})$ for the subcomplex of
$S(\GG_{0})$ generated by those simplexes in $\GG_{0}$
which lie in $U$ or in $V$. Put
$U'=(\GG|_{U})_{1}$ and $V'=(\GG|_{V})_{1}$, and denote by
$S^{\{U',V'\}}(\GG_{1})$ the subcomplex of $S(\GG_{1})$
generated by the simplexes in $\GG_{1}$ which lie in $U'$ or
in $V'$. Observe that the chain homomorphism $\omega$
restricts to a map
$S^{\{U',V'\}}(\GG_{1})\ra S^{\{U,V\}}(\GG_{0})$.
Denote by $S^{\{U,V\}}(\GG)$ the cokernel of this restriction,
$$\CD
  S^{\{U',V'\}}(\GG_{1}) \cdr{}{}
  S^{\{U,V\}}(\GG_{0}) \cdr{}{}
  S^{\{U,V\}}(\GG) \cdr{}{} 0
  \endCD\;.$$
One can easily see that $S^{\{U,V\}}(\GG)$ is a subcomplex
of $S(\GG)$ generated by the similarity classes of simplexes
in $\GG_{0}$ which lie in $U$ or in $V$.
Moreover, the sequence
$$\CD
  0 \cdr{}{} S(\GG|_{U\cap V})
  \cdr{(\iii_{\sh},-\jjj_{\sh})}{}
  S(\GG|_{U})\oplus S(\GG|_{V})
  \cdr{\iii'_{\sh}+\,\jjj'_{\sh}}{}
  S^{\{U,V\}}(\GG)
  \cdr{}{} 0
  \endCD$$
is exact and splits. Here $\iii$, $\jjj$, $\iii'$ and
$\jjj'$ are the obvious inclusions of \'{e}tale groupoids.
Using the barycentric subdivision $\gsd$ as in the proof of Theorem
\ref{Vtheo13} we can see that the inclusion of $S^{\{U,V\}}(\GG)$
into $S(\GG)$ is a chain equivalence.
\eop

\begin{theo} \label{Vtheo15}
Let $\phi:\HH\ra\GG$ be a continuous functor between \'{e}tale groupoids
which is an isomorphism on objects and surjective on morphisms. Then
$\phi$ induces isomorphisms
$$ H_{n}(\HH;A)\lra H_{n}(\GG;A) $$
and
$$ H^{n}(\GG;A)\lra H^{n}(\HH;A)\;,$$
for any abelian group $A$ and $n\geq 0$.
\end{theo}
\Proof
First we can identify $\HH_{0}$ and $\GG_{0}$. Since
$\epsilon=S(\phi)\com\epsilon$, we have
$BS_{n}(\HH)\subset BS_{n}(\GG)$. Now because of the long
exact sequences in Corollary \ref{Vcor14} it is enough to
prove that the inclusion
$$ BS(\HH)\lra BS(\GG) $$
is a chain equivalence.

As in the proof of Theorem \ref{Vtheo13} we will use the
barycentric subdivision. Because of the naturality of the
chain homomorphism $\gsd$ and the chain
homotopy $D$ both of these maps restrict from $S(\GG_{0})$
on $BS(\GG)$ and also on $BS(\HH)$. In particular, the
restrictions of $\gsd$ on $BS(\GG)$ and $BS(\HH)$ are chain
equivalences.

Let $u$ and $v$ be similar singular $n$-simplexes in $\GG_{0}$.
Hence there exists $f:\ssx{n}\ra\GG_{1}$ such that $\dom\com f=u$
and $\cod\com f=v$. Now $\phi:\HH_{1}\ra\GG_{1}$ is a surjective
local homeomorphism, hence we can choose a finite open cover
$(U_{i})$ of $f(\ssx{n})$ and sections $s_{i}:U_{i}\ra\HH_{1}$
of $\phi$. We can barycentrically subdivide
$\ssx{n}$ sufficiently many times so that each simplex
in the subdivision lies in an element of the cover $(f^{-1}(U_{i}))$
of $\ssx{n}$. This implies that the restrictions
of $u$ and $v$ on a simplex in the subdivision are $\HH$-similar.
In other words, there exists $j\geq 0$ such that
$$ \gsd^{j}(u-v)\in BS_{n}(\HH)\;.$$
Since $BS(\GG)$ is free, the standard argument yields that
$BS(\HH)\ra BS(\GG)$ is a chain equivalence.
\eop

\begin{cor}  \label{Vcor16}
Let $\GG$ be an \'{e}tale groupoid. Then
$$ \eee_{\ast}:H_{n}(\GG;A)\lra H_{n}(\Eff(\GG);A) $$
and
$$ \eee^{\ast}:H^{n}(\Eff(\GG);A)\lra H^{n}(\GG;A) $$
are isomorphisms, for any abelian group $A$ and $n\geq 0$.
\end{cor}

\begin{ex}  \label{Vex17}  \rm
(1) Let $\GG$ be an \'{e}tale $\eCe$-groupoid of dimension $q$.
There exists a Morita equivalent groupoid $\HH$ such that
$\HH_{0}=\coprod_{i\in I}\RRR^{q}$.
Using Theorem \ref{Vtheo13} and the long exact sequence from
Corollary \ref{Vcor14} we get
$$ H_{n+1}(\GG;A)\cong BH_{n}(\HH;A)\;\;\;\;\;\;\mbox{and}
   \;\;\;\;\;\; H^{n+1}(\GG;A)\cong BH^{n}(\HH;A)\;,$$
for any abelian group $A$ and any $n\geq 1$.

(2) Let $\GG$ be an \'{e}tale groupoid with $\GG_{0}$ locally
path-connected. Denote by $\rrr:\GG\ra|\GG|$ the canonical
continuous functor. Then
$\rrr_{\ast}:H_{0}(\GG;\ZZZ)\ra H_{0}(|\GG|;\ZZZ)$
is clearly surjective.
The group $H_{0}(|\GG|;\ZZZ)$
is the free abelian group on the set
of $\GG$-connected components of $\GG_{0}$.
For $a\in\GG_{0}$ denote by $\varpi(a)$ the $0$-simplex
in $\GG_{0}$ with the image point $a$.

If $a$ and $a'$ are points in the same $\GG$-connected
component of $\GG_{0}$, there exists a $\GG$-path in $\GG_{0}$
$$ \sigma_{n}\cdot g_{n}\cdot\ldots\cdot g_{1}\cdot\sigma_{0} $$
from $a$ to $a'$, see Proposition \ref{Iprop11}. Now
$$ \rob(\epsilon(\sigma_{0}+\ldots+\sigma_{n}))=
   \epsilon\varpi(a')-\epsilon\varpi(a)\;.$$
This yields that $\rrr_{\ast}$ is injective. Thus
$$ \rrr_{\ast}:H_{0}(\GG;\ZZZ)\lra H_{0}(|\GG|;\ZZZ) $$
is an isomorphism.

(3) Let $G$ be a subgroup of $\RRR^{n}$. Put the discrete topology
on $G$, and let $G$ acts on $\RRR^{n}$ by translations, i.e.
$y\cdot g=y+g$, for any $y\in\RRR^{n}$ and $g\in G$. Define
a homomorphism of groups $f:S_{0}(\RRR^{n})\ra\RRR^{n}$ on generators
by
$$ f(\varpi(y))=y\;\;\;\;\;\;\;\; y\in\RRR^{n}\;.$$
Since $G(\RRR^{n})_{0}=\RRR^{n}$, we have $BS_{0}(G(\RRR^{n}))\subset
S_{0}(\RRR^{n})$. Observe that
$$ f(BS_{0}(G(\RRR^{n})))=G\;.$$
Moreover, if $\sigma$ and $\sigma'$ are $G(\RRR^{n})$-similar
$1$-simplexes in $\RRR^{n}$, there exists $g\in G$ such that
$\sigma(x)=\sigma'(x)+g$. Hence we have
$$ f(\rob(\sigma-\sigma'))=0\;,$$
so $f$ induces a surjective homomorphism of groups
$$ \alpha:BH_{0}(G(\RRR^{n});\ZZZ)\lra G\;.$$
We will show that $\alpha$ is in fact an isomorphism.

For any $g\in G$, put
$$ \tau(g)=\varpi(g)-\varpi(0)\in BS_{0}(G(\RRR^{n}))\;.$$
Note that $f(\tau(g))=g$. The group $BS_{0}(G(\RRR^{n}))$ is
generated by the elements of the form
$$ \varpi(y+g)-\varpi(y) $$
with $y\in\RRR^{n}$ and $g\in G$. For $y,y'\in\RRR^{n}$, let
$\sigma(y',y):\III\ra\RRR^{n}$ be the singular $1$-simplex in
$\RRR^{n}$ given by
$\sigma(y',y)(t)=ty'+(1-t)y$. For any $y\in\RRR^{n}$
and $g\in G$ we have $\sigma(y+g,g)-\sigma(y,0)\in BS_{1}(G(\RRR^{n}))$
and
$$ \rob(\sigma(y+g,g)-\sigma(y,0))=
   (\varpi(y+g)-\varpi(y))-\tau(g)\;.$$
Hence $BH_{0}(G(\RRR^{n});\ZZZ)$
is generated by the homology classes of
$\tau(g)$, $g\in G$.

Next observe that for any $g,g'\in G$ we have
$\sigma(g+g',g')-\sigma(g,0)\in BS_{1}(G(\RRR^{n}))$ and
$$ \rob(\sigma(g+g',g')-\sigma(g,0))=
   \tau(g+g')-\tau(g)-\tau(g')\;.$$
This implies that $\alpha$ is an isomorphism.

Finally, from (2) we know that the map
$\uuu_{\ast}:H_{0}(\RRR^{n};\ZZZ)\ra H_{0}(G(\RRR^{n});\ZZZ)$
is an isomorphism. The
long exact sequence in Corollary \ref{Vcor14} thus implies
$$ H_{1}(G(\RRR^{n});\ZZZ)\cong BH_{0}(G(\RRR^{n});\ZZZ)\cong G\;.$$

(4) Let $\GG$ be an \'{e}tale groupoid such that the quotient map
$$ \GG_{0}\lra |\GG| $$
is a local homeomorphism. Then we can define an \'{e}tale groupoid
$\HH$ with $\HH_{0}=\GG_{0}$ and
$\HH_{1}=\GG_{0}\times_{|\GG|}\GG_{0}$, and
with the obvious structure
maps. The map $\GG_{0}\ra |\GG|$ can now be extended to a
continuous functor $\psi:\HH\ra |\GG|$ which is an essential
equivalence.

Observe that the canonical functor $\rrr:\GG\ra |\GG|$ factors
through $\psi$ as $\rrr=\psi\com\phi$, where $\phi:\GG\ra\HH$
is given as the identity on objects and by
$$ \phi(g)=(\cod g,\dom g) $$
for any $g\in\GG_{1}$. Since $\phi$ is surjective on morphisms,
Theorem \ref{Vtheo15} implies that it induces isomorphisms on
homology and cohomology. We can thus conclude that the
functor $\rrr$ induces isomorphisms
$$ H_{n}(\GG;A)\lra H_{n}(|\GG|;A) $$
and
$$ H^{n}(|\GG|;A)\lra H^{n}(\GG;A)\;,$$
for any abelian group $A$ and any $n\geq 0$.

(5) Let $\cF$ be a foliation on a finite-dimensional
$\eCe$-manifold $M$. For any $n\geq 0$ we define the
singular homology groups $H_{n}((M,\cF);\ZZZ)$ of the
transverse structure of $\cF$ as the singular homology
groups of the holonomy groupoid $\Hol(M,\cF)$,
$$ H_{n}((M,\cF);\ZZZ)=H_{n}(\Hol(M,\cF);\ZZZ)\;.$$
By Theorem \ref{Vtheo11} we have
$H_{n}((M,\cF);\ZZZ)\cong H_{n}(\Hol_{T}(M,\cF);\ZZZ)$
for any complete transversal $T:N\ra M$ of $\cF$.
Moreover, the principal $\Hol_{T}(M,\cF)$-bundle over $M$
which corresponds to $\cF$ (Example \ref{Iex9} (2))
induces a homomorphism
$$ H_{n}(M;\ZZZ)\lra H_{n}((M,\cF);\ZZZ)\;,$$
for any $n\geq 0$.
\end{ex}

\section{Bimodule Associated to a Principal Bibundle}

The aim of this section is to associate to
a Hilsum-Skandalis $\eCr$-map $E$ between object-separated
finite-dimensional \'{e}tale $\eCr$-groupoids $\HH$
and $\GG$ an isomorphism class of a $\eCrc(\GG)$-$\eCrc(\HH)$-bimodule 
$\eCrc(E)$. Here $\eCrc(\GG)$ denotes the convolution algebra of
$\eCr$-functions on $\GG_{1}$ as introduced by A. Connes
\cite{BryNis,Con78,Con82,Ren}. We show that $\eCrc$ is
a functor from the category of Hilsum-Skandalis $\eCr$-maps
between separated finite-dimensional \'{e}tale
$\eCr$-groupoids to the category
of locally unital
bimodules over the algebras with local units. In particular, this
implies that the $\eCr$-Morita equivalent separated finite-dimensional
\'{e}tale $\eCr$-groupoids have Morita equivalent algebras -- a
result announced in a more general form in \cite{BryNis}.
\vspace{4 mm}

\Not (1) A topological groupoid $\GG$ is called {\em object-separated}
if $\GG_{0}$ is Hausdorff. It is called {\em separated}
(or {\em Hausdorff}) if $\GG_{1}$ (and hence also $\GG_{0}$)
is Hausdorff. We shall denote by
$\cGpdbr$ the full subcategory of $\cGpder$ with objects all the
separated finite-dimensional \'{e}tale $\eCr$-groupoids.
Observe that if $(E,p,\xW)\in\cGpdbr(\HH,\GG)$, then 
the space $E$ is Hausdorff.

(2) Recall that an (associative complex) algebra $\cA$ has local
units if for any $x_{1},x_{2},\ldots,x_{n}\in\cA$ there exists
$x\in\cA$ such that $x_{i}x=x\,x_{i}=x_{i}$, $1\leq i\leq n$.
If $\cA$ and $\cB$ are algebras with local units, an
$\cA$-$\cB$-bimodule $\cM$ is called locally unital if for any
$m_{1},m_{2},\ldots,m_{n}\in\cM$ there exist $x\in\cA$
and $y\in\cB$
such that $x\cdot m_{i}=m_{i}\cdot y=m_{i}$, $1\leq i\leq n$.
If $\cC$ is another algebra with local units,
$\cM$ a locally unital $\cA$-$\cB$-bimodule
and $\cN$ a locally
unital $\cB$-$\cC$-bimodule, then
$$ \cM\ten{\cB}\cN $$
is a locally unital $\cA$-$\cC$-bimodule. Moreover, $\cB$ itself
can be viewed as a locally unital $\cB$-$\cB$-bimodule, and there
are canonical isomorphisms of bimodules
$$ \cM\ten{\cB}\cB\cong\cM $$
and
$$ \cB\ten{\cB}\cN\cong\cN\;.$$
Hence we have a category $\cAlg$ with the associative complex
algebras with local units as objects, with the isomorphism classes of
locally unital $\cA$-$\cB$-bimodules as morphisms in
$\cAlg(\cB,\cA)$, and with the composition induced by the
tensor product. We shall denote the isomorphism class of a
locally unital $\cA$-$\cB$-bimodule $\cM$ again by $\cM$.

(3) For a Hausdorff finite-dimensional
$\eCr$-manifold $M$, denote by $\eCrc(M)$ the
vector space of complex $\eCr$-functions on $M$ with compact support.

(4) Let $M$ be a finite-dimensional $\eCr$-manifold, not necessarily
Hausdorff. If $U\subset M$ is (the domain of) a chart for $M$ and
$f'\in\eCrc(U)$, we will denote by $f'^{M}$ the extension of $f'$ on
$M$ by zero. A function $f$ on $M$ is called a basic $\eCrc$-function
on $M$ with support in a chart $U$ of $M$
if $f|_{U}\in\eCrc(U)$ and $f=(f|_{U})^{M}$. In this case
we write
$$ \supp_{U}(f)=\supp(f|_{U})\;.$$
A $\eCrc$-function on $M$ is a function on $M$ which is a sum of
basic $\eCrc$-functions on $M$, and the vector space of such
functions will be denoted by $\eCrc(M)$
(this coincide with (3) if $M$ is Hausdorff). 
Using partitions of unity we can see that if
$\cU$ is a cover of charts of $M$, we have
$$ \eCrc(M)=\sum_{U\in\cU}\eCrc(U)^{M}\;,$$
i.e. any $\eCrc$-function on $M$ is a sum of basic $\eCrc$-functions
with supports in the elements of $\cU$.

(5) Let $u:N\ra M$ be a local $\eCr$-diffeomorphism of
finite-dimensional $\eCr$-manifolds. There is an induced linear map
$u_{\iof}:\eCrc(N)\lra\eCrc(M)$, which is given by
$$ u_{\iof}(f)(x)=\!\!\!\!\!\sum_{y\in u^{-1}(x)}\!\!\!\!\!f(y)\;.$$
We say that an open subset $U\subset N$ is
{\em elementary} for $u$ if $u|_{U}$ is injective.
Observe that if $U$ is a chart of $N$,
$i:U\ra N$ the open embedding and $f\in\eCrc(U)$, we have
$$ i_{\iof}(f)=f^{N}\;.$$

\begin{dfn}  \label{Vdfn18}
Let $\GG$, $\HH$ and $\KK$ be object-separated
finite-dimensional \'{e}tale $\eCr$-groupoids. Let
$(E,p,\xW)$ be a $\eCr$-principal $\GG$-$\HH$-bibundle, and let
$(E',p',\xW')$ be a $\eCr$-principal $\HH$-$\KK$-bibundle.
Define a bilinear map
$$ \wp=\wp_{E,E'}:\eCrc(E)\times\eCrc(E')\lra\eCrc(E\ten{}E') $$
by
$$ \wp (m,m')(e\ten{}e')=\!\!\!\!\!\!\sum_{\cod h=\xW(e)}
   \!\!\!\!\!\! m(e\cdot h) m'(h^{-1}\!\!\cdot e')\;.$$
\end{dfn}
\Rem
Since $m'$ is a $\eCrc$-function and the fibers of $\xW'$ are discrete,
the value $m'(h^{-1}\!\!\cdot e')$ is zero for all but finitely many
$h$ with $\cod h=\xW(e)$, so the sum is always finite. It is clearly
invariant on the choice of a representable for the class $e\otimes e'$ in
$E\ten{}E'$. Thus $\wp (m,m')$ is a well-defined function on
$E\ten{}E'$. We have to check that it is a $\eCrc$-function.

To see this, we can assume that $m$ and $m'$ are basic,
with supports in the charts $U$ and $U'$ respectively.
We can further assume that $U$ is elementary for $\xW$ and $U'$
is elementary for $\xW'$.

The open subset $U\times_{\HH_{0}}U'$ of $E\times_{\HH_{0}}E'$
is a chart, and elementary for the quotient projection
$q:E\times_{\HH_{0}}E'\ra E\ten{}E'$ (which is a local
$\eCr$-diffeomorphism). In particular,
$q(U\times_{\HH_{0}}U')=U\ten{}U'$ is a chart.
It is $\eCr$-diffeomorphic to $V=U'\cap p'^{-1}(\xW(U))$,
and the diffeomorphism $\delta:V\ra U\ten{}U'$
is given by
$$ \delta(e')=(\xW|_{U})^{-1}(p'(e'))\ten{}e'\;.$$
We will show that $\wp (m,m')$ is a basic $\eCrc$-function 
with support in $U\ten{}U'$. 

First, it is clear that $\wp (m,m')$ is zero outside
$U\ten{}U'$. Next, let $e'\in V$. Then we have
$$ \wp (m,m')(\delta(e'))=m((\xW|_{U})^{-1}(p'(e'))) m'(e')=
   \xW_{\iof}(m)(p'(e')) m'(e')\;.$$
Note that the right hand side, which we will denote by $f'(e')$,
is well-defined for any $e'\in U'$, and $f'\in\eCrc(U')$.
Put $S=\supp \xW_{\iof}(m)$. Now $(p'|_{U'})^{-1}(S)$ is
closed in $U'$, therefore
$$ W=(p'|_{U'})^{-1}(S)\cap\supp_{U'}(m') $$
is compact. But clearly $W\subset V$ and $f'$ is zero outside $W$.
This proves that
$$ f'|_{V}=\wp (m,m')\com \delta\in\eCrc(V)\;,$$
hence $\wp (m,m')\in\eCrc(E\ten{}E')$.

\begin{prop}  \label{Vprop19}
Assume that $\GG$, $\HH$, $\KK$ and $\LL$ are
object-separated finite-dimensional \'{e}tale $\eCr$-groupoids. Let
$(E,p,\xW)$ be a $\eCr$-principal $\GG$-$\HH$-bibundle, let
$(E',p',\xW')$ be a $\eCr$-principal $\HH$-$\KK$-bibundle, and let
$(E'',p'',\xW'')$ be a $\eCr$-principal $\KK$-$\LL$-bibundle.
Then for any $m\in\eCrc(E)$, $m'\in\eCrc(E')$
and $m''\in\eCrc(E'')$ we have
$$ \wp (\wp (m,m'),m'')=
   \wp (m,\wp (m',m''))\;.$$
\end{prop}
\Rem
Here we have identified $(E\ten{}E')\ten{}E''$ and 
$E\ten{}(E'\ten{}E'')$ in the natural way.
\vspace{4 mm}

\Proof
Take any $e\ten{}e'\ten{}e''\in(E\ten{}E')\ten{}E''\cong
E\ten{}(E'\ten{}E'')$. Then
\begin{eqnarray*}
\lefteqn{\wp (\wp (m,m'),m'')
         (e\ten{}e'\ten{}e'') }                     \hspace{25 mm} \\
& = &
\!\!\!\!\!\!\!\!\sum_{\cod k=\xW'(e')}\!\!\!\!\!\!
\wp (m,m')(e\ten{}e'\cdot k) m''(k^{-1}\!\!\cdot e'') \\
& = &
\!\!\!\!\!\!\!\!\sum_{\cod k=\xW'(e')}\,
                \sum_{\cod h=\xW(e)}\!\!\!\!\!
m(e\cdot h) m'(h^{-1}\!\!\cdot e'\cdot k) m''(k^{-1}\!\!\cdot e'') \\
& = &
\!\!\!\!\!\!\!\sum_{\cod h=\xW(e)}\!\!\!\!\! m(e\cdot h)
\wp (m',m'')(h^{-1}\!\!\cdot e'\ten{}e'')                  \\
& =&
\wp (m,\wp (m',m''))(e\ten{}e'\ten{}e'')\;.
\end{eqnarray*}
\eop
\Rem
The algebra $\eCrc(\GG)$ of $\eCrc$-functions on an
object-separated finite-dimensional \'{e}tale $\eCr$-groupoid
$\GG$ \cite{Con82} is defined as the vector space $\eCrc(\GG_{1})$,
equipped with the product
$$ (x x')(g)=\!\!\!\sum_{g=g'\com g''}\!\!\! x(g') x'(g'')\;,$$
for any $x,x'\in\eCrc(\GG_{1})$ and $g\in\GG_{1}$.
The algebra $\eCrc(\GG)$ is clearly associative with local units.
But note that
$$ (x x')(g)=\wp_{\GG_{1},\GG_{1}}(x,x')(1\ten{}g)=
             \wp_{\GG_{1},\GG_{1}}(x,x')(g\ten{}1)\;.$$
In an analogous way we can define a structure of a
$\eCrc(\GG)$-$\eCrc(\HH)$-bimodule on $\eCrc(E)$, for any
$\eCr$-principal $\GG$-$\HH$-bibundle $E$, where $\GG$ and $\HH$ are
object-separated finite-dimensional \'{e}tale $\eCr$-groupoids, as
follows:

\begin{dfn}  \label{Vdfn20}
Let $\GG$ and $\HH$ be object-separated finite-dimensional \'{e}tale 
$\eCr$-groupoids, and let $(E,p,f)$ be a $\eCr$-principal 
$\GG$-$\HH$-bibundle.
Define a $\eCrc(\GG)$-$\eCrc(\HH)$-bimodule
structure on $\eCrc(E)$ by
$$ (x\cdot m)(e)=\wp (x,m)(1\ten{}e)=
   \!\!\!\!\!\sum_{\cod g=p(e)}\!\!\!\!\!
   x(g) m(g^{-1}\!\!\cdot e) $$
and
$$ (m\cdot y)(e)=\wp (m,y)(e\ten{}1)=
   \!\!\!\!\!\sum_{\cod h=\xW(e)}\!\!\!\!\!
   m(e\cdot h) y(h^{-1}) $$
for any $x\in\eCrc(\GG)$, $m\in\eCrc(E)$, $y\in\eCrc(\HH)$ and
$e\in E$.
\end{dfn}
\Rem
Proposition \ref{Vprop19} implies that this gives indeed a
bimodule structure on $\eCrc(E)$.
For example,
\begin{eqnarray*}
      ((a\cdot m)\cdot b)(e)
& = & \wp (\wp (a,m),b)
      (1\ten{}e\ten{}1)                                     \\
& = & \wp (a,\wp (m,b))
      (1\ten{}e\ten{}1)                                     \\
& = & (a\cdot(m\cdot b))(e)\;.
\end{eqnarray*}
It is obvious that $\eCrc(E)$ is locally unital.
\vspace{4 mm}

Let $\GG$, $\HH$ and $\KK$ be object-separated
finite-dimensional \'{e}tale $\eCr$-gro\-up\-oids, let
$E$ be a $\eCr$-principal $\GG$-$\HH$-bibundle, and let
$E'$ be a $\eCr$-principal $\HH$-$\KK$-bibundle.
Proposition \ref{Vprop19} implies
$$ \wp (\wp (m,y),m')(e\ten{}1\ten{}e') =
   \wp (m,\wp (y,m'))(e\ten{}1\ten{}e') $$
for any $m\in\eCrc(E)$, $y\in\eCrc(\HH)$, $m'\in\eCrc(E')$ and
$(e,e')\in E\times_{\HH_{0}}E'$. In other words, after
identifications
$E\ten{}\HH_{1}\cong E$ and $\HH_{1}\ten{}E'\cong E'$, we have
$$ \wp (m\cdot y,m') = \wp (m,y\cdot m')\;.$$
Therefore $\wp $ induces a linear map
$$ \mho=\mho_{E,E'}:\eCrc(E)\ten{\eCrc(\HH)}\eCrc(E')
   \lra\eCrc(E\ten{}E')\;.$$
If the groupoids are separated, this map turns out to be an
isomorphism:

\begin{theo}  \label{Vtheo21}
Let $\GG$, $\HH$ and $\KK$ be separated finite-dimensional
\'{e}tale $\eCr$-groupoids, let $(E,p,\xW)$ be a
$\eCr$-principal $\GG$-$\HH$-bibundle and let
$(E',p',\xW')$ be a $\eCr$-principal $\HH$-$\KK$-bibundle.
Then the map
$\mho :\eCrc(E)\ten{\eCrc(\HH)}\eCrc(E')
\ra\eCrc(E\ten{}E')$ given by
$$ \mho (m\ten{}m')(e\ten{}e')=
   \!\!\!\!\!\!\sum_{\cod h=\xW(e)}\!\!\!\!\!\! 
  m(e\cdot h) m'(h^{-1}\!\!\cdot e') $$
is an isomorphism of $\eCrc(\GG)$-$\eCrc(\KK)$-bimodules.
\end{theo}

Before we prove the theorem, we have to study the elements
of the tensor product
$\eCrc(E)\ten{\eCrc(\HH)}\eCrc(E')$ more closely.
\vspace{4 mm}

\Not
Let $\GG$, $\HH$ and $\KK$ be separated finite-dimensional
\'{e}tale $\eCr$-groupoids, let $(E,p,\xW)$ be a
$\eCr$-principal $\GG$-$\HH$-bibundle and let
$(E',p',\xW')$ be a $\eCr$-principal $\HH$-$\KK$-bibundle.
If $m\in\eCrc(E)$ and $m'\in\eCrc(E')$, we denote by
$m\ten{}m'$ the tensor in $\eCrc(E)\ten{\eCrc(\HH)}\eCrc(E')$.
Let $U$ be a chart of $E$ elementary for
$\xW$ and $U'$ a chart of $E'$ elementary for
$\xW'$. A {\em basic pair} (of $\eCrc$-functions)
on $(U,U')$ is a pair
$(m,m')$, with $m\in\eCrc(E)$ and $m'\in\eCrc(E')$, such that
$m$ has support in $U$, $m'$ has support in $U'$ and
$p'(U')\subset \xW(U)$. A basic pair $(m,m')$ on $(U,U')$
is called normalized if
$$ \xW_{\iof}(m)|_{p'(\supp m')}\equiv 1\;.$$

\begin{lem}  \label{Vlem22}
Let $\GG$, $\HH$ and $\KK$ be separated finite-dimensional
\'{e}tale $\eCr$-groupoids, let $(E,p,\xW)$ be a
$\eCr$-principal $\GG$-$\HH$-bibundle and let
$(E',p',\xW')$ be a $\eCr$-principal $\HH$-$\KK$-bibundle.
Let $U$ and $U'$ be charts of $E$ respectively $E'$,
elementary for $\xW$ respectively $\xW'$.
\begin{enumerate}
\item [(i)]   Let $m$ and $m'$ be $\eCrc$-functions on
              $E$ respectively $E'$ with supports in $U$ respectively
              $U'$. Then there exists an open subset $V\subset U'$ and
              a $\eCrc$-function $m'_{1}$ on $E'$ with support in $V$
              such that $(m,m'_{1})$ is a basic pair on $(U,V)$ and
              $$ m\ten{}m'=m\ten{}m'_{1}\;.$$
\item [(ii)]  For every basic pair $(m,m')$ on $(U,U')$ there
              exists a normalized basic pair $(m_{1},m'_{1})$
              on $(U,U')$ such that 
              $$ m\ten{}m'=m_{1}\ten{}m'_{1}\;.$$
\item [(iii)] If $(m,m')$ and $(m_{1},m')$ are normalized basic
              pairs on $(U,U')$, then 
              $$ m\ten{}m'=m_{1}\ten{}m'\;.$$
\end{enumerate}
\end{lem}
\Proof
(i) Take $V=U'\cap p'^{-1}(\xW(U))$, and choose
$y\in\eCrc(\HH)$ with $\supp y\subset \xW(U)$ and
$y|_{\supp(\xW_{\iof}(m))}\equiv 1$.
First observe that $m\cdot y=m$, hence
$$ m\ten{}m'=m\cdot y\ten{}m'=m\ten{}y\cdot m'\;.$$
For any $e'\in E'$ we have
$$ (y\cdot m')(e')=y(p'(e')) m'(e')\;.$$
As in the remark to Definition \ref{Vdfn18} we can conclude that
$m'_{1}=y\cdot m'$ has the support in $V$.

(ii) Choose $m_{1}\in\eCrc(E)$ with support in $U$ such that
$m_{1}|_{\supp m}\equiv 1$, and put $y=(\uni\com\xW)_{\iof}(m)$.
Then we have $m_{1}\cdot y=m$ and hence
$$ m\ten{}m'=m_{1}\cdot y\ten{}m'
   =m_{1}\ten{}y\cdot m'\;.$$
Observe that $(m_{1},y\cdot m')$ is a normalized basic pair
on $(U,U')$.

(iii) Let $y=(\uni\com\xW)_{\iof}(m)$ and
$y_{1}=(\uni\com\xW)_{\iof}(m_{1})$. Note that
$m\cdot y_{1}=m_{1}\cdot y$ and
$y\cdot m'=y_{1}\cdot m'=m'$. Hence
\begin{eqnarray*}
m\ten{}m' & = & m\ten{}y_{1}\cdot m' =
                m\cdot y_{1}\ten{}m'                \\
          & = & m_{1}\cdot y\ten{}m' =
                m_{1}\ten{}y\cdot m'                \\
          & = & m_{1}\ten{}m'\;.
\end{eqnarray*}
\eop

\begin{lem}  \label{Vlem23}
Let $\GG$, $\HH$ and $\KK$ be separated finite-dimensional
\'{e}tale $\eCr$-groupoids, let $(E,p,\xW)$ be a
$\eCr$-principal $\GG$-$\HH$-bibundle and let
$(E',p',\xW')$ be a $\eCr$-principal $\HH$-$\KK$-bibundle.
Let $U$ and $U'$ be
charts of $E$ respectively $E'$,
elementary for $\xW$ respectively $\xW'$.
Let $(m,m')$ be a basic pair on $(U,U')$ and let
$V$ be a chart of $E'$, elementary for $\xW'$,
such that $\xW'(U')\subset \xW'(V)$. Then
there exist charts $U_{i}$ of $E$ elementary for
$\xW$, open subsets $U'_{i}$ of $V$ and
normalized basic pairs $(m_{i},m'_{i})$ on $(U_{i},U'_{i})$,
$i=1,2,\ldots,n$, such that
$$ m\ten{}m'=\sum_{i=1}^{n} m_{i}\ten{}m'_{i}\;.$$
\end{lem}
\Proof
Since $E'$ is principal there is a continuous map $\theta:U'\ra\HH_{1}$
such that $\cod\theta(e')=p'(e')$ and $\theta(e')^{-1}\!\!\cdot e'\in V$
for any $e'\in U'$. Put $S=\supp m'\subset U'$.

First we assume that $\theta(S)\subset W$ for a chart $W$ of
$\HH_{1}$ elementary for $\dom$. Choose $y\in\eCrc(\HH)$ with support
in $W$ such that $y|_{\theta(S)}\equiv 1$.
Note that $\alpha=(\xW'|_{V})^{-1}\com \xW'|_{U'}:U'\ra V$ is an open
$\eCr$-embedding, and
$$ \alpha(e')=\theta(e')^{-1}\!\!\cdot e' $$
for any $e'\in U'$.
Define now $f\in\eCrc(E')$ with support in $\alpha(U')$ by
$f(e')=m'(\alpha^{-1}(e'))$, for any $e'\in\alpha(U')$.
It is now easy to see that
$$ y\cdot f=m'\;.$$
Now we can write $m\cdot y$ as a sum
of $\eCrc$-functions $m_{i}$ on $E$ with supports in $U_{i}$,
$i=1,2,\ldots,n$, where $U_{i}$ are charts elementary
for $\xW$. Hence we have
$$ m\ten{}m'=m\ten{}y\cdot f=m\cdot y\ten{}f=
   \sum_{i=1}^{n}m_{i}\ten{}f\;. $$
Finally we use Lemma \ref{Vlem22} (i) to find open subsets
$U'_{i}\subset\alpha(U')$ and $\eCrc$-functions $m'_{i}$ with supports
in $U'_{i}$ such that $(m_{i},m'_{i})$ is a basic pair on
$(U_{i},U'_{i})$ and $m_{i}\ten{}f=m_{i}\ten{}m'_{i}$,
for any $i=1,2,\ldots,n$. By Lemma \ref{Vlem22} (ii) we can normalize
those basic pairs. This completes the proof in the special case.

Now we will deal with the general case. Since $\theta(S)$ is a
compact subset of $\HH_{1}$, we can choose a finite cover
$(W_{j})_{j=1}^{k}$ of $\theta(S)$ consisting of charts of
$\HH_{1}$ elementary for $\dom$. Put $V_{j}=\theta^{-1}(W_{j})$. Now
$(V_{j})_{j=1}^{k}$ is an open cover of $U'$, and we can choose a
partition of unity $(u_{j})_{j=1}^{k}$ on $S\subset U'$ such that
$u_{j}:E'\ra [0,1]$ is a $\eCrc$-function with support in
$V_{j}\cap U'$.
In particular, $\sum_{j=1}^{k}u_{j}(e')=1$ for any $e'\in S$.
Put $u=\sum_{j=1}^{k}u_{j}$. Now we have
$$ m'=m'\cdot (\uni\com\xW')_{\iof}(u)=
   \sum_{j=1}^{k}m'\cdot (\uni\com\xW')_{\iof}(u_{j})\;,$$
hence
\begin{equation}  \label{Vformula1}
  m\ten{}m'=\sum_{j=1}^{k}m\ten{}m'\cdot
  (\uni\com\xW')_{\iof}(u_{j})\;.
\end{equation}
For each $j=1,2,\ldots,k$, the function
$m'\cdot (\uni\com\xW')_{\iof}(u_{j})$
has support in $V_{j}\cap U'$, and
$$ \theta(\supp (m'\cdot (\uni\com\xW')_{\iof}(u_{j})))
   \subset W_{j}\;.$$
We can now use the first part of the proof on each of the summands
of the equation (\ref{Vformula1}).
\eop

\noindent {\bf Proof of Theorem \ref{Vtheo21}.}
Proposition \ref{Vprop19} implies that $\mho$ is a homomorphism
of $\eCrc(\GG)$-$\eCrc(\HH)$-bimodules.
First we shall prove that the map $\mho$ is surjective.
Let $m''$ be a $\eCrc$-function on $E\ten{}E'$. We have to prove
that $m''$ is in the image of $\mho$. We can assume without
loss of generality that $m''$ has support in $U\ten{}U'$, where
$U$ and $U'$ are charts of $E$ respectively $E'$
elementary for $\xW$ respectively $\xW'$, and such that
$p'(U')\subset \xW(U)$. Thus we have the $\eCr$-diffeomorphism
$f:U'\lra U\ten{}U'$
given by $f(e')=(\xW|_{U})^{-1}(p'(e'))\ten{}e'$. Define
$m'\in\eCrc(E')$ with support in $U'$ by $m'(e')=m''(f(e'))$
for $e'\in U'$. Since
$$ S=(\xW|_{U})^{-1}(p'(\supp m')) $$
is a compact subset of $U$, we can choose $m\in\eCrc(E)$ with
support in $U$ such that $m|_{S}\equiv 1$. It is easy to see that
$$ \mho (m\ten{}m')=m''\;.$$

Finally, we have to prove that $\mho$ is injective.
Let $\eta\in\eCrc(E)\ten{\eCrc(\HH)}\eCrc(E')$ such that
$\mho (\eta)=0$. By Lemma \ref{Vlem22} we can write
$$ \eta=\sum_{i=1}^{n}m_{i}\ten{}m'_{i}\;,$$ 
where $(m_{i},m'_{i})$ is a normalized basic pair on
$(U_{i},U'_{i})$, for some charts $U_{i}$
and $U'_{i}$ for $E$ respectively $E'$ elementary for $\xW$
respectively $\xW'$. We have to prove that $\eta=0$.

(a) First assume that $n=1$. Then
$$ 0=\mho (\eta)(e\ten{}e')=m_{1}(e) m'_{1}(e') $$
for any $e'\in U'_{1}$ and $e\in E$ such that $\xW(e)=p'(e')$.
Now if $m'_{1}(e')\neq 0$ for some $e'\in U'_{1}$, we can
choose $e\in \xW^{-1}(p'(e'))$ such that
$m_{1}(e)=1$ since $(m_{1},m'_{1})$ is normalized. This is
a contradiction, hence we must have $m'_{1}=0$ and thus
$\eta=0$.

(b) Assume now that there exists a chart $V$ of $E$
elementary for $\xW$ such that $U'_{i}\subset V$ for any
$i=1,2,\ldots,n$. We will now prove that $\eta=0$ by induction on
$n$. The case $n=1$ is checked with (a). For the induction step,
denote $S'_{i}=\supp m'_{i}$, $S_{i}=(\xW|_{U_{i}})^{-1}(p'(S'_{i}))$,
$A_{i}=U_{n}\cap U_{i}$ and $C_{i}=S_{n}\cap S_{i}$,
for any $i=1,2,\ldots,n$. Put $A=\bigcup_{i=1}^{n-1}A_{i}$ and 
$C=\bigcup_{i=1}^{n-1}C_{i}$. First we will prove that
\begin{equation}  \label{Veq10}
  p'(S'_{n})\subset \xW(C)\subset \xW(A)\;.
\end{equation}
Since $E$ is separated, the set $C$ is compact and hence $\xW(C)$ is
closed in $\HH_{0}$. Thus it is enough
to see that $p'(e')\in \xW(C)$ for any $e'\in S'_{n}$ with
$m'_{n}(e')\neq 0$. Take such an $e'\in S'_{n}$, and let
$e=(\xW|_{U_{n}})^{-1}(p'(e'))$.
Observe first that $e\in S_{n}$. Now since $\mho (\eta)=0$
we have
$$ \mho (\eta)(e\ten{}e')=\sum_{i=1}^{n}m_{i}(e) m'_{i}(e')=0\;.$$
The assumption $m'_{n}(e')\neq 0$ and the fact that the pair
$(m_{n},m'_{n})$ is normalized imply $m_{n}(e)=1$. Thus there exists
$0\leq i\leq n-1$ such that
$$ m_{i}(e) m'_{i}(e')\neq 0\;.$$
Therefore $e\in A_{i}$, and moreover, since $m'_{i}(e')\neq 0$,
we should have $e\in S_{i}$. This yields $e\in C$, and proves the
equation (\ref{Veq10}).

Using Lemma \ref{Vlem22} (iii) we can thus
assume without loss of generality that
$\supp m_{n}\subset A\subset U_{n}$.
Choose now a partition of unity $(u_{i})_{i=1}^{n-1}$
on $p'(S'_{n})\subset \xW(A)$ such that
$u_{i}:\HH_{1}\ra [0,1]$ is a $\eCrc$-function with support in
$\xW(A_{i})$. In particular, $\sum_{i=1}^{n-1}u_{i}(b)=1$
for any $b\in p'(S'_{n})$. Put 
$u=\sum_{i=1}^{n-1}u_{i}$. We have
$$ m'_{n}=u\cdot m'_{n}=\sum_{i=1}^{n-1}u_{i}\cdot m'_{n}\;,$$
and hence
$$ m_{n}\ten{}m'_{n}=\sum_{i=1}^{n-1}m_{n}\cdot u_{i}\ten{}m'_{n}\;.$$
Note that $m_{n}\cdot u_{i}$ has support in $A_{i}$.
Using Lemma \ref{Vlem22} (i) and (ii) we can 
find open subsets $B_{i}\subset U'_{n}$ and
normalized basic pairs $(f_{i},f'_{i})$ on
$(A_{i},B_{i})$ such that
$$ f_{i}\ten{}f'_{i}=m_{n}\cdot u_{i}\ten{}m'_{n}\;,$$
for any $i=1,2,\ldots,n-1$. But since $A_{i}\subset U_{i}$, 
$(f_{i},f'_{i})$ is a normalized basic pair on
$(U_{i},B_{i}\cup U'_{i})$ as well.
By Lemma \ref{Vlem22} (iii) we can
assume without loss of generality that
$m_{i}|_{\supp f_{i}}\equiv 1$, by changing
$m_{i}$ if necessary. In particular,
$f_{i}\ten{}f'_{i}=m_{i}\ten{}f'_{i}$. Now we have
$$ \eta=\sum_{i=1}^{n-1}m_{i}\ten{}m'_{i} +
        \sum_{i=1}^{n-1}m_{i}\ten{}f'_{i}
       =\sum_{i=1}^{n-1}m_{i}\ten{}(m'_{i}+f'_{i})\;.$$
Observe that $(m_{i},m'_{i}+f'_{i})$ is a normalized basic pair
on $(U_{i},B_{i}\cup U'_{i})$, and $B_{i}\cup U'_{i}\subset V$.
The induction hypothesis thus gives $\eta=0$.

(c) Assume that there exists a chart $V$ of $E'$ elementary
for $\xW'$ such that $\xW'(U'_{i})\subset \xW'(V)$. Using Lemma
\ref{Vlem23} we can then reduce the proof to the case in (b).

(d) Finally we deal with the general case. Choose a partition
of unity $(v_{i})_{i=1}^{n}$ on 
$$ \bigcup_{i=1}^{n}\xW'(\supp m'_{i})\subset
   \bigcup_{i=1}^{n}\xW'(U'_{i}) $$
of $\eCrc$-functions $v_{i}:\KK_{1}\ra [0,1]$
with supports in $\xW'(U'_{i})$, $i=1,2,\ldots,n$. In
particular, $\sum_{i=1}^{n}v_{i}(c)=1$
for any $c\in\bigcup_{i=1}^{n}\xW'(\supp m'_{i})$.
Put $v=\sum_{i=1}^{n}v_{i}\in\eCrc(\KK)$.
Hence we have
\begin{equation}  \label{Vformula2}
  \eta=\eta\cdot v=\sum_{j=1}^{n}\eta\cdot v_{j}=
       \sum_{j=1}^{n}\sum_{i=1}^{n}m_{i}\ten{}m'_{i}\cdot v_{j}
\end{equation}
and $\mho (\eta\cdot v_{j})=\mho (\eta)\cdot v_{j}=0$.
Let $0\leq j\leq n$. Now $m'_{i}\cdot v_{j}$ has support
in $U'_{ij}=U'_{i}\cap \xW'^{-1}(\xW'(U'_{j}))$, and
$\xW'(U_{ij})\subset \xW'(U'_{j})$, for any $i$. Therefore
the summand
$$ \eta\cdot v_{j}=\sum_{i=1}^{n}m_{i}\ten{}m'_{i}\cdot v_{j} $$
of the equation (\ref{Vformula2}) satisfies the condition of the
case (c), by taking $V=V_{j}=U'_{j}$. Thus $\eta\cdot v_{j}=0$
for any $0\leq j\leq n$, hence $\eta=0$.
\eop
 
\begin{theo}  \label{Vtheo24}
The map $\eCrc:\cGpdbr\!\!\lra\cAlg$ is a functor. In particular,
if $\GG$ and $\HH$ are $\eCr$-Morita equivalent separated
finite-dimensional \'{e}tale $\eCr$-groupoids, then the algebras
$\eCrc(\GG)$ and $\eCrc(\HH)$ are Morita equivalent.
\end{theo}
\Proof
The functoriality of $\eCrc$ follows from
Theorem \ref{Vtheo21}.
\eop

\chapter*{References \markboth{REFERENCES}{REFERENCES}}
\pagestyle{myheadings}
\addcontentsline{toc}{chapter}{References}
\startXchapterskip

\begin{list}{[\arabic{enumi}]}
{\usecounter{enumi}\settowidth\labelwidth{[99]}
\leftmargin\labelwidth\advance\leftmargin\labelsep}

\small

\bibitem{AtiBot} M. F. Atiyah, R. Bott, The moment map and equivariant
                 cohomology. {\em Topology} 23 (1984), 1--28.

\bibitem{Bot}    R. Bott, Characteristic classes and foliations. {\em Springer 
                 Lecture Notes in Math.} 279 (1972), 1--94.

\bibitem{BryNis} J.-L. Brylinski, V. Nistor, Cyclic cohomology of \'{e}tale
                 groupoids. {\em K-Theory} 8 (1994), 341--365.

\bibitem{Bro}    R. Brown, From groups to groupoids: a brief survey.
                 {\em Bull. London Math. Soc.} 19 (1987), 113--134.

\bibitem{Bun}    M. Bunge, An application of descent to a classification
                 theorem. {\em Math. Proc. Cambridge Phil. Soc.} 107 (1990),
                 59--79.

\bibitem{Cam}    C. Camacho, A. Neto, {\em Geometric Theory of Foliations}.
                 Birkh\"{a}user, Boston (1985).

\bibitem{Con78}  A. Connes, The von Neumann algebra of a foliation.
                 {\em Springer Lecture Notes in Physics} 80 (1978),
                 145--151.

\bibitem{Con82}  A. Connes, A survey of foliations and operator algebras,
                 Operator algebras and applications. {\em Proc. Symposia
                 Pure Math.} 38 Part I (1982), 521--628.

\bibitem{Con85}  A. Connes, Non-commutative differential geometry. {\em Publ.
                 Math. I.H.E.S.} 62 (1985), 41--144.

\bibitem{Con86}  A. Connes, Cyclic Cohomology and the Transverse Fundamental
                 Class of a Foliation. {\em Pitman Research Notes in
                 Math.} 123, (1986), 52--144.

\bibitem{Ehr}    C. Ehresmann, Sur les espaces localement homog\`{e}nes.
                 {\em Enseignement Math.} 35 (1936), 317--333.

\bibitem{EhrRee} C. Ehresmann and G. Reeb, Sur les champs d'\'{e}l\'{e}ments
                 de contact de dimension $p$ compl\`{e}tement int\'{e}grables.
                 {\em C. R. Acad. Sc. Paris} 218 (1944), 995--997.

\bibitem{Est}    W. T. van Est, Rapport sur les $S$-atlas. {\em
                 Ast\'{e}risque} 116 (1984), 235--292.

\bibitem{Est1}   W. T. van Est, Alg\`{e}bres de Maurer-Cartan et holonomie.
                 {\em Ann. Facult\'{e} Sci. Toulouse} 97 (1989), 93--134.

\bibitem{FacSka} T. Fack and G. Skandalis, Some properties of the
                 $C^{\ast}$-algebra associated with a foliation.
                 {\em Proc. Symposia Pure Math.} 38 Part I (1982), 629--635. 
 
\bibitem{Hae}    A. Haefliger, Structures feuillet\'{e}es et cohomologie
                 \`{a} valeur dans un faisceau de groupo\"{\i}des. 
                 {\em Comment. Math. Helv.} 32 (1958), 248--329.

\bibitem{Hae84}  A. Haefliger, Groupo\"{\i}des d'holonomie et classifiants.
                 {\em Ast\'{e}risque} 116 (1984), 70--97.

\bibitem{HilS}   M. Hilsum and G. Skandalis, Stabilit\'{e} des alg\`{e}bres de
                 feuilletages. {\em Ann. Inst. Fourier Grenoble} 33 (1983),
                 201--208.

\bibitem{HilSka} M. Hilsum and G. Skandalis, Morphismes K-orientes 
                 d'espaces de feuilles et functorialite en theorie de Kasparov.
                 {\em Ann. Scient. Ec. Norm. Sup.} 20 (1987), 325--390.

\bibitem{Kui}    N. H. Kuiper, Sur les surfaces localement affines. {\em
                 G\'{e}om\'{e}trie Diff\'{e}rentielle, Colloques Internationaux
                 du Centre National de la Recherche Scietifique, Strasbourg}
                 (1953), 79--87.

\bibitem{Lan}    S. Lang, {\em Differential Manifolds}. Addison-Wesley, Reading
                 (1972).

\bibitem{Lod}    J.-L. Loday, {\em Cyclic Homology}. Springer-Verlag, Berlin
                 (1992).

\bibitem{Moe3}   I. Moerdijk, The classifying topos of a continuous groupoid,
                 I. {\em Transactions A.M.S.} 310 (1988), 629--668.

\bibitem{Moe4}   I. Moerdijk, Toposes and groupoids. {\em Springer Lecture
                 Notes in Math.} 1348 (1988), 280--298.

\bibitem{Moe1}   I. Moerdijk, Classifying toposes and foliations.
                 {\em Ann. Inst. Fourier, Grenoble} 41, 1 (1991), 189--209.

\bibitem{Mol}    P. Molino, {\em Riemannian Foliations}. Birkh\"{a}user,
                 Boston (1988).

\bibitem{Mrc4}   J. Mr\v{c}un, On stability of foliations invariant
                 under a group action. {\em Utrecht University,
                 Preprint} 900 (1995).

\bibitem{Mrc5}   J. Mr\v{c}un, Lipschitz spectrum preserving mappings
                 on algebras of matrices. {\em Linear Algebra Appl.}
                 215 (1995), 113--120.

\bibitem{Mrc3}   J. Mr\v{c}un, An extension of the Reeb stability
                 theorem. {\em Topology Appl.} 70 (1996), 25--55.

\bibitem{Pra}    J. Pradines, Morphisms between spaces of leaves viewed as
                 fractions. {\em Cahiers Top. G\'{e}om. Diff. Cat.} XXX-3
                 (1989), 229--246.

\bibitem{Qui88}  D. Quillen, Algebra cochains and cyclic cohomology.
                 {\em Publ. Math. I.H.E.S.} 68 (1988), 139--174.

\bibitem{Reeb}   G. Reeb, Sur certaines propri\'{e}t\'{e}s topologiques des
                 vari\'{e}t\'{e}s feuillet\'{e}es. {\em Actual. Sci. Ind.}
                 1183, Hermann, Paris (1952).

\bibitem{ReeSch} G. Reeb and P. Schweitzer, Un theoreme de Thurston etabli au
                 moyen de l'analyse non standard. {\em Springer Lecture Notes
                 in Math.} 652 (1978), 138--138.

\bibitem{Ren}    J. Renault, A groupoid approach to $C^{\ast}$-algebras.
                 {\em Springer Lecture Notes in Math.} 793 (1980).

\bibitem{RenMW}  J. Renault, P. S. Muhly and D. Williams, Equivalence and
                 isomorphism for groupoid $C^{\ast}$-algebras. {\em J.
                 Operator Theory} 17 (1987), 3--22.

\bibitem{Scha}   W. Schachermayer, Une modification standard de la
                 demonstration non standard de Reeb at Schweitzer. {\em 
                 Springer Lecture Notes in Math.} 652 (1978), 139--140.

\bibitem{Seg}    G. Segal, Classifying spaces and spectral sequences,
                 {\em Publ. Math. I.H.E.S.} 34 (1968), 105--112.

\bibitem{Spa}    E. H. Spanier, {\em Algebraic Topology}. McGraw-Hill,
                 New York (1966).

\bibitem{Thu}    W. P. Thurston, A generalization of the Reeb stability
                 theorem. {\em Topology} 13 (1974), 347--352.

\bibitem{Win}    H. Winkelnkemper, The graph of a foliation. {\em Ann.
                 Global Anal. Geom.} 1 (1983), 51--75.

\normalsize

\end{list}

\chapter*{Samenvatting \markboth{SAMENVATTING}{SAMENVATTING}}
\addcontentsline{toc}{chapter}{Samenvatting}
\startXchapterskip

\small

In dit proefschrift bestuderen wij Hilsum-Skandalis-afbeeldingen,
in het bijzonder Reeb-stabi\-li\-teit en algebra\"{\i}sche invarianten
van deze afbeeldingen.

Het eerste hoofdstuk is een overzicht van enkele basisbegrippen
betreffende foliaties, topologische groepo\"{\i}den en
Haefliger-structuren,
die de achtergrond voor ons werk vormen. Het belangrijkste
voorbeeld van een Haefliger-structuur is een foliatie op een vari\"{e}teit,
maar elke Haefliger-structuur op een topologische ruimte $X$ met waarden
in een \'{e}tale groepo\"{\i}de kan worden gezien als een gegeneraliseerde
foliatie op $X$: het induceert een verdeling van $X$ in bladen
voorzien van een fijnere topologie (bladtopologie) en met de
holonomiegroep. De holonomiegroepo\"{\i}de, of zijn Morita-equivalentie
klasse, is een goede representant voor de transversale structuur van een
foliatie, in tegenstelling tot de bladenruimte die te weinig informatie
bevat. Een foliatie $\cF$
op een vari\"{e}teit $M$ kan worden voorgesteld door een Haefliger-structuur
op $M$ met waarden in de holonomiegroepo\"{\i}de van $\cF$ gereduceerd tot
een complete transversaal. Deze structuur neemt de rol over van de
quoti\"{e}nt-projectie van $M$ op de transversale structuur van $\cF$,
en is een voorbeeld
van een Hilsum-Skandalis-afbeelding tussen topologische groepo\"{\i}den.
De Hilsum-Skandalis-afbeeldingen vormen een cate\-gorie waarin twee topologische
groepo\"{\i}den isomorf zijn precies indien zij Morita-equivalent zijn.

In het tweede hoofdstuk laten wij zien dat elke
Hilsum-Skandalis-afbeelding
tussen topologische groepo\"{\i}den $\HH$ en $\GG$ kan
worden gezien als een gegeneraliseerde foliatie op $\HH$: het induceert
een verdeling van de ruimte $\HH_{0}$ van de objecten van $\HH$ aan
$\HH$-invariante deelverzamelingen -- bladen -- toegerust met
de bladentopologie, holonomiegroep en een werking van $\HH$. Indien $\GG$
\'{e}tale is, is de holonomiegroep van een blad $L$ het beeld van een
(holonomie-) homomorfisme gedefinieerd op de fundamentaalgroep
$\pi_{1}(\HH(L))$ van
de groepo\"{\i}de $\HH(L)$ verbonden met de werking van $\HH$ op $L$.

De Reeb-stabiliteitsstelling zegt dat een compact blad van een foliatie
met eindige holonomiegroep een omgeving van compacte bladen heeft. In het
derde hoofdstuk breiden wij deze stelling uit, en ook zijn generalisaties
van Thurston en Haefliger, tot Hilsum-Skandalis-afbeeldingen. Om de
Reeb-Thurston-stelling uit te breiden introduceren wij het
lineaire holonomie-homo\-mo\-rfisme $d\cH_{L}$ van een blad $L$ van een
Hilsum-Skandalis-afbeelding $E$, van een topologische groepo\"{\i}de $\HH$,
naar een \'{e}tale $\eCe$-groepo\"{\i}de $\GG$ gemodelleerd op een
Banachruimte. Dit homomorfisme is gedefinieerd als een representatie
van de fundamentaalgroep $\pi_{1}(\HH(L))$ aan de raakruimte van de
vari\"{e}teit $\GG_{0}$ in een geschikte punt. Indien $E$ voldoet aan geschikte
voorwaarden, bewijzen wij: Laat $L$ een blad van $E$ zijn met compacte
banenruimte, zodat
\begin{enumerate}
\item [(i)]   $\Ker d\cH_{L}$ een eindig voortgebrachte groep is,
\item [(ii)]  $\Hom(\Ker d\cH_{L},\RRR)=\{0\}$, en
\item [(iii)] het beeld van $d\cH_{L}$ eindig is.
\end{enumerate}
Dan is de holonomiegroep van $L$ eindig, en $L$ heeft een $\HH$-invariante
omgeving van bladen van $E$ met compacte banenruimten.

In het vierde hoofdstuk illustreren wij de resultaten van het derde hoofdstuk
voor het geval van een equivariante foliatie op een vari\"{e}teit $M$ met een
werking van een discrete groep $G$. Aan een zodanige foliatie verbinden wij een
Hilsum-Skandalis-afbeelding van $G(M)$ naar de groepo\"{\i}de van kiemen van
diffeomorfismen van $\RRR^{q}$. Op deze manier verkrijgen wij de equivariante
versies van de Reeb- en de Reeb-Thurston-stabiliteitsstelling.

In het laatste, vijfde hoofdstuk presenteren wij een aantal
algebra\"{\i}sche
invarianten van Hilsum-Skandalis-afbeeldingen. Eerst geven wij een definitie
van de hogere homotopiegroepen van topologische groepo\"{\i}den met behulp
van de Hilsum-Skandalis-afbeeldingen, en bewijzen wij dat
die invariant zijn onder de
Morita-equivalentie. Onder bepaalde voorwaarden op een topolo\-gi\-sche
groepo\"{\i}de $\GG$ zijn deze groepen precies de hogere homotopiegroepen
van de classificerende ruimte en de classificerende topos van $\GG$.
Dan laten wij zien dat de singuliere homologie van topologische ruimten
met behulp van Hilsum-Skandalis-afbeeldingen
op een natuurlijke manier uitgebreid kan worden tot
de topologische groepo\"{\i}den.
Deze homologie is invariant onder de Morita-equivalentie en is veel eenvoudiger
dan de singuliere homologie van de classificerende ruimten van de 
groepo\"{\i}den. Voor een \'{e}tale groepo\"{\i}de $\GG$ 
drukken wij de singuliere homologie van $\GG$ uit als de homologie van
een quoti\"{e}nt van de singuliere complex van de ruimte $\GG_{0}$, en wij
laten zien dat de ``effect-functor'' een isomorfisme induceert tussen de
singuliere homologie van $\GG$ en die van $\Eff(\GG)$. Tot besluit
construeren wij voor elke Hilsum-Skandalis-afbeelding $E$ tussen
Hausdorff \'{e}tale eindig-dimensionale $\eCr$-groepo\"{\i}den
$\HH$ en $\GG$ een $\eCrc(\GG)$-$\eCrc(\HH)$-bimodule $\eCrc(E)$, waarin
$\eCrc(\GG)$ de Connes-algebra van $\eCrc$-functies op $\GG$ is. Op deze
manier krijgen wij een functor van een deelcategorie van
Hilsum-Skandalis-$\eCr$-afbeeldingen naar de categorie van bimodulen
over de algebra's met locale
eenheden. E\'{e}n van de gevolgtrekkingen is dat de
$\eCr$-Morita-equivalente Hausdorff \'{e}tale eindig-dimensionale
$\eCr$-groepo\"{\i}den Morita-equivalente algebra's hebben.

\normalsize

\chapter*{Curriculum Vitae \markboth{CURRICULUM VITAE}{CURRICULUM VITAE}}
\addcontentsline{toc}{chapter}{Curriculum Vitae}
\startXchapterskip

\small

The author of this thesis was born on April 13, 1966 in Ljubljana,
Slovenia. In 1986 he began the study of theoretical mathematics
at the University of Ljubljana, and he got his master degree in 1990.
For his master thesis he received the Pre\v{s}eren award of the
Faculty of Natural Science and Technology, University of Ljubljana.
As a ``Young Researcher'' of the Mathematical Institute, University of
Ljubljana, he completed his post-graduate ``magister'' study in 1993
at the Mathematical Department, University of Ljubljana.

In 1992 he started his research as an ``Assistent in Opleiding'' at the 
Mathematical Department, Utrecht University under the supervision of
Prof. dr. I. Moerdijk, which has led to this thesis.

\normalsize
  

\end{document}